\documentclass[12pt]{amsart}
\pdfoutput=1
\usepackage{hyperref}
\usepackage{latexsym, amsmath, amssymb,amsthm,amsopn,amsfonts}
\usepackage{version}
\usepackage{epsfig,graphics,color,graphicx,graphpap}
\usepackage{amssymb}
\usepackage{amsxtra}

\ifx\pdfoutput\undefined
  \DeclareGraphicsExtensions{.eps}
\else
  \ifx\pdfoutput\relax
    \DeclareGraphicsExtensions{.eps}
  \else
    \ifnum\pdfoutput>0
      \DeclareGraphicsExtensions{.pdf}

    \else
      \DeclareGraphicsExtensions{.eps}
    \fi
  \fi
\fi

\newcommand{\la}{\langle}
\newcommand{\ra}{\rangle}

\newcommand{\p}{\partial}
\newcommand{\Lop}{\mathcal L}
 
\newcommand{\CC}{{\mathbb C}}

\newcommand{\RR}{{\mathbb R}}

\renewcommand{\Re}{\mathop{\rm Re}\nolimits}

\newcommand{\gammaH}{\sqrt{\gamma'}H}
\newcommand{\Hgamma}{H_{\sqrt{\gamma'}} }
\newcommand{\rhoH}{{\rho H}}
\newcommand{\Hrho}{ H_\rho}

\theoremstyle{plain}

\newtheorem{thm}{Theorem}
\newtheorem{prop}{Proposition}[section]
\newtheorem{cor}[prop]{Corollary}
\newtheorem{lem}[prop]{Lemma}

\theoremstyle{definition}

\newtheorem{rem}{Remark}[section]
\newtheorem{defn}[prop]{Definition} 

\numberwithin{equation}{section}

\def\squarebox#1{\hbox to #1{\hfill\vbox to #1{\vfill}}} 
 
\newcommand{\sech}{\textnormal{sech}}

\title
[The Quartic KdV Equation]
{Small data scattering and soliton stability in $\dot{H}^{-\frac16}$ for the quartic KdV Equation}

\author[H. Koch]
{Herbert Koch}

\author[J.L. Marzuola]
{Jeremy L. Marzuola}

\address{Mathematics Institute, Bonn University \\
Endenicher Allee 60, D-53115 Bonn, Germany}
\email{koch@math.uni-bonn.de}

\address{Applied Mathematics Department, Columbia University \\
200 S. W. Mudd, 500 W. 120th St., New York City, NY 10027, USA}
\email{jm3058@columbia.edu}

\makeindex

\begin{document}

\begin{abstract}
  In this note we prove scattering for perturbations of solitons in the
  scaling space appropriate for the quartic nonlinearity, namely
  $\dot{H}^{-\frac16}$.  The article relies strongly on refined estimates for a KdV
  equation linearized at the soliton.  In contrast to  the work of Tao
  \cite{Tao} we are able to work purely in the scaling space without
  additional regularity assumptions, allowing us to construct wave operators 
and  a weak version of  inverse wave operators.
\end{abstract}

\maketitle

\section{Introduction and Statement of Results}
\label{sec:intro}

The generalized KdV equation
\begin{eqnarray}
\label{eqn:kdv}
\left\{ \begin{array}{c}
\p_t \psi + \p_x (\p_x^2 \psi + \psi^p) = 0, \ t,x \in \RR \\
\psi(0,x) = \psi_0 (x)
\end{array} \right.
\end{eqnarray}
has an explicit  soliton solution  
\[ \psi_c(x,t) = Q_{p,c,c^2t+x_0}(x):= c^{\frac{2}{p-1}} Q_p(c(x-(x_0+c^2 t))) \]
with $c>0$, $x_0 \in \mathbb{R} $ and 
\begin{equation} \label{Qps}  
 Q_p = \left( \frac{p+1}{2} \right)^{\frac{1}{p-1}} \sech^{\frac{2}{p-1}} \left( \frac{p-1}{2} x \right). 
\end{equation} 

Well-posedness of the generalized KdV equation was established by
Kenig-Ponce-Vega \cite{KPV} in $H^s$ for some $s$ depending on $p$. The
case $p=4$ (quartic KdV) is particularly interesting as it is the only
 subcritical
power nonlinearity that does not lead to a completely integrable
system.  The critical space for the quartic KdV equation is
$H^{-\frac16}$. Gr\"unrock \cite{Gr} obtained local wellposedness in $H^s$, $s>
-1/6$ and the endpoint $\dot H^{-\frac16}$ was reached by Tao \cite{Tao}.
Though wellposedness is not the main focus of this note, we will return to this
question in section \ref{sec:wellposedness} and use spaces of bounded $p$ variation
and their predual (see Appendix \ref{app:u2v2} and \cite{HHK}) to simplify and
strengthen Tao's wellposedness result in the critical space.

The solutions $Q_{c,y}$ are called traveling waves or solitons. These
are minimizers of the constrained variational problem 
\begin{equation}\label{constraintL2}  
\min \{ E(w) : w \in H^1, \ \Vert w \Vert_{L^2}  = \mu >0 \},
\end{equation}  
where 
\[ E(u) = \int \left[ \frac12  u_x^2 - \frac1{p+1} u^{p+1}  \right] dx. \] 
Minimizers also are extremals of the Lagrangian
\begin{equation}\label{globalmin}  
S(u) = E(u) + \frac{\lambda}2 \int u^2 dx,
\end{equation} 
where $\lambda$ is a Lagrangian multiplier.  
Existence of the minimizer has been shown by Berestycki and Lions
\cite{BeLi} using the constrained minimization problem
\begin{eqnarray*}
\min \{ T(w) : w \in H^1, \ V(w) = \tilde \mu \},
\end{eqnarray*}  
where
\begin{eqnarray*}
T(w) & = & \int w_x^2 dx , \\
V (w) & = & \frac{\lambda}2 \int w^2 dx - \frac1{p+1} \int w^{p+1}  dx . 
\end{eqnarray*}

The function $Q$ in \eqref{Qps} is the unique
positive even solution to the Euler-Lagrange equation
\begin{equation} \label{euler}  - Q_{xx} - Q^p +   Q  = 0 \end{equation}   
 to
\eqref{globalmin} with $\lambda=1$. It is a critical point of  $S(u)$ again with $\lambda=1$, a minimizer of $E$ with constraint $\Vert u \Vert_{L^2} = \mu$,  
\begin{equation} \label{mul2}  
\mu^2 = 
\Vert Q_p (x) \Vert^2_{L^2} = \left(\frac{p+1}2\right)^{\frac2{p-1}}
   \frac{ \Gamma(\frac{p+1}{p-1})\sqrt{\pi}}{\Gamma(\frac{p+3}{2(p-1)})}
 \end{equation}
and hence the quadratic form 
\begin{equation} \label{nonnegative} 
 K (\psi):= \int \frac12 {w'}^2 + \frac12 w^2 - \frac{p}2 Q^{p-1} w^2 dx \ge 0 \qquad
 \text{ for } \langle w,Q \rangle =0 \end{equation} 
is nonnegative on the tangent space i.e. the functions orthogonal to $Q$. 

The stability of solitons for generic KdV equations has been studied in several seminal works.  Orbital stability was first effectively established in the work of Weinstein \cite{W}.  Then asymptotic stability of solitons for KdV was first observed by Pego-Weinstein \cite{PW}, who proved that solitons for KdV are stable under perturbations in exponentially weighted spaces.  Later, Martel-Merle \cite{MM1,MM2,MM-Gafa} and Martel \cite{Martel} refined this result to observe that solitons for generalized KdV equations are indeed stable under perturbations in the energy space, but measured within a moving reference frame.  As mentioned above, for the case $p=4$, building on the multilinear estimates of Gr\"unrock \cite{Gr} and the work of Martel-Merle, Tao \cite{Tao} assumes smallness in 
 $H^1 \cap \dot H^{-\frac16}$ and obtains scattering in $\dot H^{-1/6}$. 
 We will give a more thorough introduction to previous stability results including rigorous definitions of stability in Section \ref{sec:stabrev}.

In the sequel we will focus on the case $p=4$ and omit $p$ in the
notation.  It seems that any further progress is tied to an understanding of the
linearization, or more precisely of the linear equation
\begin{equation}  \label{eqn:tdlinu} u_t  + \partial_x \mathcal{L} u = 0 
\end{equation} 
and its adjoint
\begin{equation} \label{eqn:tdlinv}
 v_t + \mathcal{L} \partial_x v = 0 ,
\end{equation} 
which have the explicit solutions (with $\tilde Q = c \partial_c Q_c|_{c=1} $)
\[ u= a ( \tilde Q +2t Q') + b Q', \   v= c Q, \]
 where 
\begin{equation}\label{tildeQ}   \tilde Q := \left. c \frac{d}{dc}
    c^{\frac{2}{p-1}} Q_p(cx)\right|_{c=1} = \frac{2}{p-1} Q_p + x Q',
\end{equation} 
usually evaluated at $c=1$.

Thus both
equations \eqref{eqn:tdlinu} and \eqref{eqn:tdlinv} have linearly
growing solutions.  It is one of the first contributions of this paper
that both equations are uniformly $L^2$ bounded once we take into
account these modes, and, moreover, there are local energy estimates
global in time once we remove these modes. In particular the
assumption of Pego and Weinstein on the absence of embedded
eigenvalues holds.

The goal of this work is to build on the arguments of Weinstein \cite{W} and
Martel-Merle \cite{MM1,MM2} to establish some type of asymptotic soliton 
stability for generalized KdV equations by a  direct analysis of the equation 
itself.  We apply a variant of Weinstein's and  Martel and Merle's arguments to the 
linear equations \eqref{eqn:tdlinv} and \eqref{eqn:tdlinv} and their relatives 
with variable scale and velocity, and control nonlinear terms through estimates for linear equations.  

Specifically,  we define projection
operators  related to the spectrum of
$\mathcal{L}$:
\begin{equation}\label{projection} 
P^\perp_{Q'} \psi = \psi - \frac{\langle \psi, Q' \rangle}{\langle Q',
  Q' \rangle } Q' , \quad \tilde{P} \psi = \psi - \frac{\langle \psi,
   Q \rangle }{\langle Q, \tilde Q \rangle} \tilde Q.
 \end{equation} 

We obtain the main linear estimates which in their simplest form can be written as  
\begin{thm}
\label{thm:linmain}
Let $S$ be the solution operator for \eqref{eqn:tdlinu} and $S^*$ the
solution operator for \eqref{eqn:tdlinv}.  Then, we have
\begin{equation}
\sup_t  \Vert  S(t)  \tilde{P}^* u_0 \Vert_{L^2}
+ \Vert \sech(x) \partial_x  P^\perp_{Q'} S(t)  \tilde{P}^* u_0  \Vert_{L^2(\mathbb{R}^2)}
  \lesssim  
\| u_0 \|_{L^2},
\end{equation} 
\begin{equation}  
\sup_t  \Vert S^*(t) P^\perp_{Q'}  v(t)  \Vert_{L^2}
+ \Vert \sech(x) \partial_x   \tilde P  S^* (t) P^\perp \Vert_{L^2(\mathbb{R}^2)}
  \lesssim  \| v_0 \|_{L^2}.
\end{equation} 
\end{thm}
We note that the linear estimates presented in the sequel may be
generalized to any subcritical power $p< 5$. We provide variants of
Theorem \ref{thm:linmain} for linearization at solitons with variable
scale and velocity as well as estimates in scales of Banach spaces
similar to estimates for the Airy equation.

Even near the trivial solution dominating the nonlinear part globally
by the linear parts requires to work in a scale invariant space
similar to $\dot H^{-\frac16}$. On the positive side it will lead to
scattering for perturbations of a soliton in $\dot H^{-\frac16}$,
without the smallness  condition of Tao in the energy space \eqref{energysmall}.  The study
of the linear equation will lead to a fairly precise understanding of its
properties  which seems to be new and we hope that it will provide a model
for many other questions on the stability of solitons.

As is standard in the study of stability, we take
\begin{eqnarray*}
\psi (x,t) = Q_{c(t)} (x - y (t)) + w(x,t).
\end{eqnarray*}
Then, we have
\begin{equation} 
\begin{split} 
\label{eqn:coupledsystem-w}
\p_t w + \p_x (\p_x^2 w + 4 Q_c^3 w)  = & - \dot{c} (\partial_c Q_{c}) (x-y) + \dot{y} ( Q_{c}') (x-y) \\
& - \p_x ( \p_x^2 Q_c - c^2 Q_c + Q_c^4) - c^2 (Q_c' ( x-y)) \\
& - \p_x (6 Q_c^2 (x-y) w^2 + 4 Q_c (x-y) w^3 + w^4).
\end{split} 
\end{equation}

 The standard choice of  $\dot c$ and $\dot y$ ensures 
 orthogonality conditions for $w$. Due to low time regularity 
we are forced  to relax the orthogonality conditions to
\begin{eqnarray}
\label{eqn:coupledsystem-c}
\frac{\dot{c}}{c} \langle Q_c, \tilde{Q}_c \rangle &=&  \langle w, Q_c \rangle, \\
\label{eqn:coupledsystem-x}
(\dot y -c^2) \langle Q_c',Q_c' \rangle & = & - \kappa \langle w , Q_c' \rangle,
\end{eqnarray} 
where $\kappa \gg 1$.

From an
implicit function theorem argument similar to that in the proof of Proposition
$1$ of \cite{MM-Gafa} there exist unique $c(0)$ and $y(0)$ so that $w(.,0)$
is orthogonal to $Q_{c(0)}(.-y(0))$ and $Q'_{c(0)} (.-y(0))$ provided the distance of $\psi$ to the set of solitons is small in a suitable norm.

 We consider the
equations above as ordinary differential equations for $c$ and $y$,
coupled with the partial differential equation.

Using the decomposition and linear estimates, in Sections \ref{sec:multlinest}
and \ref{sec:apps} we can prove (refering to later sections for the definition of the function spaces, with $\dot B^{-1/6,2}_\infty$ slightly larger than $\dot H^{-1/6}$) the following global result
\begin{thm}
\label{thm:main-result}
There exists  $\epsilon > 0$ and $c >0$ such that given \eqref{eqn:kdv} with initial data of the form
\begin{eqnarray*}
\min_{c_0,y_0} \| \psi_0 - Q_{c_0} (x - y_0) \|_{\dot B^{-\frac16,2}_{\infty}} \leq \epsilon,
\end{eqnarray*}
there exist unique functions $c$ and $y$ with 
\begin{eqnarray*} 
\langle w(0), Q_{c(0)} \rangle = \langle w(0), Q'_{c(0)} \rangle = 0 ,
\\
\dot c \in L^1 \cap C^0,\quad \dot y -c^2 \in L^2 \cap C^0, \\
\end{eqnarray*}
and a function $w(x,t) \in \dot X^{-\frac16}_\infty$ such that
\begin{eqnarray*}
\psi (x,t) = Q_{c(t),y(t)} (x) + w(x,t) 
\end{eqnarray*}
satisfies the quartic KdV equation, and $w$, $c$ and $y$ satisfy
\eqref{eqn:coupledsystem-c},\eqref{eqn:coupledsystem-x} and \eqref{eqn:coupledsystem-w}. 
Moreover, 
\[ 
\| \dot c \|_{L^1 \cap C^0} + \| \dot y - c^2 \|_{L^2 \cap C^0}
+ \Vert w \Vert_{\dot X^{-\frac16}_{\infty}} \le c  \Vert w_0  \Vert_{\dot
  B^{-\frac16,2}_\infty}.   
\]
In addition, there exists a function $z_0 \in \dot{B}^{-\frac16,2}_\infty$ such that
\begin{eqnarray*}
\| w(t) - e^{-t \p_x^3} z_0 \|_{\dot{B}^{-\frac16,2}_\infty} \to 0
\end{eqnarray*}
and 
\[ 
\|  w(.) - e^{-.\p_x^3} z_0 \|_{X^{-1/6}_\infty((t, \infty))} \to 0 \quad \text{ as } t \to \infty 
\]
if $w(0)$ is in the closure of $C^\infty_0$. 
\end{thm}
In fact, we prove a far stronger result than this, though Theorem
\ref{thm:main-result} captures the main ideas. 
 Finally, in Section \ref{sec:inverse} we show for a function $v$, there exists a quantity $J(v)$ defined in \eqref{J} such that we have the following
\begin{thm}  
\label{thm:inverse}
Let $v_0$ be in the closure of $C^\infty_0$ in $\dot B^{-\frac16,2}_{\infty}$, $c_\infty>0$, $y_0 \in \mathbb{R}$. Let $v$ be the solution to the linear homogeneous $KdV$ equation.  Assume that 
\[ J(v) \le \delta \] 
for some $\delta = \delta (\| v_0 \|_{\dot{B}^{-\frac16,2}_{\infty}})$. 
Then there exists a solution $\Psi$ to the quartic KdV equation,
 a function $y\in C^1([0,\infty))$, $c \in C^1([0,\infty), (0,\infty))$
such that $w = \Psi-Q_{c,y}$, $c$ and $y$ 
 satisfy equations
 \eqref{eqn:coupledsystem-c},\eqref{eqn:coupledsystem-x},
 \eqref{eqn:coupledsystem-w}, and
\[ \langle w(0), Q_{c(0)}(.-y(0))\rangle =  \langle w(0),
Q'_{c(0)}(.-y(0))\rangle  =0 , \]
\[ c(t) \to c_\infty,\quad  y(0) = y_0, \quad    w(t) - v(t) \to 0 \text{ in }
\dot{B}^{-\frac16,2}_\infty \text{ as } t \to \infty. \]
Moreover, if in addition $ v_0 \in L^2$, then $\Psi \in C(\mathbb{R}, L^2(\mathbb{R}))$ and
\[ \Vert v_0 \Vert_{L^2}^2 + \Vert Q_{c_\infty,0} \Vert_{L^2}^2 = \Vert \Psi(t) \Vert_{L^2}.  \]
There exists $\varepsilon>0$ such that the assumptions are satisfied if $\Vert v_0 \Vert_{\dot{B}^{-\frac16,2}_{\infty}} \le \varepsilon$.
\end{thm}

\begin{rem}
\label{rem:sobreg}
The conclusions in Theorems \ref{thm:main-result} and
\ref{thm:inverse} hold as well in the spaces
$\dot{B}^{-\frac16,2}_{\infty} \cap \dot H^s \cap H^\sigma$ for any
$-1 < s \le 0$ and $\sigma \geq 0$, allowing one to prove uniform
bounds in higher Sobolev norms, see Section \ref{variants}.  In
particular, given initial data in $\dot{B}^{-\frac16,2}_{\infty} \cap
\dot H^s \cap H^\sigma$, $J$ small will imply stability and scattering
in $\dot{B}^{-\frac16,2}_{\infty} \cap \dot H^s \cap H^\sigma$.
Specifically, we note one can prove boundedness and scattering in the
energy space $H^1$, intersected with $\dot{B}^{-\frac16,2}_\infty$.
\end{rem}

To motivate the construction of our nonlinear iteration spaces, in Section \ref{sec:linkdv} we first derive some refined estimates for the linear KdV equation
\begin{equation}
\label{eqn:linkdv}
\left\{ \begin{array}{c}
\p_t u + \p_x^3  u = f , \\
u(0,x) = u_0 (x) .
\end{array} \right.
\end{equation}
Then, in Section \ref{sec:lin-subcrit} we discuss the spectral and
mapping properties of the
operator $\mathcal{L}$
and derive linear estimates for the systems \eqref{eqn:tdlinu}
and  \eqref{eqn:tdlinv} and their relatives 
\[ u_t + u_{xxx} + (Q_{c(t)} (x-x(t) )u)_x = f . \] 
 In Section \ref{sec:virial}, we
combine local smoothing estimates as for \eqref{eqn:linkdv}, where we
treat the $Q$ terms as error terms with the virial identity and energy conservation for
\ref{eqn:tdlinu} to prove uniform bounds for a projection of the
solution $v$ assuming orthogonality of the initial data to $Q'$.

With this first result at hand we pursue a standard though nontrivial path 
and  employ pseudodifferential techniques and duality to derive similar 
estimates in a full scale of function spaces. The Littlewood-Paley 
decomposition at low frequencies is severely affected by the 
 term containing $Q$. This is done in Section \ref{sec:lin-constant}
with main result Proposition  \ref{littlewood}.

Theorems \ref{thm:main-result} and \ref{thm:inverse} are proven in the final two sections by
combining the well-posedness arguments and the linear estimates.

{\sc Acknowledgments} 
J.L.M. was funded by a Hausdorff Center Postdoc
at the University of Bonn and by a National Science Foundation
Postdoctoral Fellowship. H.K. was partially supported by the DFG
through Sonderforschungsbereich 611.  The authors wish to thank Axel
Gr\"unrock and Yvan Martel for helpful comments on an early version of
the result.

\section{Review of Previous Soliton Stability Results}
\label{sec:stabrev}

To begin, we consider the linearized operator
\[ \mathcal{L} \psi = - \psi'' - p Q^{p-1} \psi + \psi \] 
associated
to the Euler-Lagrange equation \eqref{euler} 
of \eqref{globalmin} with $\lambda=1$
respectively the constraint variational problem \eqref{constraintL2}
with Lagrange multiplier $1$.  It is one of the remarkable operators
for which almost everything is know about the spectrum and scattering,
see Lamb \cite{Lamb} Section 2.4 and 2.5, and Titchmarsh
\cite{Titchmarsh} subsection 4.19. The operator
\[ \mathcal{L}_M \psi = -\psi_{xx} - M sech^2 (x) \psi \] 
has the continuous
spectrum $[0,\infty)$ and the ground state $\psi_0(x)=\sech^\alpha (x)$ with
eigenvalue $\alpha^2$ provided $M= \alpha(\alpha+1)$, $\alpha>0$.  The
other eigenvalues are $(\alpha -j)^2$ for $1\le j < \alpha$ together with the 
eigenfucntions can be
obtained as follows: Let $\psi_{0,M}$ be the ground state with the
constant $M$. Then,
\[\psi_{j, (\alpha+j)(\alpha+j+1)}(x) = \prod_{l=1}^j ( \frac{d}{dx}
-(\alpha+l)\tanh (x)) \sech^\alpha(x) 
\]
is the $j$ eigenfunction to the potential with $M =(\alpha+j)(\alpha+j+1) $. 
We consider this information useful, and we will use  these results, 
even if the arguments could easily adapted to a much larger class of 
nonlinearities.

Clearly,
\[ \mathcal{L} Q'  = 0 \]
and  a short calculation or a comparison with the results above shows that
$Q^{\frac{p+1}{2}}$ is the ground state with eigenvalue
$1-\left(\frac{p+1}{2}\right)^2$. There is no other eigenvalue if $p \ge 3$, 
but there are other eigenvalues in $(0,1)$ if $p < 3$. 
 As an immediate consequence 
$ K (\psi) \ge \Vert \psi \Vert_{L^2}^2$
if  $ \langle \psi, Q' \rangle = \langle \psi, Q^{\frac{p+1}{2}} \rangle =0$.

We recall that $K$ is positive definite on the orthogonal complement of $Q$. 
We follow Weinstein \cite{W} and use this  bound to establish a lower bound on a different 
codimension $2$ subspace if $p < 5$. There exists $\delta>0$ such that 
\begin{equation} \label{contlowerbound}  K(\psi)  \ge \delta  \Vert \psi
  \Vert_{H^1}^2 \qquad \text{ for all } \psi  \text{ with } \langle \psi,
  Q^{p-1} Q' \rangle = \langle \psi  , Q \rangle =0 .
\end{equation} 
It suffices to verify this statement independently for odd and even
functions. For odd functions the quadratic form is nonnegative, with a
null space spanned by $Q'$. Positivitiy  follows from $\langle Q', Q^{p-1}
Q' \rangle \ne 0$. The argument for even functions is harder, but again
the quadratic form is nonnegative since $Q$ is a local minimizer of
the constraint variational problem.  

Let $\psi_j$ a minimizing sequence with $\Vert \psi_j \Vert_{H^1}=1$.
Suppose that the left hand side of \eqref{contlowerbound} converges to $0$.  The sequence
maximizes $\int Q^{p-1} \psi_j^2 dx$.  There exists a weakly converging
subsequence 
which convergences against a nontrivial even limit $\psi$ 
since $\psi \to \int Q^{p-1} \psi^2 dx >0$ is weakly lower semicontinuous.
Moreover  $\langle \psi, Q \rangle =0$ and $\Vert \psi \Vert_{H^1} \le 1$. Rescaling if necessary we  see that $\Vert \psi \Vert_{H^1} = 1$. 
 
We want to show that $K(\psi) >0$ and argue by contradiction. Suppose that 
$K(\psi)=0$. Then by \eqref{nonnegative} $\psi$ is a minimizer of $K$ under the
sole constraint $\langle Q, \psi \rangle =0$ and hence it satisfies the Euler Lagrange equations 
\[  \mathcal{L} \psi = \lambda Q . \]
But then $\psi$ is a multiple of $\tilde Q$ since 
\[ \mathcal{L} \tilde Q = - 2 Q \]
is the unique symmetric function with this property.  
However, $\psi$ is orthogonal to $Q$, but $\tilde Q$ and $Q$ are linearly independent  if $p \ne 5$, hence  $\psi=0$, which contradicts our construction
and thus  implies the existence of $\delta >0$ with
\[ K (\psi) \ge \delta \Vert \psi \Vert_{H^1}. \]
Observe that here the subcriticality condition $p <5$ enters crucially.

Given $\psi$ we define the parameters  
 $c_0$ and $x_0$  by the variational problem  
\[ \Vert \psi - Q_{c_0,x_0} \Vert_{H^1}^2 = \inf_{c,x } \Vert \psi -
Q_{c,x} \Vert_{H^1}^2. \] 
Following  Weinstein \cite{W} we claim 
\begin{equation} \label{lowerbound} 
 \Vert \psi - Q_{c_0,x_0} \Vert_{H^1}^2 \le c (E(\psi) - E(Q_{c})) 
\end{equation} 
provided the left  hand side is sufficiently small.
This  is a  consequence of   the  lower bound for the quadratic form
\eqref{contlowerbound}.

Lyapunov stability of solitons has been shown in the seminal work of
Weinstein \cite{W}, Theorem 4: Let $\varepsilon >0$. There exists $\delta >0$
such that, if
\[ \Vert \psi_0 - Q_1 \Vert_{H^1} \le \delta, \]
then 
\[ \inf_{x_0}   \Vert \psi(t)- Q_1(x-x_0) \Vert_{H^1} \le \varepsilon. \]
This is a direct consequence of 
the  conservation of the $L^2$ norm and the energy, plus \eqref{lowerbound}.

The study of asymptotic stability began with Pego-Weinstein \cite{PW}
in spaces with growing exponential weights. The effect
of the weight is twofold. First, there is not much the soliton could
interact with on its path to the right. Secondly, small solitons which
are slow and prevent asymptotic stability in $L^2$ carry a weight which
makes them exponentially decreasing in time. A key assumption is the absence 
of embedded eigenvalues of $\partial_x \mathcal{L}$, other than $0$ with 
eigenfunction  $Q'$ and the generalized eigenfunction $\tilde Q$. 
Pego and Weinstein verify this assumption for $p=2$ and $p=3$ and show that 
it fails at at most a finite number of values for $p$ between $2$ and $5$. 
It is a consequence of the virial identity below that there are  no nonzero purely imaginary  eigenvalues of $\partial_x \mathcal{L}$.

The exponential weight pushes the continuous spectrum of $\partial_x
\mathcal{L} $ to the left, makes the problem more parabolic,  and 
allows the use of techniques from smooth
dynamical systems, in particular of a center manifold reduction which is a restriction of the flow to a two dimensional manifold.

Martel and Merle (\cite{MM1,MM2}) and Martel \cite{Martel} introduced
a virial identity or monotonicity formula for the adjoint problem
\eqref{eqn:tdlinv} as well as for nonlinear problems. 
 Let 
 \begin{eqnarray*}
 \eta (x) = -\frac{p+1}{p-1} \frac{Q'}{Q}= \frac{p+1}2 \tanh \frac{p-1}2  x 
\end{eqnarray*}
and suppose that $v$ satisfies the equation
\eqref{eqn:tdlinv}.  By direct computation we have
\begin{equation} 
-\frac{d}{dt} \int \eta  v^2 dx  =   \langle (3 (\mathcal{L}+ \frac{(p+1)^2}{4}-1)  Q^{\frac{p-1}2} v, Q^{\frac{p-1}2} v \rangle ,
\end{equation}
where the quadratic form is nonnegative and it has by the spectral
theory of Schr\"odinger operators with $\sech^2(x)$ potentials a one
dimensional null space spanned by $Q$.  There are two consequences:
the quantity on the left hand side is monotonically decreasing, and
the right hand side controls the $H^1$ norm of $Q^{\frac{p-1}2}v$
provided $v$ is orthogonal to a vector $\bar Q$ with $\langle \bar Q,
Q \rangle \ne 0$.  Hence, if $v(0)$ is orthogonal to $Q'$ and $\tilde Q$ -
which is preserved under the evolution - 
\[ \Vert Q^{(p-1)/2} v \Vert_{H^1} \le c \sup_t \Vert v(t) \Vert_{L^2}. \]
The left hand side   is controlled provided we obtain a bound on
$\sup_t \Vert v(t) \Vert_{L^2}$.
  Martel and Merle (\cite{MM1,MM2}) use this and related
observations together with the a priori control on the deviation of the
solution to the set of solitons in ingenious ways for indirect arguments:
The existence of a solution $H^1$ close to solitons, but not
asymptotically converging to the soliton 'on the right' leads to the 
existence of impossible objects.

Later, C\^ote \cite{Cote} constructed solutions with specific
asymptotic conditions including many soliton solutions for positive
time. This shows that $L^2$ convergence to a soliton
will not be true without restricting the set were convergence is
studied.

Already $L^2$ conservation precludes asymptotic stability of the
trivial solution. The relevant notion instead of asymptotic stability is 
for  unitary problems the notion of  scattering.
Suppose that $\psi(0)$ is close to a soliton. We seek a function 
$w$ satisfying the Airy equation as well as $c(t)$ and $y(t)$ and a Banach space $X$ so that 
\[ \Vert \psi - Q_{c(t)}(x-y(t))- w(t) \Vert_X \to 0 \ \
 \text{ as } t \to \infty . \]
Tao \cite{Tao} verifies scattering  in the following sense: Suppose that 
\begin{equation} \label{energysmall}  \Vert \psi(0) - Q(0) \Vert_{H^1} + \Vert \psi(0) - Q \Vert_{\dot
  H^{-\frac16}} \ll 1. \end{equation} 
Then scattering holds with $X=\dot H^{-\frac16}$. Tao
relies on the work of Martel and Merle, and in particular on Weinstein's 
a priori estimate of the difference to the soliton.

\section{The Airy  Equation}
\label{sec:linkdv}

For purposes of understanding and motivating dispersive estimates for the
linearized KdV equation, here we study and collect results for the Airy  \index{Linear KdV} equation 
\begin{equation} 
\label{eqn:lkdv}
\left\{ \begin{array}{c}
v_t + v_{xxx} = 0 , \\
v(x,0) = v_0 (x).
\end{array} \right.
\end{equation}
The solution operator defines a unitary group $S(t)$ with the kernel \index{$Ai$!fundamental solution}
\[ K(t,x) = t^{-\frac13} Ai(x t^{-\frac13}) , \]
where as $x \to \infty$ the Airy function is roughly
\[    x^{-\frac14} e^{ -  x^{\frac32} } \]
and as $x \to -\infty$ the Airy function is roughly
\[  \text{Re} (  x^{-\frac14} e^{ -i  x^{\frac32} }). \]
Strichartz estimates for solutions,
\begin{equation}
\label{strichartz}
\Vert u \Vert_{L^p L^q} \le c \Vert |D|^{\frac{1}{p}} u_0 \Vert_{L^2} 
\end{equation}
where $L^p L^q$ is the standard space time norm such that the
$L^p$ norm in time of the $L^q$ norm in space and
\[
\frac2p + \frac1q = \frac12,
\]
follow as an immediate consequence.  Of particular interest for this work are the homogeneous Strichartz pair $(p,q) = (6,6)$ as well as the endpoint Strichartz pair $(p,q) = (4,\infty)$.  For an overview of Airy function asymptotics,  see Fedoryuk \cite{Fed}.

Local smoothing estimates for \eqref{eqn:lkdv} go back to the work of Kato \cite{Kato}. Here we are interested in a more general version of them. Let $\gamma(t,x)\ge 1 $ be a smooth bounded increasing function. We calculate
\begin{equation}
\frac{d}{dt} \int \gamma u^2 dx = \int (\gamma_t + \gamma^{(3)})u^2 - 3 \gamma' u_x^2 dx
\end{equation}
and search for conditions ensuring that the right hand side is nonpositive. 
We assume 
\begin{equation}
\label{dyn}
\partial^3_x \gamma \le -\frac23 \partial_t \gamma 
\end{equation}
with the  easiest case being  $\gamma(t,x) = \gamma_0(x-t)$, \index{$\gamma$} for which we assume
\begin{equation}
\label{stat}
\gamma^{(3)}_0 \le \frac23 \gamma'_0.
\end{equation}
We get
\begin{equation}  \frac{d}{dt} \int \gamma u^2 dx + \int \gamma' (u_x^2 + \frac13 u^2) dx \le 0.
\end{equation}
Let us fix a particular example, \index{$\gamma_0$}
\begin{eqnarray}
\label{eqn:gamma0}
\gamma_0(x)  = 1 +   \int_{-\infty}^x (1+|y|^2)^{-\frac{1+\varepsilon}2} dy.
\end{eqnarray}
It satisfies the criteria and, provided $\varepsilon$ is sufficiently small,
a straightforward calculation gives \eqref{stat}. Next, it is instructive to consider
a scaling. For $\mu >0$ and $\gamma_0$ as above we define 
\[   \gamma_\mu (t,x) =  \gamma_0( \mu^{-1}(x-\mu^{-2} t) ). \]
Then,
 \begin{equation}
\label{energy}
\frac{d}{dt} \int \gamma_\mu u^2 dx + \int \gamma_\mu' (u_x^2 +
\frac1{3\mu^2} u^2) dx \le 0.
\end{equation}
One may easily generalize this inequality by choosing $t \to y(t)$ with $\dot y \ge \frac18 \mu^{-2}$, and setting $\gamma(t,x) = \gamma_0 (\mu^{-1} (x-y(t)))$.
In the sequel we will always restrict ourselves to $\mu=1$. 

The virial identity clearly generalizes to functions spaces with different
regularity.  To see this, we first define the space $\Hrho^s$ (and
similarly $L^2_\rho$) by the norm \index{$H^s_\rho$}
\[
 \Vert u \Vert_{\Hrho^s}=  \int |\langle D \rangle^s u|^2 \rho^2 (x) dx < \infty ,
\]

where $\rho > 0$ with uniformly bounded derivatives of order up to $k$ for
some $k \ge  |s|$ and $\langle D \rangle^s$ is defined through the Fourier
multiplication $(1+ |\xi|^2)^{s/2}$  
Similarly  we define $\rhoH^s$ \index{$\rho H^s$} where $u \in \rhoH^s$ if and only if
\begin{eqnarray*}
u = \rho f \ \text{for} \ f \in H^s, \qquad \Vert  u \Vert_{\rhoH^s} = \inf_{u= \rho f} \Vert f \Vert_{H^s}.
\end{eqnarray*}
The function $\rho$ will often depend on $t$. Given a Banach space $X$ we denote the space of $X$ valued $L^2$ functions by $L^2X$, and, with the obvious meaning $L^2 \rhoH^s$ and $L^2 \Hrho^s$. 
Such spaces will be explored further in Section \ref{sec:lin-subcrit}.

\begin{rem}
We note that $\rho H^s = H^s_{\rho^{-1}}$, if $\rho$ is nonnegative, up to equivalent norms. 
However as we wish to highlight the use of duality throughout the
linear analysis and construction of iteration spaces, we adopt the
$\rho H^s$ convention.
\end{rem}

If $\gamma$ satisfies the assumptions above and
\begin{equation}
\label{inhomo}
\left\{ \begin{array}{c}
u_t - u_{xxx} = f, \quad f \in L^2 \gammaH^{-1}  , \\
u(0,x) = u_0 (x), \quad u_0 \in L^2,
\end{array} \right.
\end{equation}
we obtain by an obvious modification of the argument above 
\begin{equation}
\label{eqn:inhomo-est}
\Vert u \Vert_{L^\infty L^2}
+ \Vert  u \Vert_{L^2\Hgamma^1}
\le c \left( \Vert u(0) \Vert_{L^2} + \Vert  f \Vert_{L^2\gammaH^{-1}} \right).
\end{equation}

We turn to a useful  technical result. 
\begin{lem}
\label{lem:commutator}
Let $m \in C^\infty(\mathbb{R})$ satisfy
\[ |m^{(j)}(\xi)| \le c_j \langle \xi \rangle^{s-j} \]
for $j \ge 1$
and let $m(D)$ be the Fourier multiplier defined by $m$.  Suppose that 
\[   \gamma \in C^\infty, \]
\[  | \gamma^{(j)}(x) | \lesssim \gamma(x) \qquad \text{ for } j \ge 0  \]
and 
\[
    \left| 1-\gamma(x)/\gamma(y) \right|    \lesssim c (|x-y| + |x-y|^N ) 
 \text{ for some  } N.     \]
For any $a \in \RR$ we have
\[ \Vert \gamma^{-a} [ m(D) , \gamma^a ]\langle D\rangle^{1-s}  f \Vert_{L^2} 
+ \Vert [ m(D) , \gamma^a ]   \gamma^{-a} \langle D\rangle^{1-s}  f \Vert_{L^2} 
\le c_{s,a} \Vert  f \Vert_{L^2}  \]
and 
\[ \Vert \langle D \rangle^{1-s} [ m(D) , \gamma^a ]  \gamma^{-a} f \Vert_{L^2}+\Vert \langle D \rangle^{1-s}  \gamma^{-a} [ m(D) , \gamma^a ]  f \Vert_{L^2} \le c_{s,a} \Vert f \Vert_{L^2} . \]
\end{lem}

The most important example of $m$ is the Fourier multiplier $\langle D \rangle
^s$
defined by the function $(1+|\xi|^2)^{\frac{s}{2}}$. \index{$\la D \ra$}

\begin{proof}
We begin with the estimate of the first term in the first inequality, the second term being similar. 
We decompose $ m(D) = m_0(D) + m_1(D)$ where the convolution kernel
$m_0(x)$ of $m_0(D)$ is supported in $|x|\le 2$, and the one for
$m_1(D)$ is supported in $|x| \ge 1$. The convolution kernel  $m_1(x)$
together with its derivatives decays exponentially. 

 The integral kernel of
\begin{eqnarray*}
\gamma^{-a} [m_1(D) , \gamma^a]
\end{eqnarray*}
is
\[ K_1 (x,y)  = m_1 (x-y) \left( 1- \left(\frac{\gamma(x)}{\gamma(y)}\right)^{a} \right). \] 
The kernel and its derivatives decay like
$\langle x - y \rangle^{-N}$, which implies  
\[\Vert   \gamma^{-a}[m_1(D), \gamma^a]  f \Vert_{H^N} \le c_N \Vert f \Vert_{H^{-N}}   \]
for all $N>0$ by Schur's lemma. It remains to prove
\[ \Vert \gamma ^{-a} [ m_0(D) , \gamma^a ]\langle D\rangle^{1-s}  f \Vert_{L^2} \le c_{s,a} \Vert  f \Vert_{L^2}.  \]
We decompose 
\[ \langle D\rangle^s = D_0 + D_1. \]
The bound for 
\[ \gamma^{-a} [ m_0(D) , \gamma^a ] D_0  \] 
follows from standard pseudodifferential calculus. The bound for the term with 
$D_1$ follows from 
\[ \Vert \gamma^{-a}[m_0(D),\gamma^a] f \Vert_{L^2} \le c_N \Vert f \Vert_{H^N}, \]
which again follows easily by standard pseudodifferential calculus. 
\end{proof}

\begin{lem} Suppose that 
\begin{equation}
\label{inhomo-sob}
\left\{ \begin{array}{c}
u_t + u_{xxx} = f, \quad f \in \gammaH^{s-1}  , \\
u(0,x) = u_0 (x), \quad u_0 \in H^s.
\end{array} \right.
\end{equation}
Then 
\begin{equation}
\label{sobest}
\Vert u \Vert_{L^\infty H^s}
+ \Vert u  \Vert_{L^2 \Hgamma^{s+1}}
\le c \left( \Vert u(0) \Vert_{H^s} +   \Vert f \Vert_{L^2 \gammaH^{s-1}}\right) .
\end{equation}
Moreover, if 
\begin{equation}
\label{eq}
\left\{ \begin{array}{c}
u_t + u_{xxx} =   ( \sech^{2}(x-x(t)) f)_x + \p_x g , \\
u(0,x) = u_0 (x),
\end{array} \right.
\end{equation}
with $\dot x \ge \delta $, then 
\begin{equation} \label{sobest2}  \Vert u \Vert_{L^\infty {\dot H}^{-1}} +
\Vert  u \Vert_{L^2  \Hgamma^0) }
\lesssim \Vert u(0) \Vert_{\dot H^{-1}} +  
\Vert f \Vert_{L^2 H^{-1}} + \| g \|_{L^1 L^2}.
\end{equation} 
\end{lem}

\begin{proof} 
We set $v= \langle D \rangle^s u$ where $u$ satisfies \eqref{inhomo-sob}, hence 
\[ v_t + v_{xxx} = \langle D \rangle^s f, \]
and, 

\[
\begin{split} 
 \Vert u \Vert_{L^\infty H^s} + \Vert u \Vert_{L^2 \Hgamma^{s+1}} 
= & \Vert v \Vert_{L^\infty L^2} + \Vert v \Vert_{L^2 \Hgamma^1} 
\\ \le & c ( \Vert v(0) \Vert_{L^2} + \Vert \langle D \rangle^s f \Vert_{L^2 \gamma 
H^{-1} }), 
\end{split} 
\]
where the first term is equal to $\Vert u (0) \Vert_{H^s}$ and 
\[
\begin{split} 
 \Vert \langle D \rangle^s f \Vert_{L^2 \sqrt{\gamma'} 
H^{-1} } = & \Vert   \langle D \rangle^{-1} (\gamma')^{-\frac12}   
   \langle D \rangle^{s}  f \Vert_{L^2L^2} 
\\ \le &  \Vert  f \Vert_{L^2  \gammaH^{s-1} } + 
\Vert [\langle D \rangle^{-1}, (\gamma')^{-\frac12} ] \langle D \rangle^{s} f \Vert_{L^2L^2} 
\\ \le &  \Vert  f \Vert_{L^2 \gammaH^{s-1} } + 
\Vert  (\gamma')^{-\frac12}  \langle D \rangle^{s-2} f \Vert_{L^2 L^2 }. 
\end{split} 
  \]
The last inequality follows from Lemma \ref{lem:commutator} applied with 
$\gamma'$ for $\gamma$, $a= \frac12$,  and $-1$ for $s$.   
This implies the desired estimate \eqref{sobest}. 
Now suppose that $u$ satisfies \eqref{eq} and let $v$ be the antiderivative of $u$  with respect to $x$.  It satisfies 
\begin{equation}
\left\{ \begin{array}{c}
v_t + v_{xxx} =   \sech^{2}(x-x(t)) f +  g , \\
v(0,x) = v_0 (x),
\end{array} \right.
\end{equation}
hence 
\[ \Vert u \Vert_{L^\infty \dot H^{-1} } + \Vert u \Vert_{L^2 \Hgamma^0} 
\le c \Vert g \Vert_{L^1 L^2} + \Vert f \Vert_{L^2 H^{-1}}. \]
\end{proof}

\section{Properties of the  Schr\"odinger Operator}
\label{sec:lin-subcrit}

We briefly recall notions from the introduction. 
Given $p>1$ solitons of the form $Q_p(x-t)$ satisfy \eqref{euler} 
and it is not hard to verify that all bounded solutions are translates of 
$\pm Q_p$  in equation \eqref{Qps}. Similarly $Q_{p,c} =  c^{\frac{2}{p-1}} Q_p ( c x)$
satisfies 
 \begin{equation}
\label{eqn:c-sol}
\p_x^2 (Q_p)_c - c^2 (Q_p)_c + (Q_p)_c^p = 0 .
\end{equation}
We will focus on $p=4$ and omit again $p$ from the notation. Let $({}')$ denote
the differentiation with respect to $x$ and $(\dot{})$ the differentiation with
respect to time. We recall the definition of $\tilde Q$, \eqref{tildeQ}  
and $\tilde Q_c = c \partial_c Q_c $ respectively $ \tilde {\tilde Q}_c= c\partial_c c \partial_c Q_c $ the corresponding 
differentiation at $c$. There are many explicit calculations, and we collect
some of them here. 
Using the properties of $Q_c (x) = c^{\frac{2}{3}} Q (cx)$, it follows that
\begin{equation}
\label{eqn:tildeqq}  
\Vert Q_c \Vert_{L^2} = c^{\frac16} \Vert Q_1 \Vert_{L^2}, \quad \langle
\tilde Q_c, Q_c \rangle = 
\frac12 c\partial_c  \Vert Q_c \Vert_{L^2}^2  = \frac1{6} \Vert Q_c
\Vert_{L^2}^2,\end{equation} 
where the $L^2$ norm is given by \eqref{mul2}, and
\begin{equation}\label{eqn:qprime} 
 \Vert Q_c' \Vert_{L^2} = c^{\frac76} \Vert Q'_1   \Vert_{L^2}. 
\end{equation} 
In addition,
\[
\p_x Q_c  =  c^{\frac{5}{3}} Q' (cx),
\]
\[
c\p_c Q_c  =  (\frac{2}{3} Q_c + x Q_c') =  \tilde{Q}_c  =    c^{\frac23}  \tilde Q(cx).
\] 

 The operator $\Lop_c$ is defined by  
\begin{equation} \label{Lop} 
\Lop_c u = - u_{xx} + c^2 u - 4 Q^3_c u ,
\end{equation} 
where we mostly omit $y$ and $c$ if $c=1$. We recall that vitually everything 
is known about the spectrum of $\Lop$, 
see
Andrews-Askey-Roy \cite{AAR}, Lamb \cite{Lamb} and Titchmarsh \cite{Titchmarsh}.  We  summarize the findings below.  We also
refer to Martel \cite{Martel}, Weinstein \cite{W} and the references therein for
extensive discussions of these properties for more general operators
of type similar to $\Lop$.

By direct differentiation in $x$ of \eqref{euler}, we see $\Lop
Q' = 0$.  Hence, the null space of $\Lop$ consists at least of the
space $\alpha Q'$ for all $\alpha \in \RR$.  Similarly, by
differentiation in $c$ of \eqref{eqn:c-sol}, we see $\Lop (\tilde{Q})
= -2Q$, so $\p_x \Lop$ has at least a $2$-dimensional generalized null
space.  Also, since $Q' = 0$ only at $x=0$, we know from the Sturm
Oscillation Theorem that there exists some $\lambda_0 > 0$,
$\mathcal{Q}_0 > 0$ such that $\Lop \mathcal{Q}_0 = - \lambda_0
\mathcal{Q}_0$, the unique negative eigenstate of $\Lop$.  Note,
because $\Lop$ is a $\sech^2$ potential perturbation of the Laplacian,
it is possible to exactly construct $\mathcal{Q}_0 = Q^{\frac52}$
\index{$\mathcal{Q}_0$}and $\lambda_0 = \frac{21}{4}$
\index{$\lambda_0$}using standard techniques. The above analysis
summarizes the entire discrete spectral decomposition for $\Lop$.

Following the introduction resp. the analysis in Propositions 2.7 and 2.9 of Weinstein \cite{W}, if
\begin{eqnarray*}
\langle \tilde{u}, Q \rangle = 0
\end{eqnarray*}
and
\begin{eqnarray*}
\langle \tilde{u}, Q' \rangle = 0,
\end{eqnarray*}
then there exists $k_0 > 0$ such that
\begin{eqnarray}
\label{eqn:Loppos}
\langle \tilde{u} , \Lop \tilde{u} \rangle \geq k_0 \| \tilde{u} \|_{L^2}^2.
\end{eqnarray}
Here $k_0$ depends only on the power $p=4$ in \eqref{eqn:kdv}.

We will consider $\rho= e^\nu$ with $\nu \in C^{|s|+1} $ with \index{$\nu$}
\begin{equation}  \label{nu} 
 |\nu^{(j)}(x)| \le \varepsilon 
\end{equation}
for $0 \leq j \leq |s| + 1$ and a small constant $\varepsilon$ to be chosen
later. Clearly we may regularize $\nu$ and hence $\rho=e^\nu$ without changing the
spaces. Then
\[ u \in \Hrho^s \Leftrightarrow \rho u \in H^s
\Leftrightarrow u \in \rho^{-1} H^s. \] 
 It is quite obvious that the dual space of $\Hrho^s$ is $\rhoH^{-s}$ with isometric norms, and this statement does not depend on the
regularity of $\rho$. We recall the definition of the projectors \eqref{projection}.   
\begin{lem}
\label{lem:LopXYbounds}
For all $s \in \mathbb{R}$, there exists $C >0$ such that:
\[
\Vert P_{Q'}^\perp u \Vert_{H^{s+2}}  \le C \Vert \Lop u \Vert_{H^s} ,
\]
\[
\Vert P_{Q'}^\perp u \Vert_{\rhoH^{s+2}} \le C \Vert \Lop u \Vert_{\rhoH^{s}} ,
\]
\[
\Vert P_{Q'}^\perp u \Vert_{ \Hrho^{s+2}} \le C \Vert \Lop u \Vert_{ \Hrho^s}  .
\]
\end{lem}

\begin{proof}
The first inequality is an immediate consequence of the nature of the spectrum
described above along with ellipticity. The second and the third statement are equivalent because  $H^s_{\rho} = \rho^{-1} H^s$, with equivalent norms.

Fix $\mu = 1- \frac{(p+1)^2}4$,  where $p=4$.
For $\lambda= \lambda_0 + i \lambda_1$ in the complex half plane left of $\mu$ we obtain the following
resolvent estimate
\[ |\lambda -\mu| \Vert u \Vert_{L^2} \le \Vert (\Lop - \lambda)
u \Vert_{L^2} \]
and also for some $1>\kappa > 0$, we have
\begin{eqnarray*}
\Re \int u \overline{ ( \Lop- \lambda) u} dx & \ge & 
| \lambda - \mu| \| u \|_{L^2}^2 + \langle (\Lop - \mu) u, u \rangle \\
& \ge & \frac{1}{2} |\lambda -mu| \Vert  u \Vert_{L^2}^2 
+ \kappa \|  u_x \|_{L^2}^2  +  (1-\kappa) \langle (\Lop-\mu)u, u \rangle  \\
& & + \left(\frac12 |\lambda-\mu| + \kappa(1-\mu)  - 4\kappa \Vert Q \Vert_{L^\infty}^3)\right) \Vert u \Vert_{L^2}^2  \\
& \ge & \frac{1}{2} |\mu -\lambda_0| \Vert  u \Vert_{L^2}^2 + 
\min \{ \frac{|\mu-\lambda_0|}{8 \Vert Q \Vert_{L^\infty}^3},\frac12 \} 
 \| u_x \|_{L^2}^2
\end{eqnarray*}
by the obvious choice of $\kappa$. 

We obtain the estimate for $\lambda$ with real part at most $\mu$,
\[  |\mu-\lambda| \Vert u \Vert_{L^2} + \min \{ |\lambda- \mu|,1\} 
\Vert u_x \Vert_{L^2} 
 \le C  \Re \langle (\Lop - \lambda)u, u \rangle 
\le C \Vert  \langle (\Lop - \lambda)u \Vert_{H^{-1}} 
\Vert u \Vert_{H^1} .
\]
These estimates imply that the resolvent $(\Lop - \lambda)^{-1}$ defines a
continuous uniformly bounded map (for $\Re \lambda \le \lambda_0 <  \mu$) from $H^{-1}$ to $H^1$.
Moreover, 
\[ \Vert u \Vert_{L^2} \le |\lambda-\mu|^{-1} \Vert \Lop u \Vert_{L^2} \]
and 
\[ \Vert u_x \Vert_{L^2}^2 \le |\lambda-\mu|^{-1}\max \left\{ \frac12, \frac{8 \Vert Q
  \Vert_{L^\infty} }{ |\lambda-\mu|} \right\} \Vert (\Lop -\lambda)u \Vert_{L^2}^2. \]

We turn to the weighted estimates and calculate formally
\[ e^{\nu} (\Lop - \lambda) e^{-\nu} = \Lop - \lambda - |\nu'|^2 +
\partial_x \nu' + \nu' \partial_x, \]
and hence, since $ \p_x \nu' + \nu' \p_x$ is antisymmetric,
\[
\begin{split}
\Re \int u \overline{e^{\nu} (\Lop - \lambda) e^{-\nu} u )}
dx = &    \Re \int u \overline{(\Lop -\lambda) u} dx - \Vert \nu' u \Vert_{L^2}^2
\\
\ge  \frac{|\lambda_0-\mu|}2 \Vert u \Vert_{L^2}^2
\end{split}
\]
if $ \varepsilon \le \sqrt{\frac{|\lambda_0-\mu|}{2}}$, which we assume in the
sequel. As above  we obtain  with an explicit constant $C$
\begin{equation}\label{resolvent1} 
\Vert u \Vert_{H^1} \le C \Vert e^{\nu} (\Lop - \lambda) e^{-\nu} u
\Vert_{H^{-1}}.
\end{equation} 
It follows from these estimates that given $\delta>0$ 
there is a single resolvent family (for
$\Re \lambda < \mu-\delta$) mapping $ \rhoH^{-1} \to \rhoH^1 $ and from
$\Hrho^{-1} \to \Hrho^1$ provided $\varepsilon$ is sufficiently small. 

Recall that $ \Lop $ has a zero eigenvalue with eigenfunction $Q'$
and a single negative eigenvalue $-\lambda_0$ with a ground state $\mathcal{Q}_0$. Let $P$ be the orthogonal projection to the
orthogonal complement of these two eigenfunctions. The remaining
spectrum is contained  in $[\rho,\infty)$ where $\rho >0$ is either $1$ (if $p
\ge 3$), or the next positive eigenvalue, which can be easily be calculated. 
  Moreover, $\Lop$ is selfadjoint. The resolvent
  $R(\lambda)=(\Lop-\lambda)^{-1}$ is a holomorphic map in $\mathbb{C}
  \backslash (1,\infty)$ with simple poles in $\mu$, $0$, and possibly some 
 other eigenvalues in $(0,1)$. 
In addition, $R_0(\lambda) = R(\lambda)P$ has a continuous and hence holomorphic extension to $\lambda = 0$ and $\lambda = -\lambda_0$, which is uniformly bounded in each half plane strictly left of $\rho$.

By equation \eqref{resolvent1} the resolvent is uniformly bounded on
the weighted spaces if $\lambda$ is in the half plane left of $-\mu$.
Decreasing $\varepsilon$ if necessary (so that the orthogonal projection
$P^\perp_{Q'}$ along $Q'$ is bounded in the weighted space) we obtain the same statement for $R_0(\lambda)$. Now complex
interpolation implies
\[ \Vert   \Lop^{-1} P   f \Vert_{\Hrho^1} \le C \Vert Pf \Vert_{\Hrho^{-1}} . \]
This implies the desired estimates for $s = -1$.

Standard elliptic theory extends this estimate to
\begin{equation}
\label{resolvent}
\Vert u \Vert_{\Hrho^{s+2}} \le C \Vert (\Lop - \lambda) u \Vert_{\Hrho^s}
\end{equation}
and 
\begin{equation}
\label{resolvent2}
\Vert u \Vert_{\rhoH^{s+2}} \le C \Vert (\Lop - \lambda) u \Vert_{\rhoH^s}
\end{equation}
first to all $s \ge -1$, and then, by duality, to all $s \in \mathbb{R}$.
The first estimate is the special situation when $\nu$ is constant.

We conclude with the trivial observation that we may replace 
\eqref{nu} by 
\[ 
\lim_{x\to \in \infty} \nu^{j} = 0 ,
\]
which holds for $\rho (x) = (1+|x|^2)^a$ for all real numbers $a$, since in
that  case we may choose an equivalent norm which satisfies \eqref{nu}.

\end{proof}

\section{Energy Methods for the linearized  equation}
\label{sec:virial}

 We turn to a study of what we call the linear $u$-problem
\index{$u$ problem}
\begin{eqnarray}
\label{eqn:forward}
\left\{ \begin{array}{c}
u_t = \p_x (\Lop u), \\
u(0,x) = u_0,
\end{array} \right.
\end{eqnarray}
where 
\begin{eqnarray*}
\Lop u & = & ( -\p_x^2 + 1 - 4 Q^3 ) u  \\
& = & ( -\p_x^2  + 1 - 10 \sech^2 (\frac{3}{2} x) ) u .
\end{eqnarray*}
We note here that $\Lop$ \index{$\Lop$}is the operator that results from linearization of the KdV equation about $Q$ when we work in a moving reference frame or in other words make the change of variables 
\[ x \to x-t. \]
Indeed, setting $\psi(x,t) = Q(x-t) + u(x-t,t)$ and plugging into \eqref{eqn:kdv}, we get 
\begin{eqnarray*}
\p_t u & = & - \p_x(\p_x^2 u - u +(Q+u)^4-Q^4 + \p_x^2 Q -Q +Q^4) \\
& = & \p_x (\Lop u) - \p_x ( 6 Q^2 u^2 + 4 Q u^3 + u^4).
\end{eqnarray*}
For reasons that will become clear in the sequel, we also consider the linear $v$-problem \index{$v$ problem}
\begin{eqnarray}
\label{eqn:linv}
\left\{ \begin{array}{c}
v_t =  \Lop( \p_x v) , \\
v(0,x) = v_0.
\end{array} \right.
\end{eqnarray}

The two equations \eqref{eqn:linv} and \eqref{eqn:forward} are related in many ways. 
\begin{enumerate} 
\item They are dual equations of each other. 
\item If $u$ satisfies the $u$ equation then $v= \partial_x u$ satisfies the $v$ equation. 
\item If $v$ satisfies the $v$ equation then $u= \Lop v$ satisfies the $u$ equation. 
\end{enumerate} 

We observe that $u=Q'$ is a solution to the $u$ equation, and hence
 $\langle v, Q' \rangle$  is preserved by the flow for $v$. In particular 
orthogonality is preserved by the evolution. 
Similarly $v=Q$ is a solution to the $v$ equation and $\langle u,Q \rangle $ 
is preserved by the $u$ flow. Moreover $u= a Q' + b (\tilde Q +2 t Q')$
satisfies the $u$ equation for all coefficients $a$ and $b$. As a consequence 
both equations admit solutions which grow linearly with time. Moreover, if $v$ 
satisfies the $v$ equation then 
\[ \frac{d}{dt} \langle v ,\tilde Q \rangle + 2t \langle v, Q' \rangle = 0 \]
and $v$ is orthogonal to $\tilde Q$ and $Q'$ provided it is initially.

Inspired by a set of ideas collected from Martel-Merle
\cite{MM3} and the references therein, let us look at a virial
identity for \eqref{eqn:linv} given by \index{$I_\eta$!constant
  coefficient}
\begin{eqnarray*}
I_\eta (v) = - \int \eta (x) v^2 dx,
\end{eqnarray*}
where $\eta (x)$ will be defined in the sequel.  We have
\begin{eqnarray*}
- \frac{d}{dt} I_\eta (v) & = & - 2  \int \eta (x) v ( \Lop v_x)  dx \\
& = & - 2  \int \eta (x) v ( (-\p_x^2 + 1 - 4 Q^3) v_x)  dx \\
& = & 2 \int \eta v \p_x^3 v dx - 2 \int \eta v \p_x v dx+ 8 \int \eta Q^{3} v \p_x vdx \\
& = & -2 \int \eta' v \p_x^2 v dx - 2 \int \eta \p_x v \p_x^2 v dx + \int \eta_x v^2 dx \\
& & -4  \int \eta' Q^3 v^2 dx - 4  \int \eta \p_x \left( \sech^2 \left( \frac{3}{2} x \right) \right) v^2 dx \\
& = & 3 \int \eta' v_x^2 dx + 2 \int \eta'' v \p_x v dx + \int \eta' v^2 dx \\
& & -4  \int \eta' Q^3 v^2 dx - 12  \int \eta Q^{2} Q' v^2 dx .
\end{eqnarray*}

As in the work of Martel \cite{Martel}, we take \index{$\eta$!constant
  coefficient}
\begin{eqnarray}
\label{eqn:eta}
\eta (x) = -\frac{5}{3} \frac{Q'}{Q} = \frac53 \tanh\Big(\frac32 x\Big)  ,
\end{eqnarray}
which is similar to $x$ near $0$ and bounded at $\infty$.  Note, the
sign convention here is chosen to match that of \cite{MM3}.
By direct computation we have
\begin{eqnarray*}
\eta' (x) = Q^3 (x),
\end{eqnarray*}
\begin{eqnarray*}
\left( \frac{ \eta'' (x) }{ \eta'(x) } \right)^2 = 9 \left( 1 - \frac{2}{3} Q^3 (x) \right),
\end{eqnarray*}
\begin{eqnarray*}
\frac{\eta''' (x)}{\eta' (x)} = 9 \left( 1 - \frac{3}{5} Q^{3} (x) \right),
\end{eqnarray*}
\begin{eqnarray*}
\eta^2 (x) = \left( \frac{5}{3} \right)^2 \left( 1 - \frac{2}{5} Q^{3} (x) \right)
\end{eqnarray*}
and
\begin{eqnarray*}
( Q^3 \eta)' = -5 Q^3 + 3 Q^6, \ | \eta | \leq \frac{5}{3}.
\end{eqnarray*}

\begin{prop}
\label{prop:virial}
If $v$ satisfies the $v$-KdV equation and $v \perp \tilde{Q},Q'$, then there exists some $C > 0$ such that given $\eta$ as in \eqref{eqn:eta} we have
\begin{eqnarray*}
\frac{d}{dt} I_\eta (v) + C \|  \sech\Big(\frac32 x\Big)  v \|^2_{H^1} \le 0 .
\end{eqnarray*}
\end{prop}

\begin{proof}[Proof of Proposition \ref{prop:virial}]

Following the formalism presented above, we see
\begin{eqnarray*}
-\frac{d}{dt} I_\eta (v) & = & - 2 \int \Lop (\p_x v) v \eta dx \\
& = & 3 \int (\p_x v)^2 \eta'\, dx  + \int v^2 \left[ -\eta''' + \eta' - 4
  (Q^3 \eta)' \right]\, dx .
\end{eqnarray*}
Selecting \index{$\tilde{w}$!constant coefficient}
\begin{eqnarray*}
\tilde{w} (t,x) = v(t,x) \sqrt{\eta'(x)}
\end{eqnarray*}
we see
\begin{eqnarray*}
-\frac{d}{dt} I_\eta (v) = 3 \int (\p_x  \tilde{w})^2 dx +  \int A(x) \tilde{w}^2 dx,
\end{eqnarray*}
where \index{$A$!constant coefficient}
\begin{eqnarray*}
A(x) & = & 1 + \frac{1}{2} \frac{ \eta'''}{\eta'} - \frac{3}{4} \left( \frac{\eta''}{\eta'} \right)^2 - 4 \frac{(Q^3 \eta)'}{\eta'} \\
& = & \frac{75}{4} - 12 Q^3.
\end{eqnarray*}
Hence,
\begin{eqnarray*}
-\frac{d}{dt} I_\eta (v) = 3 \left[ \langle \Lop \tilde{w}, \tilde{w} \rangle + \frac{21}{4} \int \tilde{w}^2 dx \right].
\end{eqnarray*}
Since $\Lop \p_x Q =0$, we know that given $v = Q$, we have
\begin{eqnarray*}
-\frac{d}{dt} I_\eta (v) = 0 .
\end{eqnarray*}
However, $v=Q$ corresponds directly to $\tilde{w} = Q^{\frac{5}{2}}$, which is the ground state or $\Lop$, which has exact eigenvalue $-\frac{21}{4}$.  
Then, since $\langle Q, \tilde Q \rangle \neq 0$, our orthogonality condition
\begin{eqnarray*}
v \perp \tilde{Q}
\end{eqnarray*}
is enough to guarantee that there exists $C > 0$ such that
\begin{eqnarray*}
\tilde{B} (\tilde{w},\tilde{w}) \geq C \| \tilde{w} \|^2_{H^1} = C \| \sqrt{\eta'} v \|^2_{H^1},
\end{eqnarray*}
which is the desired result.
\end{proof}

We note in the case of more general weight functions $\eta$, virial
identity methods are still applicable even if perhaps analytic proofs
of the virial identities are more challenging.  

By choosing the multiplier $\gamma (v-v_{xx})$ with $\gamma = \gamma_0 (x-t)$ for $\gamma_0$ as in \eqref{eqn:gamma0}, we see
\begin{eqnarray}
\frac{d}{dt} \int \gamma (v^2 + v_x^2) dx & = & - 3 \int \gamma' v_x^2
dx + \int \gamma^{(3)} v^2 dx - \int \gamma' v^2 dx \notag \\
\label{eqn:vxxmult}
& & + \ 4 \int \gamma' Q^3 v^2 dx + 12  \int \gamma Q^2 Q' v^2 dx \\
& & - \ \int [ 3  \gamma' v_{xx}^2 + \gamma' v_x^2  - \gamma^{(3)}
v_x^2 ] dx  + \ \int 4 \gamma' Q^3 v_x^2 dx, \notag
\end{eqnarray}
which consists of a number of negative semidefinite terms. 
 All non negative semidefinite terms 
 are easily dominated by a multiple of $\Vert v \Vert_{H^1_{\sech(3x/2)}}^2 $,
 the term in \eqref{prop:virial}.

Finally, note that by direct computation
\begin{eqnarray}
\label{eqn:invariant}
\p_t \langle \Lop^{-1} v,v \rangle = 0.
\end{eqnarray}

Now, let us define an energy for the solution $v$ of \eqref{eqn:linv} to be
\begin{eqnarray}
\label{eqn:energy}
E (v) = \int \gamma (x) (v^2 + v_x^2) dx + \lambda_E \int \eta(x) v^2 dx + \Lambda_E \langle \Lop^{-1} v,v \rangle ,
\end{eqnarray} \index{$\lambda_E$} \index{$\Lambda_E$}
where $\eta(x)$ is chosen as in \eqref{eqn:eta}. \index{$E$!constant coefficient}

\begin{prop}\label{l2const} 
Let us assume $v$ satisfies the $v$-KdV equation and 
$v \perp \tilde{Q},Q'$.  There exist $\lambda_E, \Lambda_E, \delta >0$ such that
\begin{equation}
\label{eqn:energy1}
E(v) \sim  \| v \|^2_{H^1}
\end{equation}
 and
\begin{equation}
\label{eqn:energy2}
\frac{d}{dt} E(v) + \delta \| v \|^2_{H^2_{\sqrt{\gamma'}}}\le 0 .
\end{equation}
\end{prop}

\begin{proof}
  From \eqref{eqn:vxxmult} and the proof of Proposition \ref{prop:virial}, we
  see easily one may choose a $\lambda_E$ which depends only on $\delta$ and
  $C$
so that \eqref{eqn:energy2} holds for all $\Lambda_2>0$. We choose $\Lambda_E$
large to achieve $E(v) \geq C'' \| v \|_{H^1}^2$.
There exists some constant
  $C'$ such that $E(v) \leq C' \| v \|_{H^1}^2$.
  estimates given the orthogonality conditions on $v$. 
\end{proof}

The assertions of Proposition \eqref{l2const} are robust under suitable perturbations. We turn to  the analysis of the time dependent problem
\begin{eqnarray}
\label{eqn:vtd}
v_t -  (-\p_x^2 - 4 Q_{c(t),y(t)}^3) \p_x v = \alpha(t) Q_{c(t),y(t)} + \beta(t) Q_{c(t),y(t)}' ,
\end{eqnarray}
where
\begin{eqnarray}
\label{eqn:alpha-subtd}
\alpha (t) = - \frac{ \frac{\dot{c}}{c} \langle v, \tilde{\tilde{Q}}_{c(t),y(t)} \rangle + (\dot{y}-c^2) \langle v, \tilde{Q}_{c(t),y(t)}' \rangle }{\langle Q_{c(t),y(t)} , \tilde{Q}_{c(t),y(t)}  \rangle}
\end{eqnarray}
and
\begin{eqnarray}
\label{eqn:beta-subtd}
\beta (t) = - \frac{ \frac{\dot{c}}{c} \langle v , \tilde{Q}_{c(t),y(t)}' \rangle + (\dot{y}-c^2) \langle v, \tilde{Q}_{c(t),y(t)}'' \rangle }{\langle Q_{c(t),y(t)}' , Q_{c(t),y(t)}' \rangle }  .
\end{eqnarray}
Here,
\begin{eqnarray}
\label{eqn:2tilde}
\tilde{\tilde{Q}}_{c(t),y(t)} = \frac23  \tilde{Q}_{c(t),y(t)} + x  \tilde{Q}_{c(t),y(t)}' = c(t) \partial_c \tilde{Q}_{c(t),y(t)}.
\end{eqnarray}
For simplicity of exposition, in the sequel we suppress the $t,y$
dependence and write simply $Q_{c(t),y(t)} = Q_c$ unless we want to
stress the dependence on $y(t)$ respectively on $t$. Similarly we recall 
\begin{equation}   
\Lop_c v = \Lop_{c,y} v = -v_{xx} + c^2 v -  4 Q_{c,y}^3 v.
\end{equation} 

The terms on the right hand side ensure that $\langle v(0), Q_{c(0),y(0)}' \rangle =0 $
 implies  $\langle v(t), Q_{c(t),y(t)}' \rangle =0 $, and, in addition,  
$\langle v(0), \tilde Q_{c(0),y(0)} \rangle =0 $ implies 
$\langle v(t), \tilde Q_{c(t),y(t)} \rangle =0$.
We choose $\gamma(x,t) = \gamma_0(x-y(t))$ and we  prove the following
\begin{prop}
\label{prop:tdenergy}
There exists a $\delta, \lambda, \Lambda > 0$ such that the following is true: Suppose that
\begin{equation}\label{eqn:cysmall} 
|c(t)-1|+ |\dot{c} (t)|+ |\dot y (t) - c^2 (t)| < \delta
\end{equation}
for all $t \geq 0$ and define 
\begin{equation}
\label{eqn:energytd}
E (v) = \int  \gamma (x,t) (v^2+v_x^2) dx + \lambda \int \eta_{1,y(t)} (x) v^2 dx + \Lambda \langle \Lop_{1,y(t)}^{-1} v,v \rangle,
\end{equation}
where we suppress the dependence  of $E$ and $v$ on $t$. 
Then
\begin{equation}
  \label{eqn:tdenergy1}
E(v)  \sim  \| v \|^2_{H^1}  , 
\end{equation} 
for all $t>0$ provided   
\begin{equation} 
 \langle v, Q'_c\rangle = \langle v, \tilde Q_c \rangle =0.
\end{equation} 
Moreover, if 
 $v$ satisfies the system consisting of 
 \eqref{eqn:vtd}, \eqref{eqn:alpha-subtd} and \eqref{eqn:beta-subtd} 
and  $v(.,0) \perp \tilde{Q}_{c(0),y(0)} ,Q_{c(0),y(0)}'$ (which implies the
orthogonality for all $t$) we have  
\begin{equation} 
\label{eqn:tdenergy2}
\frac{d}{dt} E(v) + \delta \| v \|^2_{H^2_{\sqrt{\gamma'} }}\le 0 .
\end{equation}
\end{prop}
 
\begin{proof} 

Since $\langle v(t), Q'_{c(t),y(t) } \rangle =0$ we have
\begin{eqnarray*}
\langle \Lop_{c(t),y(t)}^{-1} v,v \rangle \geq C \| v \|^2_{H^{-1}} ,
\end{eqnarray*}
for some $C > 0$ as seen in \eqref{eqn:Loppos}. Here and in the 
remaining part of this section we use the Moore-Penrose 
inverse, which is by an abuse of notation given the orthogonal projection 
to the complement of $Q'$, followed by an inversion of $\Lop$ on this 
orthogonal subspace. 
Let us look at a slightly different quantity (where we replace $c$ by $1$)
given by
\begin{eqnarray*}
\langle \Lop_{1,y(t)}^{-1} v,v \rangle.
\end{eqnarray*}
Then, since $\langle v, Q_c' \rangle = \langle v, \tilde{Q}_c \rangle =0$, for $|c-1|$ small enough we have
\begin{eqnarray*}
\langle \Lop_{1,y(t)}^{-1} v,v \rangle & \geq & 2C \| P^\perp_{Q_{1,y(t)}'} v \|^2_{H^{-1}} \\
& \geq & 2C \| v \|^2_{H^{-1}} - C' |c-1| \| v \|_{H^{-1}}^2 \\
& \geq & C \| v \|^2_{H^{-1}},
\end{eqnarray*}
for some constants $C, C'>0$ and $\delta \le C/C'$. The properties are similar to the
previous proposition, but the calculations are more tedious. We
consider them to be important for the understanding of the
linearization. We recall that we suppress the dependence of $Q$ and $\Lop$ on 
$y$ in the notation below.  Then, we have
\begin{eqnarray*} 
\frac{d}{dt} \langle \Lop_{1}^{-1} v,v \rangle &  = & 2 \langle  v_t, \Lop_{1}^{-1} v \rangle - 12 \dot{y} \langle Q_{1}^2 Q_{1}' \Lop_{1}^{-1} v, \Lop_{1}^{-1} v \rangle 
+ 2 {\dot y} \frac{\langle v , Q_{1}' \rangle}{\langle Q_{1}' , Q_{1}' \rangle} \langle Q_{1,y}'' , \Lop_{1}^{-1} v \rangle \\
& = & 2I_1-12I_2 +2 I_3 ,
\end{eqnarray*} 
where $I_2$ originates from the differentiation of the inverse and $I_3$ 
from the dependence of the implicit projection on time. We have
\begin{eqnarray*} 
I_1 & = &  \langle \Lop_{c} \p_x v ,\Lop_{1}^{-1} v \rangle 
 - c^2 \langle \p_x v, \Lop_{1}^{-1} v \rangle + \alpha \langle Q'_c, \Lop^{-1}_1 v \rangle + \beta \langle Q_c , \Lop^{-1}_1 v \rangle ,
\\
\langle \Lop_{c} \p_x v ,\Lop_{1}^{-1} v \rangle &=& \langle
(\Lop_{c}- \Lop_{1}) \p_x v ,\Lop_{1}^{-1} v \rangle = (c^2-1) \langle
\p_x v, \Lop^{-1}_1 v \rangle + 4\langle (Q_c^3-Q_1^3) \p_x v, \Lop^{-1} v
\rangle ,
\\
   \left[ \Lop_1^{-1},\p_x \right] & = & -\Lop_1^{-1} [ \Lop_1,\p_x] \Lop_1^{-1}= 12 \Lop^{-1}_1 Q_1^2 Q_1' \Lop_1^{-1} + \Lop_1^{-1} [ \p_x, P_{Q'_1} ] + [\p_x, P_{Q'_1}] \Lop_1^{-1} ,   
\\
 \left[ \partial_x , P_{Q'_1} \right]  v & =& -\frac{\langle v, Q_1'\rangle}{\langle Q_1' , Q_1' \rangle} Q_1''  + \frac{\langle v, Q_1''\rangle}{\langle Q_1' , Q_1' \rangle} Q_1' ,
\\   \langle \p_x v, \Lop_1^{-1} v \rangle 
&=&  \frac12 \langle v, [ \Lop_1^{-1}, \p_x ]  v \rangle = 
  6 \langle Q_1^2 Q_1' \Lop_1^{-1} v, \Lop_1^{-1} v \rangle,
 \\
\langle v, Q_1' \rangle & = &  \langle v, Q_1' -Q_c' \rangle  
\end{eqnarray*} 
by the orthogonality conditions and  
\[ \langle \Lop^{-1}_1  v, Q_c' \rangle = \langle \Lop^{-1}_1 v, Q_c' -Q_1' \rangle , \]
\[ \langle \Lop^{-1}_1 v, Q_c  \rangle  = \langle \Lop_1^{-1} v, (\Lop_1^{-1} - \Lop_c^{-1})  Q_c \rangle
\]
because of the orthogonality conditions and since $\Lop^{-1}_c Q_c= \tilde Q_c$.    

Altogether,  and applying Lemmas \ref{lem:LopXYbounds}, we have 
\begin{equation}
|\frac{d}{dt} \langle \Lop_{1}^{-1} v,v \rangle| \leq \mathcal{O} (|c^2-1| + |\dot{y}-1| + |\dot c|) \| v \|_{H^{-\frac12}_{\sqrt{\gamma'}}}^2,
\end{equation}
which we will control by the virial identity below.

We now look at virial weights of the form \index{$\eta$!variable coefficient}
\begin{eqnarray*}
\eta(x,t) = -\frac{5}{3} \frac{Q_{1}' (x-y(t))}{Q_{1} (x-y(t))}= \frac53
\tanh( \frac32(x-y(t))) ,
\end{eqnarray*}
which has similar properties to those of $\eta(x)$ with appropriate changes for the unit scaling.

We have defined $v$ such that
\begin{eqnarray*}
v(x,t) \perp \tilde{Q}_{c(t),y(t)}
\end{eqnarray*}
and 
\begin{eqnarray*}
v(x,t) \perp Q_{c(t),y(t)}'
\end{eqnarray*}
for all $t \ge 0$.
Following the formalism presented above and in Martel \cite{Martel},
\index{$I_\eta$!variable coefficient} select \index{$\tilde{w}$!variable coefficient}
\begin{eqnarray*}
\tilde{w}(t,x) = v(t,x) \sqrt{\eta'(x)}.
\end{eqnarray*}
Then, 
\begin{eqnarray*}
-\frac{d}{dt} I_\eta (v) & = & 3 \int (\p_x \tilde{w} )^2 dx +  \int A(x) \tilde{w}^2 dx  - 2 \int \eta_{1,y} \left( \beta(t) Q_{c,y}' + \alpha (t) Q_{c,y} \right) v dx \\
&& \hspace{-19pt} + \frac{3}{2} c^2 (\dot{y} - c^2) \int \sech^2  \left( \frac{3}{2} c(x - y(t)) \right) v^2 dx  + \frac{3}{2} c^4 \int \sech^2  \left( \frac{3}{2} c(x - y(t)) \right) v^2 dx,
\end{eqnarray*}
where \index{$A$!variable coefficient}
\[
A(x,t)  =  1 + \frac{1}{2} \frac{ \eta'''}{\eta'} - \frac{3}{4} \left( \frac{\eta''}{\eta'} \right)^2 - 4 \frac{(Q^3 \eta)'}{\eta'} 
 =  \frac{75}{4} - 12 Q^3.
\]
Hence,
\begin{eqnarray*}
-\frac{d}{dt} I_\eta (v) > 3 \left[ \langle \Lop \tilde{w}, \tilde{w} \rangle + \frac{21}{4} c^2 \int \tilde{w}^2 dx \right] + \mathcal{O} (|1^2-c^2| + |\dot{y}-c^2|) \| v \|_{L^2_{\sqrt{\gamma'}}}^2.
\end{eqnarray*}
From above, we know that for $v = Q_{1,y(t)}$,
\begin{eqnarray*}
3 \left[ \langle \Lop \tilde{w}, \tilde{w} \rangle + \frac{21}{4}  \int \tilde{w}^2 dx \right]= 0 .
\end{eqnarray*}
This corresponds to $\tilde{w} = Q_{1,y(t)}^{\frac{5}{2}}$, which is the ground state or $\Lop_{1,y(t)}$.  Hence, $v=Q$ is the ground state of the quadratic form 
\begin{eqnarray*}
3 \left[ \langle \Lop \tilde{w}, \tilde{w} \rangle + \frac{21}{4} \int w^2 dx \right] .
\end{eqnarray*} 
From Lemmas \ref{lem:LopXYbounds}, our orthogonality condition
\begin{eqnarray*}
v \perp \tilde{Q}_{c(t),y(t)}
\end{eqnarray*}
is enough to guarantee there exists $\delta  > 0$ such that
\begin{eqnarray*}
\frac{d}{dt} I_\eta (v) + \|  v \|_{H^1_{\sqrt{\gamma'}}}^2 \le 0 
\end{eqnarray*}
provided $|c^2-1^2|+|\dot{y}-c^2| $ is small for all $t \ge 0 $, which follows from our assumptions on the initial perturbation.

The time dependent version of
\begin{eqnarray*}
\frac{d}{dt} \int  \gamma(x,t) (v^2 + v^2_x) dx
\end{eqnarray*}
is done in full generality in the analysis of \eqref{eq} 
in Section \ref{sec:linkdv} for the Airy equation. The terms which we have to
control are the same as for constant $c$ and $\dot y$, plus the terms coming 
from the right hand side. Those are easy to control. Namely, 
\[ 2 \left| \int \left[ (\alpha Q + \beta Q')\gamma  - \partial_x( \gamma  (\alpha Q' + \beta Q''))\right] v  dx   \right| \lesssim  
(|\dot c| +|\dot y-c^2|) \Vert |v|^2 \|_{L^2_{\sqrt{\gamma'}}} 
 \]
for $\gamma$ as in Section \ref{sec:linkdv}.
  \end{proof} 

Note, above we have always assumed the proper orthogonality conditions, but 
without them  we easily obtain the following estimate for solutions of the $v$ equation. 
\begin{equation} \label{nonorth} 
 \Vert v \Vert_{ L^\infty H^1  \cap L^2 \Hgamma^2 } \le C 
\left( \Vert v(0) \Vert_{H^1} + \sup_{t} | \langle v(.,t),Q'_{1, y(t)} \rangle | 
+ \Vert \langle v(.,t),\tilde Q_{1,y(t)}  \rangle \Vert_{L^2([0,\infty))} \right) .
\end{equation}

\section{Function Spaces and Projection Operators}
\label{sec:lin-constant}

In this section we construct the function  spaces for our
nonlinear analysis using properties of the linear evolution we have
studied in Sections \ref{sec:linkdv}-\ref{sec:virial}.
Based on the energy functional \eqref{eqn:energy} for the $v$-equation, it seems natural to look at \index{$X^1$}
\begin{eqnarray*}
v \in X^1 = L^\infty H^1  \cap L^2 \Hgamma^2  ,
\end{eqnarray*}
where $\gamma= \gamma_0(x)$ is as in \eqref{stat} and again by convention we set $L^p X$ to be the $L^p$ norm in time of the $X$ norm in space.

Then, as follows naturally from the equation, we define \index{$Y^1$}
\begin{eqnarray*}
Y^1 = L^1 H^1  +  L^2 \sqrt{\gamma'}L^2 .
\end{eqnarray*}

Generically, we define  \index{$X^{s}$}
\begin{eqnarray*}
X^{s} =  L^\infty H^{s} \cap L^2 \Hgamma^{s+1}  
\end{eqnarray*}
and \index{$Y^{s}$}
\[
Y^{s} = L^1 H^{s}  +  L^2  \gammaH^{s-1} ,
\]
where we note $(Y^s)^* = X^{-s}$.

\subsection{The Scale of  Energy Spaces}
\label{sec:est}

Let us study the $v$-equation
\begin{eqnarray}
\label{eqn:nat}
\left\{ \begin{array}{c}
(\p_t - \Lop \p_x) v = f_0 + \sqrt{\gamma'} f_1 = f , \\
v(0,x) = v_0 ,
\end{array} \right.
\end{eqnarray}
where $f_0 \in L^1 H^{s}$, $f_1 \in L^2 L^{s-1}$, $v_0 \in H^s$. We assume that 
the orthogonality conditions 
\begin{eqnarray}
\label{eqn:orth1}
v_0 \perp Q'  ,  v_0 \perp \tilde{Q} 
\end{eqnarray}
and
\begin{eqnarray}
\label{eqn:orth2}
\left( f_0 + \sqrt{\gamma'}  f_1 \right) \perp \tilde{Q} , \left( f_0 + \sqrt{\gamma'}  f_1 \right) \perp  Q' \ \ \  \text{ for all $t$}  
\end{eqnarray} 
hold. 

\begin{prop} \label{proplin0} 
There exists a unique solution $v \in X^s$ which satisfies 
\[ \Vert v \Vert_{X^s} \le c \left( \Vert v_0 \Vert_{H^s} + 
\Vert f_0 + \sqrt{\gamma'} f_1 \Vert_{Y^s} \right). 
\]
Moreover, $v(t)$ is orthogonal to $Q'$ and $ \tilde Q$. 
\end{prop} 

Note, Theorem \ref{thm:linmain} is an immediate consequence. 

\begin{proof} We begin by considering the case $s=1$. The previous section 
implies the estimate 
\[ \Vert v \Vert_{X^s} \le c \left( \Vert v_0 \Vert_{H^s} + 
\Vert f_0 \Vert_{L^1(H^s)}  \right). 
\]
if $f_1=0$  by a variation of constants argument. 
We retrace the steps  and its modifications needed for $f_1$. 
Using the multipliers from the energy inequalities, we need  the obvious 
estimates  
\[
\left| \int f_0 \gamma (v - v_{xx}) dx dt \right| +
 \left| \int f_0 \eta v dx\, dt
 \right|   + \left| \int f_0 \mathcal{L}^{-1}  v dx dt \right|
\le c  \Vert  v \Vert_{L^\infty  H^1} \Vert f_0 \Vert_{L^1H^1}  
\]
and, using Lemma \ref{lem:LopXYbounds}, 
\[
\left| \int \sqrt{\gamma'} f_1 (\gamma (v - v_{xx})dx dt \right|
 +\left|\int \sqrt{\gamma'}f_1 \eta v\,dx\,  dt \right|  
+\left| \int \sqrt{\gamma'}f_1 \Lop^{-1} v dx dt \right|
\le c \Vert f_1 \Vert_{L^2} \Vert v \Vert_{L^2 \Hgamma^2}. 
\]
It is not hard to see that $v(t)$ remains orthogonal to $Q'$ and $\tilde Q$ so that we can close the argument as in the previous section. 
We obtain the desired estimate for $s=1$: 
\[ \Vert v \Vert_{X^1} \le c \left( \Vert v_0 \Vert_{H^1} + \Vert f_0 + \sqrt{\gamma'} f_1 \Vert_{Y^1}\right).  \]
We denote  the solution operator for the inhomogeneous $v$-problem ($u$-problem) to be $S_{v}$ ($S_{u}$) \index{$S_{v}$} \index{$S_{u}$} and we write 
\begin{equation} \label{X1}
\Vert S_{v} f \Vert_{X^1} \le c \Vert f \Vert_{Y^1}.
\end{equation} 

The role of the two orthogonality conditions are different: The equation is invariant under the addition of a multiple of $Q$ to $v$, and orthogonality to 
$Q'$ is conserved. Orthogonality to $\tilde Q$ was needed for the virial 
identity of Martel-Merle, whereas orthogonality of $v$ and $Q'$ entered the control of the $H^{-1}$ norm by the Moore-Penrose inverse of $\Lop$. Without orthogonality one still obtains \eqref{nonorth}. 

Suppose now that $v$ satisfies 
\begin{equation}
\left\{ \begin{array}{c}
v_t - \Lop \p_x v =  f , \\
v(x,0) = v_0 .
\end{array} \right. 
\end{equation}

Let $\varepsilon$ be a small constant. We apply $(1+ \varepsilon^2 D^2)^{\frac{s-1}{2}}$ 
to both sides of the equation and denote $v^s = (1+\varepsilon^2 D^2)^{\frac{s-1}{2}} v$. 
It satisfies 
\[  
v^s_t - \Lop \p_x v^s =  
  (1+\varepsilon^2D^2)^{\frac{s-1}{2}}  f  
+   [(1+\varepsilon^2 D^2)^{\frac{s-1}{2}},4Q^3] \partial_x v  .
\]  
Hence, applying \eqref{nonorth} 
\[ 
\begin{split} 
\Vert v \Vert_{X^s}  \le & c_1   \Vert v^s \Vert_{X^1} 
\\ \le &  c_2 \left( \Vert   (1+\varepsilon^2D^2)^{\frac{s-1}{2}}  f  
\Vert_{Y^1} +   \Vert  [(1+\varepsilon^2 D^2)^{\frac{s-1}{2}},4Q^3] \partial_x v 
\Vert_{Y^1} \right. \\ 
& \left. + \sup_t | \langle v^s,  Q' \rangle | + 
\Vert  \langle v^s, \tilde Q \rangle \Vert_{L^2}\right). 
\end{split}  \]
and we turn to the commutator term. 

\begin{lem}\label{kernelb}  
Let $\phi\in C^\infty(\mathbb{R})$ satisfy $|\phi|+ |\phi'| \le  C e^{-|x|}$. Let 
$k(x,y)$ be the kernel of the operator 
\[  [(1+\varepsilon^2 D^2)^{\frac{s}{2}}, \phi] (1+\varepsilon^2  D^2)^{-\frac{s}{2}}.  \]
Then, 
\[ |k(x,y)| \le c 
\varepsilon |s| e^{-(|x|+|y|)/4 -  |x-y|/(4\varepsilon)}.
 \]
\end{lem} 
We postpone its proof.
 By Lemma \ref{kernelb} (with $\phi = 4 Q^3$ and $s-1$) and Schur's Lemma 
\begin{equation} \label{eqn:smalla}
\begin{split} 
\Vert  [(1+\varepsilon^2 D^2)^{\frac{s-1}{2}},4Q^3] \partial_x v 
\Vert_{Y^1} \le &   
\Vert (\gamma')^{-\frac12}  [(1+\varepsilon^2 D^2)^{\frac{s-1}{2}},4Q^3](1+\varepsilon^2 D^2)^{\frac{1-s}2}  \partial_x v_s 
\Vert_{L^2} 
\\ \le & c \varepsilon \Vert (\gamma')^{\frac12} \partial_x v^s \Vert_{L^2}, 
\end{split} 
\end{equation} 
and by Lemma \ref{lem:commutator}, after rescaling, as for the constant coefficient equation,     
\[
 \Vert (1+\varepsilon^2 D^2 )^{\frac{s-1}2} f \Vert_{Y^{1}}  \le c  \Vert  f \Vert_{Y^s}. 
\]

For all Schwartz functions,  
\[
 \Vert (1+|x|^2)^N \left[ (1+ \varepsilon^2 D^2)^{\frac{s}{2}} \phi- \phi \right] 
\Vert_{L^2} \le C \varepsilon .
\]
If $\langle v,\tilde Q \rangle = \langle v,Q' \rangle =0$, then
\begin{equation} \label{eqn:smallb}
 | \langle v^s,\tilde Q \rangle | = | \langle v,\tilde Q \rangle - \langle v^s,   \tilde Q- (1+\varepsilon^2 D^2)^{-\frac{s}{2}} \tilde Q \rangle |
\le C \varepsilon \Vert {\gamma'}^{1/2} v^s \Vert_{L^2} 
\end{equation} 
and
\begin{equation} 
\label{eqn:smallc} 
 | \langle v^s,Q' \rangle | =|   \langle v,Q' \rangle - \langle v^s,   Q'- (1+\varepsilon^2 D^2)^{-\frac{s}{2}} Q' \rangle |
\le C \varepsilon \Vert v^s \Vert_{L^2}. 
\end{equation} 
Suppose that
\[ \langle f, Q' \rangle = \langle f,\tilde Q \rangle = 0. \]
Then we obtain for all $s\in \mathbb{R}$ from  \eqref{eqn:smalla}, \eqref{eqn:smallb} and \eqref{eqn:smallc} 
\[ \Vert v^s \Vert_{X^s} \le c \left(\Vert f \Vert_{Y^s} +  \varepsilon \Vert v^s \Vert_{X^s} \right) 
\]
and hence 
 \begin{equation}\label{Hsv} 
 \Vert v \Vert_{X^s} \lesssim   \Vert f \Vert_{Y^s},   
\end{equation} 
which again implies for solutions $v$  to
\[ v_t - \Lop \partial_x v = P_{Q'} f ,\]
 given by the variaton of constants formula,  the bound
\[ \Vert \tilde P^* v \Vert_{X^s} \le C \Vert f \Vert_{Y^s} \]
or, equivalently (recall \eqref{projection}) 
 \begin{equation} \label{eqn:energy0}
\| \tilde P^* S_{v}  P^\perp_{Q'}\|_{Y^s \to X^s} \lesssim 1.
\end{equation}

Using space-time duality, we consider
\[
\left\{ \begin{array}{c}
(\p_t - \p_x \Lop) u = g , \\
u(0,x) = 0 .
\end{array} \right.
\]

The adjoint estimate to \eqref{eqn:energy0} is 
\begin{equation} 
\label{eqn:energy4}
\| P^\perp_{Q'}  S_{u} \tilde{P}  \|_{Y^{s} \to X^{s}} \lesssim 1.
\end{equation} 
\end{proof} 

\begin{proof} 
We turn to the proof of Lemma \ref{kernelb}. 

  Let $\hat \phi$ be the Fourier transform of $\phi$, which, because of the exponential 
decay extends to a holomorphic function $\hat \phi$ is the strip $\{ z : |\text{ Im } z| < 1\}$. Moreover there exists $C$ so that 
\[ \int |\hat \phi (\xi+ i \sigma) | d\xi \le C \qquad \text{ if } |\sigma| \le \frac12. \]
This estimate in turn implies exponential decay. Let $k(x,y)$ be the integral kernel of 
\[ [(1+\varepsilon^2 D^2)^{\frac{s}2},  \phi] (1+\varepsilon^2 D^2)^{-\frac{s}2}.  \] 
  We claim 
\begin{equation}\label{kernelest}  
| k(x,y) | \le c_N \varepsilon |s| e^{-\delta (|x|+|y|) }  e^{-\delta |x-y|/\varepsilon }  ,      
\end{equation} 
which implies Lemma \ref{kernelb}.  

The symplectic  Fourier transform 
\[ \hat k(\xi,\eta) = \frac1{2\pi} \int e^{-i\xi x+iy \eta} k(x,y) dx dy \]
satisfies 
\[ \hat k(\xi,\eta) = \left( \left( \frac{1+ \varepsilon^2 \xi^2}{1+\varepsilon^2 \eta^2} \right)^{\frac{s}2} -1 \right) \hat\phi(\xi-\eta). 
\]
We set $a= \varepsilon(\xi+\eta)/2$ and $b= (\xi-\eta)/2$. Then 
\[ \hat k(\xi,\eta) =  \hat g(a,b), \]
where 
\[  \hat g(a,b ) = \left( \left( \frac{1+ (\varepsilon b+a)^2}
{1+(\varepsilon b-a)^2} \right)^{\frac{s}2}  -1 \right) \hat \phi(2b)  
\]
and 
\[
\begin{split} 
  k(x,y) = & (2\pi)^{-1} \int e^{i(x\xi-y\eta)}
 \hat g(\varepsilon(\xi+\eta)/2,(\xi-\eta)/2)) d\xi d\eta \\
= &  2 \varepsilon (2\pi)^{-1} \int e^{i (\frac{x-y}{\varepsilon}  a +b(x+y))}  \hat g(a,b)    da \, db \\
=:& 2\varepsilon  g((x-y)/\varepsilon,x+y). 
\end{split} 
\]
The function $\hat g$ expands to a holomorphic function in $a$ to the strip 
$\{ z : |\text{ Im }  z| < 1/2\}$ if $\varepsilon|\text{ Im } b | < \frac12$.  
Clearly, 
\[ \frac{1+ (a+\varepsilon b)^2}{1+(a-\varepsilon b)^2} = 1+ \frac{4(\varepsilon b)^2}{1+(a-\varepsilon b)^2} + 4\varepsilon b  \frac{a-\varepsilon b}{1+(a-\varepsilon b)^2}, 
\] 
and hence we define the error term $h$ by the right hand side of 
\[  \left( \frac{1+ (\varepsilon b+a)^2}
{1+(\varepsilon b-a)^2} \right)^{\frac{s}2}  -1 
= 2s \varepsilon b \frac{ a-\varepsilon b}{ 1+ (\varepsilon b-a)^2} 
+   h(\varepsilon b, a).  
\]
It satisfies 
\[ |h(\varepsilon b, a) | \le  c s^2 \varepsilon^2 |b|^2 (1+ |\varepsilon b - a | )^{-2} , \]
if $|\varepsilon \text{ Im } b+a| \le \frac12$. 
Hence, 
\[ \left|\int e^{ i(av+bw) } h(\varepsilon b, a) \hat \phi(2b) da db\right| \le c s^2\varepsilon^2 e^{-\frac14 (|v|+|w|)} \]
by the extension of $a$ and $b$ to a suitable complex strip. 
The leading term contributing to $g$ can be calculated: 
\[ \begin{split} 
g_0(v,w) = & (2\pi)^{-1}  \int e^{i ( av +bw)}
\frac{a-\varepsilon b}{1+(a-\varepsilon b)^2}      b \hat \phi (2b)     da \, db
\\  = &  i \frac{v}{|v|} e^{-|v|}  \int e^{i (b(w + \varepsilon v) )} b \hat \phi(2b) db 
\\  = & \sqrt{\pi}2 i \frac{v}{|v|} e^{-|v|} \phi'( (w+ \varepsilon v)/2)).
\end{split} 
\]   
The leading term for $k$ is 
\[ k_0(x,y) = \sqrt{\pi}2 i \varepsilon  \frac{x-y}{|x-y|} e^{-|x-y|/\varepsilon}
\phi'(x). \]
This completes the proof.  
\end{proof} 

\subsection{$U$ and $V$ Space Estimates}
\label{sec:UV-subcrit}

In this section, we generalize and improve  Theorem \ref{thm:linmain}
using the $U^p$ and $V^p$  spaces as defined in \cite{HHK} and in the
 Appendix \ref{app:u2v2}.
For notational simplicity, let us define \index{$U^p_{KdV}$}
\index{$V^p_{KdV}$} \index{$U^p$!usage} \index{$V^p$!usage}
\begin{eqnarray*}
U^p = U^p_{KdV}, \ V^p = V^p_{KdV}.
\end{eqnarray*}

We begin with a number of estimates which we will use often in the sequel.

Let $c,y \in C^1$ satisfy  \eqref{eqn:cysmall}  
and let 
\[  \gamma (x,t) = \gamma_0 (x-y(t))  . \]
Then, 
\[\Vert  a Q'_{c(t), y(t)} + b \tilde Q_{c(t),y(t)} \Vert_{Y^0 }
\lesssim    \Vert a \Vert_{L^2+L^1} +  \Vert b \Vert_{L^2+L^1} ,
\] 
hence
\[ 
\Vert  P^\perp_{Q'} \tilde P f  \Vert_{DU^2 + L^2 \gammaH^{-1}} \lesssim  
\Vert f \Vert_{DU^2} 
+ \Vert \langle f,Q \rangle \Vert_{L^2+L^1} +  \Vert \langle f, Q' \rangle \Vert_{L^2+L^1}. 
\]
We consider 
\[ w_t + w_{xxx} = f, \qquad w(0)= u_0 . \]
Then, 
\[ \Vert w \Vert_{U^2} \lesssim \Vert u_0 \Vert_{L^2} + \Vert f \Vert_{DU^2} \]
and, since $U^2 \subset L^2 \Hgamma^1$, 
\[ \Vert \langle w, Q \rangle \Vert_{L^2} + \Vert \langle w, \tilde Q \rangle \Vert_{L^2} \lesssim  \Vert f \Vert_{DU^2}+ \Vert u_0 \Vert_{L^2}. 
\]
Hence, with $v = \tilde P  P^\perp w$, we have 
\[ \Vert v  \Vert_{L^2\Hgamma^1} 
\lesssim  \Vert f  \Vert_{DU^2} + \Vert u_0 \Vert_{L^2}. \] 
We calculate 
\[ 
\begin{split} 
(\partial_t + c^2 - \partial_x \Lop)  \left( \frac{\langle w, Q \rangle}{\langle
  Q, \tilde Q\rangle } \tilde Q + \frac{\langle w, Q' \rangle}{\langle
  Q', \tilde Q'\rangle } Q' \right) = 
& \frac{\dot c}{c} \frac{\langle w, Q \rangle}{\langle Q, \tilde
  Q\rangle }\tilde {\tilde Q} + (\dot y-c^2) \frac{\langle w, Q'
  \rangle}{\langle Q', \tilde Q'\rangle } Q''
 \\ &   + \left[ (\dot
y-c^2) \frac{\langle w, Q \rangle}{\langle Q, \tilde Q\rangle }+
\frac{\dot c}{c} \frac{\langle w, Q' \rangle}{\langle Q', \tilde Q'\rangle } \right] \tilde Q' \\
& -\tilde \alpha \tilde Q - \tilde \beta Q',
\end{split} 
\]
where $\tilde \alpha$ and $\tilde \beta$ are the time derivatives of
the coefficients of $\tilde Q$ and $Q'$ and 
\[\tilde{\tilde Q} = (x-y)\tilde Q'+ \frac23 \tilde Q. \] 
 Hence,  assuming that $u_0$ satisfies the orthogonality conditions, that $w$ and $v$ are as above, and with 
$g$ defined through the previous calculations, 
\begin{equation} \label{projv}  
 \partial_t  v + c^2 v_x - \partial_x \Lop_{c,y}  v = \tilde \alpha \tilde Q + \tilde \beta Q' + g + f ,
 \qquad   v(0) = u_0, 
\end{equation}  
where we collect the properties of $v$ and $g$ in the following
\begin{lem} \label{initialu} 
Assuming \eqref{eqn:cysmall}, we have $\langle v(t), Q \rangle = \langle v(t), Q' \rangle = 0$ and  
\begin{equation}
  \Vert v \Vert_{V^2 \cap L^2 \Hgamma^1} +   \Vert g \Vert_{L^2 \gammaH^{-1}} 
  \lesssim \Vert u_0 \Vert_{L^2} + \Vert f \Vert_{DU^2} + \Vert \langle f,Q \rangle \Vert_{L^2+L^1} + \Vert \langle f, Q' \rangle \Vert_{L^2+L^1}.  
\end{equation} 
\end{lem} 

\begin{proof} We claim that 
\[ \Vert \tilde \alpha \Vert_{L^1+L^2} + \Vert \tilde \beta \Vert_{L^1+L^2} 
\le c \left(  \Vert w_0 \Vert_{L^2} +  \Vert f \Vert_{DU^2} \right)  
\]
the proof of which we postpone. Assuming its validity we put the term
$4 \partial_x Q v$ in \eqref{projv} on the right hand side. We bound $\Vert v \Vert_{V^2}$ in terms of
$\Vert w_0 \Vert_{L^2}$ and the right hand side in $DV^2$. Since $DU^2
\subset DV^2$ and $L^2 \sqrt{\gamma'} H^{-1} \subset DV^2$ we can control all
terms on the right hand side.

 The only missing piece is the $L^2+L^1$ bound for $\alpha$ and $\beta$. There are two different arguments: either we can follow the calculation above and calculate 
$\tilde \alpha$ and $\tilde \beta$ 
above, or we can test by $Q$ and $Q'$ and use orthogonality 
to obtain the standard equations for $\alpha$ and $\beta$. We use the first
approach and recall the calculations after \eqref{tildeQ}.
Then 
\begin{equation} 
\begin{split} 
\frac{d}{dt} \langle w,Q \rangle = &  \langle f, Q \rangle + \dot y \langle w, Q'\rangle +  \frac{\dot c}{c} \langle w, \tilde Q \rangle , \\
\frac{d}{dt} \langle w,Q' \rangle = &  \langle f, Q' \rangle + \dot y \langle w, Q''\rangle +  \frac{\dot c}{c} \langle w, \tilde Q' \rangle.
\end{split} 
\end{equation} 

There is one more term entering the coefficient of $Q'$ coming from applying the linear operator to $\tilde Q$, which gives  
\[ -2 \frac{\langle w, Q\rangle}{\langle Q, \tilde Q \rangle}    c^2 Q'.  
\]
All these terms are easily controlled. 
\end{proof}

We return to  the analysis of the time dependent $v$-problem
\begin{equation}
\label{eqn:vtdf}
v_t + c^2 c_x - \Lop \p_x v = \alpha(t) Q_{c(t),y(t)} + \beta(t) Q_{c(t),y(t)}' + f  ,
\end{equation}
where
\begin{equation}
\label{eqn:alpha-subtdf}
\alpha (t) = - \frac{ \frac{\dot{c}}{c} \langle v, \tilde{\tilde{Q}}_{c(t),y(t)} \rangle + (\dot{y}-c^2) \langle v, \tilde{Q}_{c(t),y(t)}' \rangle }{\langle Q_{c(t),y(t)} , \tilde{Q}_{c(t),y(t)}  \rangle}  - \frac{\langle \tilde Q, f \rangle}{\langle Q,\tilde Q \rangle}
\end{equation}
and
\begin{equation}
\label{eqn:beta-subtdf}
\beta (t) = - \frac{ \frac{\dot{c}}{c} \langle v , \tilde{Q}_{c(t),y(t)}' \rangle + (\dot{y}-c^2) \langle v, \tilde{Q}_{c(t),y(t)}'' \rangle }{\langle Q_{c(t),y(t)}' , Q_{c(t),y(t)}' \rangle } -  \frac{\langle Q', f \rangle}{\langle Q', Q' \rangle}.
\end{equation}
with the initial data $v(x,0) = v_0(x)$ orthogonal to $\tilde Q$ and $Q'$.  Then also $v(t)$ satisfies these orthogonality conditions. 
We combine the arguments of the previous subsection with those of Proposition \ref{prop:tdenergy} and obtain

\begin{lem} \label{propvarcoefv} 
Suppose that \eqref{eqn:cysmall} holds. 
There exists a unique solution $v$ to \eqref{eqn:vtdf} and \eqref{eqn:alpha-subtdf} and \eqref{eqn:beta-subtdf} 
 which satisfies
\[ \Vert \langle D \rangle^s  v \Vert_{X^0\cap V^2} \le c \left( \Vert v_0 \Vert_{H^s} + 
\Vert \langle D \rangle^s  f \Vert_{Y^0 + DV^2} \right). 
\]
Moreover $v(t)$ is orthogonal to $Q'$ and $\tilde Q$.
 \end{lem}

\begin{proof} 
We begin with $s=0$.
We write $f= f_U + f_Y$ with $f_U \in DU^2$ and $f_Y \in Y^0$.  
Let $\tilde v$ be defined with $f=f_U$ as in Lemma \ref{initialu}. It satisfies 
\[
 \Vert \tilde v \Vert_{V^2 \cap X^0} \le C ( \Vert f_U \Vert_{DU^2} + \Vert u_0 \Vert_{L^2} ).
\] 
Let us take $v= \tilde v+ w$, where $w$ satisfies 
\[ 
w_t + w_{x} - \partial_x \Lop   w = \alpha \tilde Q + \beta Q' +
f_Y + g, \qquad w(0)=0, \] 
with $g$ as in Lemma \ref{initialu} and by Lemma \ref{propvarcoefv} 
\[ \Vert w \Vert_{X^0} \lesssim   \Vert f_Y \Vert_{Y^0} +  \Vert g \Vert_{Y^0} 
\lesssim  \Vert f \Vert_{DU^2+Y^0} + \Vert u_0 \Vert_{L^2}.
\]
We put the term $4 \partial_x (Q^3 w)$ to the right hand side, which we easily 
control in $Y^0$ as well as $\alpha$ and $\beta$ and we 
 arrive at 
 \begin{equation}\label{uV} 
\Vert v \Vert_{V^2 \cap X^0 } \le C \left(
     \| v_0 \|_{L^2} + \Vert f \Vert_{Y^0 + DU^2} + \Vert \langle f,Q
     \rangle \Vert_{L^2+L^1} + \Vert \langle f, Q' \rangle
     \Vert_{L^2+L^1}\right).
\end{equation}
The case of general $s$ follows by the same arguments as above.
\end{proof} 

Our main interest will be in similar estimates for the $u$ problem below.

We consider the $u$ equations 
\begin{equation} \label{eqn:td}
 u_t + c^2 u_x  - \partial_x (\Lop_{c,y} u) = \alpha \tilde Q  + \beta Q' + f, 
\end{equation}  
 with  initial data $u(0)= u_0 $ which satisfies  
$\langle u_0, Q\rangle = \langle u_0, Q' \rangle =0$, together with the modal equations
\begin{eqnarray}
  \label{eqn:alpha-sub-td}
\alpha (t) = - \frac{\frac{\dot c}{c} \langle u, \tilde{Q}\rangle + \langle f , Q \rangle }{\langle Q , \tilde{Q} \rangle}
\end{eqnarray}
and
\begin{eqnarray}
\label{eqn:beta-sub-td}
\beta (t) = - \frac{ (\dot y - c^2) \langle u, Q'' \rangle + \frac{\dot c}{c} \langle u, \tilde{Q}' \rangle + \langle u , \Lop Q_{xx} \rangle + \langle f , Q' \rangle }{\langle Q', Q' \rangle},
\end{eqnarray}
which again ensures the  orthogonality of $u(t)$ with $Q$ and  $Q'$. 

We obtain first the analog of Lemma \ref{propvarcoefv}. 
\begin{lem} \label{propvarcoefu} 
Suppose that \eqref{eqn:cysmall} holds. 
There exists a unique solution $u$ to \eqref{eqn:td}, \eqref{eqn:alpha-sub-td} and \eqref{eqn:beta-sub-td}   which satisfies
\[ \Vert u \Vert_{X^0\cap U^2} \le c \left( \Vert u_0 \Vert_{L^2} + 
\Vert f \Vert_{Y^0+DV^2} \right). 
\]
Moreover, $u(t)$ is orthogonal to $Q'$ and $Q$. 
\end{lem}

It is not difficult to construct solutions, however we are
interested in global estimates. Moreover we may restrict to a finite time interval and assume that all the data as well as $u$ are smooth and decay at infinity.

We set $v = \Lop u$. It satisfies the orthogonality conditions  
\[ \langle v, Q' \rangle = 0  =  \langle u , Q \rangle = \langle \Lop^{-1} v , Q \rangle =  \langle v, \tilde Q \rangle .\] 
 Moreover, $v$ satisfies 
\[ v_t + c^2 v_x - \Lop \partial_x v = -2c^2 \alpha Q + 12 Q^2( \frac{\dot c}c \tilde Q + (\dot y -c^2 ) Q') u  + \Lop f \]
and we may apply Lemma  \ref{propvarcoefv} with $s=-2$: 
\[ \Vert v \Vert_{X^{-2}} \lesssim  \Vert \Lop u(0) \Vert_{H^{-2}} +\Vert \Lop f \Vert_{Y^{-2}} + 
(|\frac{\dot c}c|+ |\dot y -c^2|) \Vert u \Vert_{L^2(H^{-3}_{\sqrt{\gamma'}})}. \]
We apply Lemma \ref{lem:LopXYbounds} several times to get 
\[ \Vert u \Vert_{X^0} \lesssim \Vert \Lop v \Vert_{X^{-2} } \lesssim  
\Vert u_0 \Vert_{L^2} +
\Vert f \Vert_{Y^0} +   \sup_t (|\frac{\dot c}c|+ |\dot y -c^2|) \Vert u \Vert_{X^{-1}}. \]
To complete the proof we observe that by \eqref{eqn:cysmall} we may subtract 
the last term on the right hand side from both sides to arrive at the desired estimate. 
The inclusion of $V^2$ and $DU^2$ works now exactly as for the $v$ equation. 

We collect the results for the case $s=0$, which is the only estimate we will need later on. 

\begin{prop} \label{proplinUV} 
Suppose that \eqref{eqn:cysmall} holds. 
There exists a unique solution $v$ to \eqref{eqn:vtdf}, \eqref{eqn:alpha-subtdf} and \eqref{eqn:beta-subtdf}  which satisfies
\[ \Vert v \Vert_{V^2\cap X^0} \le c \left( \Vert v_0 \Vert_{L^2} + 
\Vert f \Vert_{DU^2+ Y^0} \right). 
\]
Moreover, $v(t)$ is orthogonal to $Q'$ and $ \tilde Q$. 
Similarly there is a unique solution $u$ to \eqref{eqn:td}, \eqref{eqn:alpha-sub-td} and \eqref{eqn:beta-sub-td}
 which satisfies 
\[ \Vert u \Vert_{V^2\cap X^0} \le c \left( \Vert u_0 \Vert_{L^2} +  \Vert f \Vert_{DU^2+ Y^0} \right). 
\]
Moreover, $u(t)$ is orthogonal to $Q'$ and $Q$. 
\end{prop} 

\subsection{Littlewood-Paley Decomposition}
\label{sec:lp-spaces}

We consider functions $c$ and $y$ satisfying \eqref{eqn:cysmall}  
 We set $\lambda \in \Lambda_0 = 1.01^{\mathbb{N}} $ and
$P_\lambda$  the Littlewood-Paley decomposition with Fourier
multipliers supported in $\{ \xi: 1.01^{-1} \lambda \le |\xi| \le 1.01
\lambda \}$ if $\lambda>1$ and $\{ \xi: |\xi| \le 1 \}$ if $\lambda=0$. \index{$\Lambda_0$} Then, we denote
\begin{eqnarray*}
u_\lambda = P_\lambda u.
\end{eqnarray*}

The Besov spaces are defined as the set of all tempered distributions for which the norm 
\[ \Vert v \Vert_{B^{s,p}_q} =    \left\|   \lambda^s \Vert v_\lambda \Vert_{L^p} \right\|_{l^q(\Lambda_0)} 
\]
is finite. Here $s \in \mathbb{R}$ and $1\le p,q\le \infty$. 
Similarly we define the homogeneous spaces $\dot B^{s,p}_q$ with the summation over 
$\Lambda= 1.01^{\mathbb{Z}}$, where the frequency $\lambda=1$ plays no special role.  
There is an ambiguity about the meaning of $v_0$, which differs depending on whether we consider 
$B^{s,p}_q$ or the homogeneous space $\dot B^{s,p}_q$.

We define the spaces $X^s_{\infty}$ and $Y^s_{\infty}$ using the norms
\[ \Vert u \Vert_{X^s_\infty} = \sup_{\lambda \in \Lambda_0} \lambda^s \Vert u_\lambda  \Vert_{V^2 \cap X^0},  \quad 
 \Vert f \Vert_{Y^s_\infty} =  \sup_{\lambda \in \Lambda_0} 
\lambda^s   \Vert f_\lambda   \Vert_{DU^2+ Y^0}.
\]
The homogeneous spaces $\dot X^s_\infty$ and $\dot Y^s_\infty$ are defined in the same way as the 
homogeneous Besov space $\dot B^{s,p}_q$
with $\Lambda =1.01^{\mathbb{Z}}$, 
though with a slight modification for $s<0$ in the $Y$ spaces due to the $\rho$ multiplier. \index{$\Lambda$} Namely, we take
\begin{equation}
\label{eqn:homXnorm}
\begin{split} 
\Vert u \Vert_{\dot X^s_\infty} = &  \sup_{\lambda \in \Lambda} \left( \lambda^s \Vert
u_\lambda  \Vert_{V^2 \cap  X^0} \right),  \\
\Vert F \Vert_{\dot Y^s_\infty} = & \inf_{F=f+ g } \Big( \sup_{\lambda \in \Lambda}
\lambda^s \Vert f_\lambda   \Vert_{DU^2} + \sup_{\lambda \in \Lambda_0} 
\lambda^s \Vert g_\lambda \Vert_{Y^0} \Big),
\end{split} 
\end{equation}
where there is a slight abuse of notation as the operators in $f_0$ and $g_0$ 
are taking on two different meanings:  The homogeneous projection for
$f_0$ and the inhomogeneous projection for $g_0$.  

We study
\begin{equation} \label{basic} 
  \left\{ \begin{array}{c}
u_t + u_{x} + \p_x \Lop u  = \alpha \tilde Q + \beta Q_x + f + 
\p_x (\rho   g) ,  \\
u(x,0) = 0,
\end{array} \right.
\end{equation}
where $\alpha$ is given by \eqref{eqn:alpha-sub-td} and $\beta$ by 
\eqref{eqn:beta-sub-td}. As a first step we obtain a weighted $L^2$ bound for $u$ 
in \eqref{eqn:ul2} below.

Let  $f= f^++f^-$ and $g= g^++g^-$ be a decomposition into high ($|\xi| > 1$) and low ($|\xi| \leq 1$) frequencies.  We define \index{$f^\pm$} \index{$g^\pm$}
\begin{eqnarray*}
\left\{ \begin{array}{c}
v_t + c^2 v_x - \Lop v_{x}   =  \alpha_+ Q + \beta_+ Q' + (\p_x^{-1} f^+  +  \rho  g^+) ,    \\
v(x,0) = 0,
\end{array} \right.
\end{eqnarray*}
where 
\[ \langle Q, \tilde Q \rangle \alpha_+= (c^2-\dot y) \langle v, \tilde
Q'\rangle - \frac{\dot c}{c} \langle v, \tilde{\tilde Q}\rangle - \langle
\p_x^{-1} f^++ \rho g^+ , \tilde Q\rangle ,
  \]
\[\langle Q', Q' \rangle \beta_+ = (c^2-\dot y) \langle v, Q'' \rangle -
\frac{\dot c}{c} \langle v, \tilde Q' \rangle + \langle   f^++ \p_x (\rho g^+ ),  Q\rangle 
 \]
 ensure $\tilde P^* P^\perp_{Q'} v = 0$. 
Then by Proposition \eqref{proplinUV}
\[ 
\Vert v \Vert_{X^0}\lesssim 
\Vert  \p_x^{-1} f^+  +  \rho g^+ \Vert_{DU^2 + Y^0}  \lesssim  
\Vert F^+ \Vert_{\dot Y^s_{\infty} } \lesssim  \Vert F \Vert_{\dot Y^s_{\infty} }   ,
\]
where the second inequality holds for all $s > -1$.

As a simple consequence, we obtain 
\[ 
\Vert P_{Q'} \p_x v \Vert_{L^2 L^2_{\rho}} \lesssim \Vert F \Vert_{\dot Y^s}
\]
and  compute similar to arguments above 
\begin{eqnarray*}
(\partial_t - c^2 \partial_x + \partial_x \Lop ) (P_{Q'} \p_x v) &=&
\left(\alpha_+ - \frac{d}{dt} \frac{\langle v, Q'' \rangle}{\langle Q' ,
  Q'\rangle}\right) Q' + \left(\beta_+ + \frac{\langle v, Q'' \rangle}{\langle Q' ,
  Q' \rangle} (\dot y- c^2)\right) Q'' \\
&& + \frac{\langle v, Q'' \rangle}{\langle
  Q' , Q' \rangle} \frac{\dot c}{c} \tilde Q' + f_+ + \partial_x
(\rho g_+).
\end{eqnarray*}

We make the ansatz $u = P_{Q'} \p_x v + u_-$ and observe that $\langle \p_x v, Q \rangle = 0$ by construction. Then, 
\begin{eqnarray*}
\partial_t u_- + c^2 \partial_x u_- - \partial_x \Lop u_- &=& \alpha \tilde Q
+ \beta Q' + f_- + \partial_x ( \rho g_- ) \\
&&+ \left(\frac{\dot c}{c} \frac{\langle v, \tilde Q' \rangle}{\langle
    Q',Q'\rangle}  -  \frac{\langle f^++\partial_x (\rho g^+), Q \rangle}{\langle Q' , Q' \rangle} \right) Q'' 
 -   \frac{\langle v, Q'' \rangle}{\langle Q' , Q' \rangle}
 \frac{\dot c}{c}    \tilde Q'  ,
\end{eqnarray*}
where $\alpha$ and $\beta$ ensure orthogonality. Later we will need the 
obvious identity (integrate by parts in the second term)
\[ 
\left\langle f_-+\partial_x (\rho g_-) - \left(\frac{\langle f^++\partial_x (\rho g^+), Q
  \rangle}{\langle Q' , Q' \rangle} \right) Q'' , Q \right\rangle 
= \langle F,Q \rangle .
\] 

Then, $u= \p_x P_{Q'}^\perp v + u_- $ and hence with $F^+ = f^+ + \partial_x \rho
g^+$ we have
\begin{equation}\label{eqn:ul2}  \Vert u \Vert_{L^2 L^2_\rho} \lesssim   \Vert F \Vert_{Y^s}
+ \Vert \langle F,Q\rangle \Vert_{L^2+L^1} +
\Vert \langle F^+, Q \rangle \Vert_{L^2+L^1}. 
\end{equation} 
By \eqref{eqn:alpha-sub-td} we see
\begin{equation} \label{prep1} 
 \Vert \alpha \Vert_{L^1}  \lesssim  
\Vert \dot c \Vert_{L^2\cap L^\infty} (\Vert F \Vert_{Y^s} + \Vert \langle
F,Q\rangle \Vert_{L^2+L^1} ) + \Vert \langle F, Q  \rangle \Vert_{L^1} 
\end{equation} 
and, using \eqref{eqn:cysmall}  
\begin{equation} \label{prep1a} 
 \Vert \alpha \Vert_{L^2}  \lesssim  \Vert F \Vert_{Y^s} + \Vert \langle
F,Q\rangle \Vert_{L^2} 
\end{equation} 
and  by \eqref{eqn:beta-sub-td}
\begin{equation} \label{prep2} 
 \Vert \beta \Vert_{L^2}  \lesssim \Vert F \Vert_{Y^s} + \Vert \langle F, Q' \rangle \Vert_{L^2}.     
\end{equation}

We turn to the frequency localized equation
\begin{eqnarray*}
\left\{ \begin{array}{c}
(u_\lambda)_t  + (u_\lambda)_{xxx}  = - P_\lambda \p_x (4 Q^3 u)  + \alpha P_\lambda \tilde Q + \beta P_\lambda Q_x + P_\lambda f + P_\lambda \p_x (\rho  g) , \\
u^\lambda(x,0) = 0.
\end{array} \right.
\end{eqnarray*}
Observe that by using first the boundednessr of Fourier multipliers 
 on $U^2$, $DV^2$ and the dual of the embedding $U^2 \subset L^2 \Hrho^1$, then
\[ \Vert P_\lambda \p_x (4 Q^3 u) \Vert_{DV^2} \lesssim \lambda 
\Vert Q^3 u \Vert_{DV^2}  \lesssim \lambda \Vert Q^2 u \Vert_{L^2 L^2}  \lesssim \lambda 
\Big( \Vert F \Vert_{\dot Y^s_\infty} + \Vert \langle F,Q \rangle \Vert_{L^2+L^1} \Big) . \]
If $\lambda >1$, then by Lemma \ref{lem:commutator}
\[  \Vert [P_\lambda \p_x , Q^3] u \Vert_{L^2 L^2_\rho} \lesssim  \Vert u
 \Vert_{L^2 L^2_\rho}. 
\]
Repeating these estimates for the term containing $g$ and
using the estimates of the previous section we obtain for $\lambda \le 1$,
\[ \Vert u_\lambda \Vert_{V^2 \cap L^2 \Hrho^1} 
\lesssim \Vert P_\lambda f \Vert_{DU^2}  + \lambda \Big(\Vert F \Vert_{Y^s} + 
\Vert \langle F,Q \rangle \Vert_{L^2+L^1}+ \Vert \langle F^+, Q \rangle \Vert_{L^2+L^1}  \Big) + \lambda^{\frac12} \Vert \alpha \Vert_{L^1} ,
\]
since 
\[ \Vert \alpha \tilde Q \Vert_{L^1 \dot{B}^{-\frac12,2}_\infty}  \lesssim 
 \Vert \alpha \Vert_{L^1}
\] 
and, for $\lambda >1$, 
\[ \Vert u_\lambda \Vert_{V^2 \cap L^2 \Hrho^1} 
\lesssim  \Vert f_\lambda \Vert_{DU^2}+ \Vert g_\lambda \Vert_{L^2} 
 + \Vert F \Vert_{Y^s} + \Vert \langle F,Q \rangle \Vert_{L^2+L^1} + \Vert \langle F^+, Q \rangle \Vert_{L^2+L^1} 
+  \Vert \langle F,Q' \rangle \Vert_{L^2+L^1}
.
\]

As a result, we arrive at the following key
\begin{prop}\label{littlewood}
Suppose \eqref{eqn:cysmall} holds for some small $\delta$,
that $-1/2<s<0$, $F \in Y^s$  and
\begin{eqnarray*}
\left\{ \begin{array}{c} 
u_t + u_{xxx} + 4 \p_x (Q^3 u)  = \alpha \tilde Q + \beta Q_x + F , \\
u(x,0) = 0,
\end{array} \right.
\end{eqnarray*} 
where $\alpha$ and $\beta$ are defined in \eqref{eqn:alpha-sub-td}
and \eqref{eqn:beta-sub-td}.
Then,
\[ \Vert u \Vert_{\dot X^s_\infty} \lesssim 
\Vert F \Vert_{\dot Y^s_\infty}  + \Vert \langle F,Q \rangle \Vert_{L^1} 
+ \Vert \langle F, Q' \rangle \Vert_{L^2+L^1} + \Vert \langle F,Q^+ \rangle \Vert_{L^2+L^1}. \]
\end{prop}
This result will play a large role in the nonlinear analysis required
to prove asymptotic stability.

For future use, we denote by $L^p_I$, $\dot X^{-\frac16}_{\infty,I}$
etc the function spaces on the space time set $I \times \mathbb{R}$,
and specifically we set $L^p_T$, $\dot X^{-\frac16}_{\infty,T}$ for $I
= (0,T)$. \index{$L^p_T$} \index{$\dot X^{-\frac16}_{\infty,T}$}
All previous constructions carry over to finite time intervals.

\section{Local Well-posedness for the Quartic KdV Equation}
\label{sec:wellposedness}

In this section we study local well-posedness for the quartic generalized KdV equation
\begin{eqnarray}
\label{eqn:kdv4}
\left\{ \begin{array}{c}
\p_t \psi - \p_{xxx} \psi - (\psi^4)_x = 0 , \\
\psi (0,x) = \psi_0 (x). 
\end{array} \right.
\end{eqnarray}
Let $v$ be the solution to the Airy equation with the same initial data. 
\begin{equation} \label{linhom}
\left\{ \begin{array}{c}
v_t + v_{xxx} = 0 , \\
v(0,x) = \psi_0 (x),
\end{array} \right.
\end{equation}
The main local wellposedness is the contents of the next 
\begin{thm}
\label{thm:lwp}
Let $r_0>0$. There exist $\epsilon_0, \delta_0 >0$ such that, if $0<T \le \infty$,
\begin{equation} \label{initial}
\Vert \psi_0  \Vert_{\dot{B}^{-\frac16,2}_\infty} \le r_0
\end{equation} \index{$r_0$}
and
\begin{equation}
\label{smallness}
\sup_\lambda \Vert v_\lambda  \Vert_{L^6 ([0,T],\RR)}\le \delta_0 ,
\end{equation}
then there is a unique solution $\psi = v + w$ with
\[ \Vert w \Vert_{\dot X^{-\frac16}_{\infty,T} } \le \epsilon_0 . \]
Moreover, the function $w$ (and hence $\psi$)  depends analytically on the initial data.
\end{thm}

By the Strichartz estimates for linear KdV (see also \eqref{eqn:linest2-wp} and 
\eqref{eqn:bilinest-wp} below), given $v$ as in \eqref{linhom} we have
\[ \sup_\lambda \Vert  v_\lambda \Vert_{L^6} \le \kappa_0  \Vert v \Vert_{\dot
 X^{-\frac16}_{\infty,T}},  \] \index{$\kappa_0$}
and by the definition of the spaces 
\[ \Vert v \Vert_{\dot X^{-\frac16}_{\infty, T}} \le \kappa_1 \left( \Vert \psi (0) \Vert _{\dot
 B^{-\frac16 ,2}_\infty}+ \Vert \p_t v + \p_{xxx} v \Vert_{\dot
 Y^{-\frac16}_{\infty,T}}\right).  \] \index{$\kappa_1$}
Hence, we obtain global existence from Theorem \ref{thm:lwp} 
for \eqref{eqn:kdv4} if
\[ \Vert \psi \Vert_{\dot B^{-\frac16,2}_{\infty}} 
\le  \min\{1, \frac{\delta(1)}{(\kappa_0 \kappa_1)} \}  , \]
where $\delta(1)$ is the $\delta$ (which depends on $r_0$) evaluated at $r_0=1$.  \index{$\delta(r_0)$}

In any case, if condition \eqref{smallness} is satisfied for  $T= \infty$,
then, since  $\psi_\lambda \in V^2$, $e^{t\p_{xxx}} \psi_\lambda  $ is of
bounded 2-variation with values in $L^2$ (see Appendix \ref{app:u2v2}), and hence it has a
limit in $L^2$ as $t\to \infty$ .  This
implies that $\sum_\lambda \lim_{t\to \infty} e^{t\p_{xxx}} \psi_\lambda  =: S(\psi_0) $ exists and is the
scattering state. If in addition $\psi_0$ is in the closure of
$C^\infty_0$ in $\dot B^{-\frac16,2}_\infty$, then we may exchange the summation and
limit.

Under the same assumptions we can solve the
initial value problem with initial data $\psi_0 (T) =  e^{-T\p_{xxx}} \psi_0$, which,
by an easy limit as $T \to \infty$, gives the inverse of the map $S$. We will later
see similar constructions for perturbation of the soliton.

It is not hard to see that if $\psi_0$ is in the closure of $C^\infty_0$ in
$\dot B^{-\frac16 ,2}_\infty$, then we can  achieve  condition \eqref{smallness} by
choosing $T$ small. This implies local existence with smooth dependence on
initial data.  Moreover, since we obtain smooth dependence on the initial data, if we have any global solution $\psi (t)$ in the closure of $C^\infty_0$ and perturb the initial data by an amount
$\varepsilon$, we obtain a solution at least with a life span
\[ T = - c \ln \varepsilon \]
by easy perturbation arguments. In particular, if the initial datum lies in an
$\varepsilon$ neighborhood of a soliton, then the solution exists at least
until time $\sim |\ln  \varepsilon| $ and remains in a small neighborhood until
that time.

Before turning to the proof we remark that in this section we work with 
the weaker norms
\[ \Vert u \Vert_{\dot X^{-\frac16}} = \sup_\lambda \lambda^{-\frac16} \Vert u_\lambda \Vert_{V^2} \]
and
\[ \Vert f \Vert_{\dot Y^{-\frac16}} = \sup_\lambda \lambda^{-\frac16} \Vert f_\lambda \Vert_{DU^2}. \]
On the other hand, since the results remain trivially true for the original definition of the spaces we keep the notation.

\begin{proof}

First, we recall some estimates for $u \in U^2_{KdV}$.  Let $m(\xi,\xi_1)=m(\xi,\xi-\xi_1)$. Then 
\begin{equation}
\label{eqn:linest2-wp}
\| u \|_{L^6_t L^6_x}  \lesssim  \| |D|^{-\frac{1}{6}} u(0) \|_{L^2} \ (\text{$L^6$ Strichartz Estimate}) , 
\end{equation} 
\begin{equation} 
\label{eqn:bilinest-wp}
\begin{split} 
\Vert \int_{\mathbb{R}} 
m(\xi,\xi_1)  |\xi_1^2 - (\xi-\xi_1)^2|^{\frac{1}{2}} \hat{u}_1(\xi_1) 
\hat{u}_2(\xi-\xi_1) d\xi_1 \Vert_{L^2(\mathbb{R}^2)} & \\ & \hspace{-5cm}  \lesssim  
 \sup \frac{|m(\xi, \xi_1)|^2}{|\xi_1^2-(\xi-\xi_1)^2|^{\frac12}}
\|  u_1 (0) \|_{L^2} \| u_2 (0) \|_{L^2} \ (\text{Bilinear Estimate}).
\end{split} 
\end{equation}

The bilinear estimate is a variant of standard estimates as in Gr\"unrock \cite{Gr}. 
The most important choice is $m = |\xi_1^2-(\xi-\xi_1)^2|^{\frac12}$. 

Let $m(\xi, \xi_1)$ be a function which satisfies $m(\xi,\xi-\xi_1)=m(\xi,\xi_1)$. Then, 
\[
\begin{split} 
  \Vert \int m(\xi,\xi_1) e^{it(\xi_1^3+(\xi-\xi_1)^3)} d\xi_1
  \Vert_{L^2}^2 \hspace{-3cm} & \\   = &  \int m(\xi,
  \xi_1)m(\xi,\eta_1)e^{3it(\xi_1^2-\xi\xi_1-\eta_1^2+\xi\eta_1)} \hat
  u(\xi_1)\hat u(\xi_1) \overline{ \hat u(\xi-\eta_1) \hat
    u(\eta_1) } dt d\xi_1 d\eta_1 d\xi
  \\
  = & \int \frac{|m(\xi,\eta_1)|^2}{|\eta_1^2-(\xi-\eta_1)^2|}
  |u(\eta_1)|^2 |u(\xi-\eta_1)|^2 d\xi d\eta_1 
\\ \le & \sup
  \frac{|m(\xi, \xi_1)|^2}{|\xi_1^2-(\xi-\xi_1)^2|} \Vert u_1
  \Vert_{L^2}^2\Vert u_2 \Vert_{L^2}^2
  \end{split} 
\]  
since 
\[ \phi(\xi_1)= \xi \xi_1^2 -\xi^2  \xi_1 - \xi \eta_1^2 + \xi^2 \eta_1 =
\xi (\xi_1-\eta_1)(\xi_1+\eta_1-\xi)
\]
 vanishes if $\xi_1= \eta_1$ or $\xi_1 = \xi - \eta_1$ and 
\[ \phi'(\eta_1)= \xi(2\eta_1-\xi)= \xi(\eta_1 -(\xi-\eta_1)), \qquad 
\phi'(\xi-\eta_1) = \xi((\xi-\eta_1) - \xi_1).
 \]

These results immediately imply (see Appendix \ref{app:u2v2} for more information) for
$\lambda \ge 1.1 \mu$ the estimates
\begin{equation}\label{strichartzU} 
\Vert u_\lambda \Vert_{L^6_T} \lesssim \lambda^{-\frac16} \Vert u_\lambda \Vert_{U^2_T},
\end{equation}
\begin{equation}\label{bilinear1} 
\Vert u_\lambda  u_\mu \Vert_{L^2_T}  \lesssim  \lambda^{-1} \|  u_\lambda \|_{U^2_T} \| u_\mu \|_{U^2_T}.
\end{equation}
By interpolating the bilinear estimate and the Strichartz estimate, if $2\le p <3$,  
\begin{equation}\label{bilinear1b} 
\Vert u_\lambda  u_\mu \Vert_{L^p_T}  \lesssim  \lambda^{-1}( \mu^{-\frac16} \lambda^{\frac56} )^{ \frac{3p-6}p}  \|  u_\lambda \|_{U^2_T} \| u_\mu \|_{U^2_T}
\end{equation}
and, if $\rho << \mu \sim \lambda$,  
\begin{equation}\label{bilinear2} 
\Vert (u_\lambda  u_\mu)_\rho  \Vert_{L^2_T}  \lesssim  \lambda^{-\frac12}\rho^{-\frac12} 
 \|  u_\lambda \|_{U^2_T} \| u_\mu \|_{U^2_T}.
\end{equation}
Interpolating once again, we have
\begin{equation}\label{bilinear2b} 
\Vert (u_\lambda  u_\mu)_\rho  \Vert_{L^{p}_T}  \lesssim   
 \lambda^{-\frac12}\rho^{-\frac12}(\lambda^{\frac16} \rho^{\frac12})^{\frac{3p-6}p}   
\|  u_\lambda \|_{U^{p}_T} \| u_\mu \|_{U^{p}_T}.
\end{equation}

We proceed with a standard fixed point argument, which requires bounds on the nonlinearity.  The solution $\psi = v + w$ is constructed by studying
\begin{eqnarray}
\label{eqn:weq}
\left\{ \begin{array}{c} 
w_t + w_{xxx} + (v+w)^4_x = 0 , \\
w(0) = 0,
\end{array} \right.
\end{eqnarray}
where again
\begin{eqnarray*}
\left\{ \begin{array}{c}
v_t+ v_{xxx}=0, \\
v(0)= \psi_0.
\end{array} \right.
\end{eqnarray*}
Then, the key estimate is contained in the following
\begin{lem} 
\label{lem:wp}
There exists $r>0$ independent of $T$ such that given $v_k \in \dot X^{-\frac16}_{\infty,T}$ for $k = 1,2,3,4$ we have
\begin{equation} 
\label{eqn:em1-wp} 
\Vert \p_x (v_1 v_2 v_3 v_4) \Vert_{\dot Y^{-\frac16}_{\infty,T}} \le r  \prod_{k = 1}^4 \Vert v_k \Vert_{\dot X^{-\frac16}_{\infty,T}} ,
\end{equation} 
and, with $v$, $w$ defined by \eqref{linhom} and \eqref{eqn:weq} respectively
\begin{equation}
\label{eqn:em2-wp} 
\Vert \p_x (v^3 w) \Vert_{\dot Y^{-\frac16}_{\infty,T}} \le r 
\sup_\lambda \Vert v_\lambda \Vert_{L^6} \Vert \psi_0  \Vert_{B^{-\frac16,2}_\infty}^2
\Vert w \Vert_{X^{-\frac16}_{\infty,T}}.
\end{equation} 
\end{lem}

We apply these estimates to $ v^4 + 4 v^3 w + 6v^2 w^2 + 4 v w^3 + w^4$.
Either we may choose to estimate one factor $v$ in $L^6$ or the dependence on $w$ is at least quadratic. Suppose that $\Vert w \Vert_{\dot X^{-\frac16}_{\infty,T}} \le \mu$.  We obtain
\[ \Vert \p_x (v+w)^4 \Vert_{\dot Y^{-\frac16}_{\infty,T}}
\le 6 r (\kappa_1^3 \delta r_0^3 + \kappa_1^2 \delta \mu r_0^2 
+  \kappa_1^2 r_0^2 \mu^2 +  \kappa_1 r_0 \mu^3 + \mu^4).
\]
If $  \mu \le  \kappa_1 r_0 $,
then the right hand side is bounded by 
\[ 20 r  ( \kappa_1^3 \delta r_0^3 +  \kappa_1^2  \mu^2 r_0^2).   \]
Suppose that
\begin{equation} 
\mu \le  \min\left\{\kappa_1 r_0, \frac{1}{40 r \kappa_1^2 r_0^2}\right\}
\end{equation} 
and 
\begin{equation}  
\delta \le \frac{\mu}{40 r\kappa_1^4 r_0^3}.
\end{equation}
If $w$ solves
\begin{eqnarray*}
\left\{ \begin{array}{c}
w_t + w_{xxx} + (v+ W)^4_x = 0 , \\
w(0) = 0
\end{array} \right.
\end{eqnarray*}
and $\Vert W \Vert_{\dot X^{-\frac16}_{\infty,T}} \le \mu$, then $w$ exists and satisfies $\Vert w \Vert_{\dot X^{-\frac16}_{\infty,T}} \le \mu$.

Standard arguments then allow one to construct a unique solution satisfying the contraction assumption, possibly after decreasing  $\mu$ by an absolute multiplicative factor.
\end{proof}

It remains to prove Lemma \ref{lem:wp}. By duality, it suffices to verify
\[ 
\lambda \left|\int v_1 v_2 v_3 v_4 u_\lambda  \, dx \, dt \right|\le
C  \lambda^{\frac16} \Vert u_\lambda  \Vert_{V^2}    
\prod_{k = 1}^4 \Vert v_k \Vert_{\dot X^{-\frac16}_{\infty,T}} 
\] 
and 
\[ \lambda \left|\int v^3 w u_\lambda  \, dx \, dt\right| \le
C  \lambda^{\frac16} \Vert u_\lambda  \Vert_{V^2} 
  \sup_{\mu} \Vert v_\mu \Vert_{L^6} 
\left( \sup_{\mu}\mu^{-\frac16}  \Vert v_\mu  \Vert_{U^2}\right)^2 
\Vert w \Vert_{\dot X^{-\frac16}_{\infty,T}}  ,
\]
where $u_\lambda \in V^2$ is frequency localized at frequency $\lambda$.

By summation, the statement of the lemma holds provided we can prove the following
\begin{lem}\label{quartic1} 
We have for $\lambda_1 \le \lambda_2 \sim \lambda_3 \sim \lambda_5 \sim \lambda_5$ and $\varepsilon>0$
\begin{equation} 
\label{eqn:multilin1-wp}
\lambda_{5} \int v_{1,\lambda_1} v_{2,\lambda_2} v_{3,\lambda_3} v_{4,\lambda_4} v_{5,\lambda_5} dx dt
 \lesssim
\lambda_{min}^{-\frac13} \lambda^{-\frac16} \left( \frac{\lambda_{max}}{\lambda_{min}} \right)^{\varepsilon} 
   \prod_{k=1}^5 \Vert v_{k,\lambda_k} \Vert_{V^2} 
\end{equation}
and
\begin{equation} 
\label{eqn:multilin2-wp}
\begin{split} 
\lambda_{max} \int v_{\lambda_1} v_{\lambda_2} v_{\lambda_3} u_{\lambda_4} w_{\lambda_5} dx dt
 \lesssim &
 \left( \frac{ \lambda_{min} }{\lambda_{max}} \right)^{\frac16}
 \left( \frac{\lambda_{max}}{\lambda_{min}} \right)^{\varepsilon} 
\\ & \times \sup_{\mu} \Vert v_\mu \Vert_{L^6}
\left( \sup_{\mu}\mu^{-\frac16}  \Vert v_\mu  \Vert_{U^2}\right)^2
\Vert u_{\lambda_4}  \Vert_{V^2} \Vert w_{\lambda_5}   \Vert_{V^2}  ,
\end{split}
\end{equation}
where $\lambda_{max}$ respectively $\lambda_{min}$ is the maximal or the
minimal $\lambda_j$ and $\lambda_{med}^3$ is the product of the second,
third and fourth largest $\lambda_j$.
\end{lem}

\begin{proof}

We claim that
\begin{equation} 
\label{qu2}
\left| \int v_{1,\lambda_1} v_{2,\lambda_2} 
v_{3,\lambda_3} v_{4,\lambda_4} v_{5,\lambda_5}
\, dx\,  dt \right| 
\le C \lambda_{max}^{-1} \Vert  v_{1,\lambda_1} \Vert_{U^2} \Vert v_{2,\lambda_2} \Vert_{U^2} 
\Vert v_{3,\lambda_3}  \Vert_{L^6} \Vert v_{4,\lambda_4}  \Vert_{L^6} \Vert v_{5,\lambda_5}  \Vert_{L^6}
\end{equation}
provided
\begin{equation}
|\lambda_1 -\lambda_2| \ge \frac1{10} \lambda_{max}
  \label{lambda} .
\end{equation}
This estimate is a consequence of
H\"older's inequality and the bilinear estimate
\eqref{eqn:bilinest-wp}. We recall that
\[
\Vert v_\lambda \Vert_{L^6} \lesssim  \lambda^{-\frac16} \Vert v_\lambda
\Vert_{U^6} \lesssim \lambda^{-\frac16} \Vert v_\lambda
\Vert_{V^2}.
\]
In order to obtain a nontrivial integral there have to
be elements in the support of the Fourier transforms which add up to
zero. Unless there is at least one pair of $(\lambda_j,\lambda_k)$
satisfying \eqref{lambda}, the integral is zero.  Hence, we would obtain
\eqref{eqn:multilin1-wp} if we were allowed to replace the $V^2$ norms
there by $U^2$ norms for the first two factors. Observe that we may reorganize 
the factors as we wish.

Let us assume $\lambda_1 \lesssim \lambda_2 \lesssim \lambda_3 \lesssim 
\lambda_4 \lesssim \lambda_5 $. We consider first the case when 
$\lambda_4 \le 1.05 \lambda_1$. Then, if there are elements in the support 
of the truncations on the Fourier side adding up to zero - otherwise the integral vanishes -  either 
\[ 0.8 \frac{\lambda_5}{4} \le \lambda_1  \le \lambda_4 \le 1.1 \lambda_1 \le 1.2 \frac{\lambda_5}{4} \] 
or 
\[  0.6 \frac{\lambda_5}{2} \le \lambda_1  \le \lambda_4 \le 1.1 \lambda_1 \le 1.4 \frac{\lambda_5}{2}. \] 

In this case we can replace the $U^2$ norms by $V^2$ norms as follows.  
 We decompose into low and high modulation
\[ v_{j,\lambda_j} = v^l_{j,\lambda_j} + v^h_{j,\lambda_j}, \] 
where $v^l_{j} $ is defined by the Fourier multiplier projecting to $|\tau -
\xi^3| \le \lambda_5^3/1000 $. Then, we have
\[ \Vert v^l_{j,\lambda_j} \Vert_{V^2} + \Vert v^h_{j,\lambda_j} \Vert_{V^2} 
\le \Vert v_{j,\lambda_j} \Vert_{V^2} \]
and 
\[ \Vert v^h_{j,\lambda_j} \Vert_{L^2} \lesssim \lambda_5^{-\frac32} \Vert v_{j,\lambda_j} \Vert_{V^2}. \]
We refer to Appendix A and \cite{HHK} for more information.

We expand the product. The integral over the product of the five
$v^l_{j,\lambda}$ vanishes because of the support of the Fourier
transforms. Hence at least one term has high modulation. We estimate
it in $L^2$, put another term into $L^\infty$ and the others into
$L^6$ using H\"older's inequality. We estimate the $L^\infty$ norm through
energy and Bernstein's inequality.

 Hence 
\begin{equation} \label{modulation} 
\left| \int v_{1,\lambda_1} v_{2,\lambda_2} 
v_{3,\lambda_3} v_{4,\lambda_4} v_{5,\lambda_5}
\, dx\,  dt \right| 
\lesssim \lambda_{5}^{-\frac32} 
  \prod_{j=1}^5 \Vert  v_{j,\lambda_j} \Vert_{V^2} ,
\end{equation}
which implies the desired estimate. 

It remains to study  $\lambda_1 \lesssim \lambda_2 \lesssim \lambda_3 \lesssim 
\lambda_4 \lesssim \lambda_5 $, $\lambda_4 \ge 1.05 \lambda_1$. The most difficult case is $\lambda_5 \le 1.02 \lambda_2$ since otherwise we apply two stronger
bilinear estimates. For simplicity we consider $\lambda_1 << \lambda$ where $\lambda_2 = \lambda_5= \lambda$. We have to bound 
\[ \int_{\sum \xi_j=0}   \prod \hat u_{j, \lambda_j}(\xi_j) d\xi_2 d\xi_5 dt  \]
with $\xi_1 = - \sum_{j=2}^5 \xi_j $. 
We may restrict the integration to $\sum_{j=2}^5 \xi_j \sim \lambda_1$ and $\xi_j \sim \lambda$. By symmetry it suffices to consider 
\[   \int_{\sum \xi_j=0} \chi_{| |\xi_3|- |\xi_2| | \sim \lambda_1 }    \prod \hat v_{j, \lambda_j}(\xi_j) d\xi_2 d\xi_5 dt.  \]
We choose $\varepsilon>0$ small, $p,q$ so that $1/p= (1-\varepsilon)/2 + \varepsilon/3$, $1/q = \varepsilon/2 + (1-\varepsilon)/3$. By H\"older's inequality
\[ 
\begin{split} 
\left| \int v_{1,\lambda_1}  (v_{2,\lambda} v_{3,\lambda})_{\lambda_1}  
 v_{4,\lambda} v_{5,\lambda}
\, dx \,  dt \right| 
\lesssim  &  
\lambda^{-1} (\lambda^{\frac56} \lambda_1^{-\frac16})^{\varepsilon} 
\lambda^{-\frac12} \lambda_1^{-\frac12} ( \lambda^{\frac16} \lambda_1^{\frac12})^{1-\varepsilon} 
\lambda_3^{-\frac16}     \prod_j  
\Vert v_{j,\lambda_j} \Vert_{U^{\frac{12}5}}.
\\ \lesssim & 
\lambda^{-\frac32} (\lambda^{\frac23} \lambda_1^{-\frac23})^{\varepsilon} 
 \prod_j   \Vert v_{j,\lambda_j} \Vert_{U^{\frac{12}5}}.
\end{split} 
\] 

For the second part we would like to put one $v$ term into $L^6$, and
up to $2$ into $U^2$. This can be easily be done if there are two
frequencies of $v$ which differ by a small constant times
$\lambda_{max}$.  If not it is not hard to see that in the argument
above we can put one term into $L^6$.

\end{proof}

\subsection{Variants and Extensions of Well-posedness for the Quartic KdV Equation}
\label{variants}

The arguments of the last sections have implications for well-posedness
questions in other function spaces. Given $\omega  \in
C^1((0,\infty),(0,\infty))$,  $1\le p \le \infty$ and $T \in (0,
\infty]$ we define the function space
$X^\omega_{p,T}$ as the set of all distributions for which the norm \index{$X^\omega_{p,T}$}
\begin{equation}
\Vert u \Vert_{X^\omega_{p,T}}^p = \sum_{\lambda} ( \omega(\lambda)
\Vert u_\lambda \Vert_{V^2})^p ,
\end{equation}
with obvious modifications if $p=\infty$ is finite. We will always assume that
\begin{equation} 
\label{moderate}  \sup \frac{|\omega'|}{\omega} < \infty
\end{equation}
and \index{$\omega(\lambda)$}
\begin{equation}
\label{growth}   \inf \frac{\omega'}{\omega} > -1 .
\end{equation}

This is a Banach space
provided for some $C > 0$ we have
\begin{equation}
\label{eqn:omega-banach}
\liminf_{\lambda \to 0} \omega(\lambda) \lambda^{\frac12} > C,
\end{equation}
otherwise we obtain a Banach space of equivalence classes of functions.
Similarly, we define the Banach space \index{$Y^\omega_{p,T}$}
\begin{equation}
\Vert f \Vert_{Y^\omega_{p,T}}^p = \sum_{\lambda} ( \omega(\lambda)
\Vert f_\lambda \Vert_{DU^2})^p .
\end{equation}
The definition of $B^{\omega,p}_q$ follows the same pattern. 
It is not hard to see that
\[ \int u f dx dt \lesssim \Vert u \Vert_{X^{\omega}_{p,T}} \Vert f
\Vert_{Y^{\omega^{-1} }_{p',T}}
\]
and
\[ \Vert f \Vert_{Y^{\omega^{-1} }_{p',T}} \lesssim \sup_{ \Vert u
\Vert_{X^{\omega}_{p,T}}\le 1} \int uf dx dt .
\]
Moreover, we may expand the inner product into dyadic pieces and apply
uniformly elliptic pseudodifferential operators to the pieces. In
particular, we may replace differentiation by multiplication on the
dyadic pieces and vice versa.

\begin{prop}
\label{prop:Yomegapt}
The following estimate holds
\begin{equation}
\Vert \p_x (u^4) \Vert_{Y^\omega_{p,T}} \le
C \sup_{\lambda} \Vert u_\lambda \Vert_{L^6} \Vert u
\Vert_{\dot X^{-\frac16}_{\infty,T} }^2
\Vert u \Vert_{X^\omega_{p,T}} .
\end{equation}
\end{prop}

\begin{proof}
Given $v \in X^{\omega^{-1}}_{p',T}$, we expand $ \int \p_x (u^4) v dx dt $ into dyadic pieces, to which we apply
the arguments and estimate \eqref{eqn:multilin2-wp} from the previous section. By symmetry
\[ \sum_{\lambda_j} \lambda_5 \Big| \int u_{\lambda_1} u_{\lambda_2}   
 u_{\lambda_3} u_{\lambda_4} v_{\lambda_5} dx\,  dt \Big|
\lesssim \sum_{\lambda_1 \le \lambda_2 \le \lambda_3 \le \lambda_4,
\lambda_5 } \lambda_5 \Big| \int u_{\lambda_1} u_{\lambda_2}   
 u_{\lambda_3} u_{\lambda_4} v_{\lambda_5} dx\,  dt \Big|.
\]
If $\lambda_5 \sim \lambda_4$ we obtain
\[
\begin{split}
 \sum_{\lambda_1 \le \lambda_2 \le \lambda_3 \le \lambda_4 \sim
  \lambda_5 } \lambda_5 \Big| \int u_{\lambda_1} u_{\lambda_2}
u_{\lambda_3} u_{\lambda_4} v_{\lambda_5} dx\, dt \Big|
\lesssim  \sup_{\mu} \Vert u_\mu \Vert_{L^6}
 \Big( \sup_{\mu} \mu^{-\frac16} \Vert u_\mu \Vert_{U^2} 
\Big)^2
 \\  \times \sum_{\lambda_1 \le \lambda_2 \le \lambda_3 \le \lambda_4 \sim
\lambda_5}
     \lambda_1^{\frac16}  \lambda_5^{-\frac16}  
\left( \frac{\lambda_5}{\lambda_1} \right)^{\varepsilon} 
\Vert u_{\lambda_4} \Vert_{V^2}
     \Vert v_{\lambda_5}  \Vert_{V^2} ,
\end{split}
\]
which is bounded by
\[
\sup_{\mu} \Vert u_\mu \Vert_{L^6} \left( \sup_{\mu} \mu^{-\frac16} \Vert u_\mu
\Vert_{U^2} \right)^2
  \Vert u \Vert_{X^{\omega}_{p,T}} \Vert v \Vert_{X^{\omega^{-1}}_{p',T}}.
\]
The other extreme is
\[
\begin{split}
\sum_{\lambda_5 \le \lambda_1 \le \lambda_2 \le \lambda_3 \sim 
\lambda_4 } \lambda_5 \Big| \int u_{\lambda_1} u_{\lambda_2}
u_{\lambda_3} u_{\lambda_4} v_{\lambda_5} dx\, dt \Big|
\le   \sup_{\mu} \Vert u_\mu \Vert_{L^6}
 \left( \sup_{\mu} \mu^{-\frac16} \Vert u_\mu
\Vert_{U^2} \right)^2
 \\ \times   \sum_{\lambda_5 \le \lambda_1 \le \lambda_2 \le \lambda_3
\sim  \lambda_4} \lambda_5   \lambda_4^{-1}  \left( \frac{\lambda_4}{\lambda_5}\right)^\varepsilon \Vert u_{\lambda_4} \Vert_{V^2} \Vert v_{\lambda_5}  \Vert_{V^2} ,
\end{split}
\]
which satisfies the same  estimate provided
\[ \sum_{\lambda \le \mu} \lambda \omega(\lambda) \lesssim \mu 
\omega_{\mu} .\]
However, this is ensured by \eqref{growth}. The remaining cases are similar and the result follows. 
\end{proof}

From Proposition \ref{prop:Yomegapt}, we can prove the following corollary to Theorem \ref{thm:lwp}.

\begin{cor}\label{omega}  Suppose that $\omega$ satifies \eqref{moderate},
\eqref{growth} and \eqref{eqn:omega-banach}.  
If $\psi_0 \in \dot{B}^{-\frac16,2}_{\infty} \cap B^{\omega,2}_{p}$ is the initial data for a solution  of \eqref{eqn:kdv4} and
$v$ satisfies \eqref{smallness}, then the solution $\psi$ of Theorem
\ref{thm:lwp} is in
$ X^{\omega}_{p,T}$ and satisfies
\[ \Vert \psi \Vert_{X^\omega_{p,T}} \le C \Vert \psi_0 \Vert_{B^{\omega,2}_{p}} . \]
\end{cor}

In addition, we can show the following
\begin{cor} Suppose that $\psi_0$ lies in the closure of $C^\infty_0$ in
$\dot{B}^{-\frac16}_{\infty,T}$. Then, it follows that
\[  (t\to  \psi(t))  \in C([0,T] , \dot{B}^{-\frac16,2}_{\infty}) . \]
If $T=\infty$, then $e^{t\partial_{xxx}^3} \psi$ converges 
 to the scattering data as $t \to \infty$ in $\dot B^{-\frac16}_{\infty,T}$.
Moreover, if in addition $\psi_0 \in L^2$, then
\[  (t\to  \psi(t))  \in C([0,T] , L^2 )\]
and $e^{t\partial_{xxx}^3} \psi$ converges also in $L^2$. 
\end{cor}

There exists $\omega$ satisfying the assumptions above,
$\omega(\lambda) \lambda^{-\frac16} \to \infty$ as $\lambda \to \infty $ and
$\lambda \to 0$ and $\Vert \psi_0 \Vert_{B^{\omega,2}_{\infty}} <
\infty$. By corollary \ref{omega} the $X^\omega_{\infty,T}$ is
controlled by the initial data. Hence 
\[ \lambda^{-\frac16} \Vert v_\lambda \Vert_{DV^2\cap X^0} \to 0 \]
as $\lambda \to \infty$ or $\lambda \to 0$. By the previous argument the 
deviation of the solution to the linear solution tends to zero as the considered interval shrinks to zero. This implies continuity. Continuity at infinity 
always holds in $V^p$. 

The second part requires an obvious specialization of corollary \ref{omega} 
to the case $\omega=1$, Plus a repetition of the argument for scattering.

Particular examples for $\omega$ are $\la \lambda\ra ^s$ for $s \ge -\frac16$
and $ \lambda^s + \lambda^\sigma $ for $-\frac12 \le s \le -\frac16 \le \sigma$. It is not hard to see that we can replace the homogeneous
spaces
by inhomogeneous ones if we restrict to finite $T$ and allow the
constants to depend on $T$.

\section{Stability and Scattering for Perturbations of the Soliton}
\label{sec:nonlin-subcrit}
\subsection{Setup and main result}
We return now to the full nonlinear problem \eqref{eqn:kdv4}.  Let us take
\begin{eqnarray*}
\psi (x,t) = Q_{c(t)} (x - y (t)) + w(x,t).
\end{eqnarray*}
Then, we have
\[
\begin{split}  
\p_t w + \p_x (\p_x^2 w + 4 Q_c^3 w)  = & - \dot{c} (\partial_c Q_{c}) (x-y) + \dot{y} ( Q_{c}') (x-y) \\
& - \p_x ( \p_x^2 Q_c - c^2 Q_c + Q_c^4) - c^2 (Q_c' ( x-y)) \\
& - \p_x (6 Q_c^2 (x-y) w^2 + 4 Q_c (x-y) w^3 + w^4).
\end{split}
\]
Hence, 
\begin{equation}\label{deviation} 
\begin{split} 
\p_t w + \p_x (\p^2_{xx}  w + 4 Q^3_c w)  = & - \frac{\dot{c}}{c}  \tilde Q_c (x-y)  + (\dot{y} - c^2) Q_c' (x-y)  \\
& - \p_x (6 Q_c^2 (x-y) w^2 + 4 Q_c (x-y) w^3 + w^4).
\end{split} 
\end{equation}
In order to use the dispersive estimates proved in Section \ref{sec:lin-constant}, we wish to have
\begin{eqnarray}
\label{eqn:orthdes}
w \perp Q_c (x-y), \ w \perp Q_c' (x-y).
\end{eqnarray}

To get more regularity for $y$ and $c$  we ask for \eqref{eqn:orthdes} only asymptotically and hence take as in \eqref{eqn:coupledsystem-c} and \eqref{eqn:coupledsystem-x}
the modal equations
\begin{eqnarray}
\label{modal1}
\frac{\dot{c}}{c} \langle Q_c, \tilde{Q}_c \rangle =  \langle w, Q_c \rangle
\end{eqnarray}
and
\begin{eqnarray}
\label{modal2}
(\dot y -c^2) \langle Q_c',Q_c' \rangle = -\kappa \langle w , Q_c' \rangle  ,
\end{eqnarray}
where $\kappa > 0$ is taken to be large.  \index{$\kappa$}

We calculate
\begin{eqnarray*}
 \frac{d}{dt} \langle w,Q \rangle & = & \langle w_t,Q \rangle + \dot y \langle w,Q' \rangle  + \frac{\dot c}{c}  \langle w, \tilde Q \rangle \\
& = & \langle w, \Lop Q' \rangle - \frac{\dot c}{c} \langle Q,\tilde Q \rangle + \langle 6Q^2 w^2 + 4Q w^3+w^4, Q' \rangle \\
&&  + (\dot y - c^2) \langle w,Q' \rangle  + \frac{\dot c}{c}  \langle w, \tilde Q \rangle
\end{eqnarray*}
and
\begin{eqnarray*}
 \frac{d}{dt} \langle w,Q' \rangle & = & \langle w_t, Q' \rangle + \dot y \langle w,Q'' \rangle + \frac{\dot c}{c} \langle w, \tilde Q' \rangle \\
& = & \langle w, \Lop Q'' \rangle + (\dot y-c^2) \langle Q',Q' \rangle + \langle 6 Q^2 w^2+ 4Q w^3+w^4,Q'' \rangle \\
&& + (\dot y-c^2) \langle w,Q'' \rangle + \frac{\dot c}{c} \langle w, \tilde Q' \rangle.
\end{eqnarray*}
Hence,
\begin{equation} 
\label{mod1} 
 \frac{d}{dt} \langle w,Q \rangle + \langle w,Q \rangle 
= -\kappa \frac{\langle w,Q' \rangle^2}{\langle Q', Q' \rangle} - \frac{\langle w,Q \rangle \langle w, \tilde Q \rangle}{\langle Q,\tilde Q \rangle} + \langle 6Q^2 w^2 + 4Q w^3+w^4, Q' \rangle
\end{equation} 
and 
\begin{eqnarray}
\label{mod2} 
 \frac{d}{dt} \langle w,Q' \rangle &+& \kappa  \langle
w,Q' \rangle - \langle w, \Lop Q'' \rangle = \\
&& -\kappa \frac{\langle w,Q'\rangle \langle w,Q'' \rangle}{\langle Q', Q' \rangle} + \frac{\langle w,Q \rangle  \langle w, \tilde Q' \rangle }{\langle Q,\tilde Q \rangle}
 + \langle 6Q^2 w^2+ 4Qw^3+w^4,Q'' \rangle. \notag
\end{eqnarray}
The right hand sides are at least quadratic in $w$, 
and, as we shall see, small compared to $\Vert w_0 \Vert$ in a suitable sense.
As a consequence  the orthogonality conditions are approximately satisfied
for large $t$.  In addition, $\dot c $ and $\dot y - c^2 $ are small and continuous.  

We study the initial value problem $w(0) = w_0$. Let again $v$ be the solution
to the linear problem. We will prove scattering for small perturbations of the
soliton in $\dot B^{-\frac16,2}_\infty$. It will be important for the reverse
problem that we will achieve something slightly stronger.

Using the notation 
\begin{equation} \Gamma=\Big\{ y\in C([0,\infty)): y(0)=0, |\dot y -1| \le
  \frac{1}{10} \Big\} , \end{equation} we define for any interval, $I$, the
quantity
\begin{equation} \label{J}
J_I (v)  = \sup_\lambda \left(  
\Vert v_\lambda \Vert_{L^6_I}+  \lambda^{\frac14-\frac16} \Vert v_\lambda
\Vert_{L^4_I L^\infty}  + \sup_{ y \in \Gamma}
\lambda^{-\frac16}  \int_{\mathbb{R} \times I}  \gamma_0' (x-y(t)) 
 (v^2_\lambda+ (\partial_x v_\lambda)^2)  dx
dt \right) .
\end{equation}
\begin{prop} 
\label{smalllimit}
Let $v$ be a solution of \eqref{eqn:lkdv} with initial data $v_0 \in \dot
B^{-\frac16,2}_{\infty}$. Then,
\[ J_{[0,\infty)} (v) \lesssim \Vert v_0 \Vert_{\dot B^{-\frac16,2}_{\infty}}. \]
Moreover, if $v_0$ is in the closure of $C^\infty_0$ in $\dot B^{-\frac16,2}_{\infty}$ then
\[ \lim_{t\to \infty} J_{[t,\infty)} (v) = 0 . \]
\end{prop}

\begin{proof} The first statement is an immediate consequence of the
 Strichartz  estimate and local smoothing. For the second statement we fix
$\varepsilon>0$. There are at most finitely many $v_{0,\lambda}$ of norm
larger than $\varepsilon/c$. Hence it suffices to verify the statement for a
single $\lambda$. Since $v_\lambda \in L^6L^6$ and $L^4L^\infty$, we have
\[ \lim_{t\to \infty} \Vert v_\lambda  \Vert_{L^6(\mathbb{R}\times (t,\infty))} =
\lim_{t\to \infty} \Vert v_\lambda \Vert_{L^4_{[t,\infty)}L^\infty} =0. \]
Let $I$ be a bounded interval. Then the map
\[ \Gamma \to \lambda^{-\frac16}  \int_{\mathbb{R} \times I}
\gamma_0' (x-y(t)) (v_\lambda^2+ (\partial_x v_\lambda)^2)  dx\, dt  \]
is continuous with respect to uniform convergence, hence it assumes its
maximum. Given $j\ge 1$, let  $y_j : [2^j, 2^{j+1}] \to \mathbb{R}$ be the path
for which this quantity is maximal. We choose two paths $y_{o}, y_{e}$  with $\gamma(0) =0$ and the difference between $1$ and the derivative at most $.2$,
 one which coincides with $y_j$ for $j$ even on the corresponding intervals, and
one which does so for $j$ odd.  For both paths we have the local
smoothing estimate. But this implies the claim.
\end{proof}

Let $y \in \Gamma$. The function spaces $\dot X^s_{\infty;T}$ and
$\dot Y^s_{\infty;T}$ depend on $y$ but not on $c$. This dependence is
not reflected in the notation.  In addition, let $c \in
C^1([0,\infty))$. We assume \eqref{eqn:cysmall}, $ \dot c \in L^2 $ and
$\dot y -c^2 \in L^2 $ in this section, which we have to verify for the solutions we study,  and turn to a study of a priori
estimates for solutions to \eqref{deviation}, \eqref{modal1} and
\eqref{modal2} and recall \eqref{mod1} and \eqref{mod2}.  Because of
translation and scaling invariance we may restrict ourselves to a
study for $y(0)=0$ and $c(0)=1$. Moreover, we may and do assume that
the orthogonality conditions hold at time $0$, i.e.
\[ \langle w_0,Q \rangle = \langle w_0,Q' \rangle = 0. \] 
The main result is the following sharpened version of Theorem
\ref{thm:main-result}.

\begin{prop}
\label{prop:linear-ap}
Let $C>0$. There exist $\varepsilon>0$ and $K >0$ such that for
$\Vert w_0 \Vert_{\dot B^{-\frac16,2}_\infty} < C $ and $J_{[0,T)} (v) \le \varepsilon$ for $v$ a solution of \eqref{eqn:lkdv} with initial data $w_0$,
the solution $w$ in the system of equations \eqref{eqn:coupledsystem-c} -
\eqref{eqn:coupledsystem-w} satisfies \eqref{eqn:cysmall}, 
\[  \Vert w  \Vert_{\dot X^{-\frac16}_{\infty,T}} \le K  J^{\frac12}_{[0,T)} (v)  \]
with $K$ depending on $C$ but not on time.  
Moreover, if $J_{(0,\infty)}(v) \le \varepsilon$ then there exists a unique $\eta \in \dot B^{-\frac16,2}_{\infty} $ such that
\[ \lim_{t\to \infty} e^{t\p_{xxx}^3} w_\lambda (t) = \eta_\lambda , \]
with convergence in $L^2$. In addition
\[ \lim_{t\to \infty} \Vert w(t) \Vert_{\dot B^{-\frac16}_{\infty,T}} = \Vert \Psi \Vert_{\dot B^{-\frac16,2}_\infty}. \]
\end{prop} \index{$\Psi$}

\begin{rem}
Variants in the spirit of Corollary \ref{omega} can be easily obtained by including the arguments there, which will establish Theorems \ref{thm:main-result} and \ref{thm:inverse} with higher Sobolev regularity as stated in Remark \ref{rem:sobreg}.
\end{rem}

The proof consists of three step, a preliminary part consisting of an important initialization, 
multilinear estimates that are less critical variants of those of the last 
section, and a priori estimates for the nonlinear equation using the
multilinear estimates and the linearized equation.

We recall that $v$ satisfies $v_t + v_{xxx}=0$ with initial data $v(0)=
w_0$.   We want to control the difference between $v$  and the solution $\nu$ to 
\[ \nu_t + c^2 \partial_x \nu - \partial_x \Lop \nu = \alpha \tilde Q + \beta
Q' \]
with initial data $\nu(0) = w_0$ with $\alpha$ and $\beta$ ensuring
\[ \langle \nu, Q \rangle = \langle \nu, Q' \rangle = 0, \]
which we assume to hold  initially. We recall that \eqref{eqn:cysmall} is a standing assumption.

 For simplicity, let us define
\begin{eqnarray*}
J = J_{[0,\infty)} (v).
\end{eqnarray*}
The following result is the first step of the proof. 
\begin{lem}\label{lem:nu}  
Suppose that $w_0 \in B^{-\frac16,2}_\infty$ satisfies the orthogonality
conditions. Then 
\[ \Vert \nu \Vert_{\dot X^{-\frac16}_{\infty}} \lesssim \Vert w_0
\Vert_{\dot B^{-\frac16,2}_{\infty}}  \]
and 
\begin{equation} \label{eqn:nu-v} 
  \Vert \nu-  \tilde P P^\perp_{Q'} v  \Vert_{\dot X^{-\frac16}_{\infty}}+ \Vert \alpha \Vert_{L^1\cap L^2} +
  \Vert \beta \Vert_{L^2}    \lesssim J.
\end{equation} 
\end{lem} 
 
\begin{proof} The first bound on $\nu$ is an immediate consequence of
  Proposition \ref{littlewood}. The second statement is more delicate. 
As a first step we consider  $u = \tilde P P^{\perp} v$. It satisfies 
\begin{equation} 
  \sup_{\lambda}  \left( \Vert u_{\lambda} \Vert_{L^6}  
+ \lambda^{-\frac16+\frac14} \Vert u_{\lambda} \Vert_{L^4L^\infty} 
+ \lambda^{-\frac16} \Vert u_{\lambda} \Vert_{L^2 H^1_\rho}\right) \lesssim J,
\end{equation} 
since 
\[ \Vert \langle v,Q \rangle \Vert_{L^2 \cap L^6} + \Vert \langle v, Q'
\rangle \Vert_{L^2\cap L^6} \lesssim J.   \]
We calculate  
\begin{equation} \label{projnu} 
 \partial_t  u  + c^2 u_x - \partial_x \Lop_{c,y}  u =  G  , \qquad   u(0) = w_0 , 
\end{equation}  
where
\[ 
\begin{split} 
G = &  - 4 \partial_x (Q^3 u) - \left(\frac{d}{dt}\frac{ \langle v,Q\rangle}
{\langle  Q,\tilde Q \rangle } \right) \tilde Q 
 -\left(\frac{d}{dt} \frac{ \langle v,Q'\rangle}{\langle Q', Q' \rangle} \right) Q'
\\ & \quad  -  \frac{ \langle v,Q\rangle}{\langle
  Q,\tilde Q \rangle } (\partial_t + c^2 \partial_x - \partial_x \Lop ) \tilde
Q 
- \frac{ \langle v,Q'\rangle}{\langle
  Q',Q' \rangle } (\partial_t + c^2 \partial_x - \partial_x \Lop ) \tilde
Q'. 
\end{split} 
\]
We consider the terms separately.   Any derivative falling on
$\langle Q, \tilde Q \rangle$ or $\langle Q', Q' \rangle$ can be computed using 
\eqref{eqn:tildeqq} and \eqref{eqn:qprime}, yielding a factor $\frac{\dot c}{c}$. 
Next, 
\[ \frac{d}{dt}  \langle v, Q \rangle = 
\langle \partial_t v + c^2 v' , Q \rangle 
+ (c^2 -\dot y)\langle v, Q' \rangle + \frac{\dot c}{c} \langle v, \tilde Q
\rangle 
= -\langle Q^4_x , v \rangle+    (c^2 -\dot y)\langle v, Q' \rangle + \frac{\dot c}{c} \langle v, \tilde Q
\rangle     \]
and
\[\begin{split}  \frac{d}{dt} \langle v, Q' \rangle = & \langle \dot v + c^2 v', Q' \rangle 
+  (c^2 - \dot y ) \langle v, Q'' \rangle + \frac{\dot c}{c} \langle v, \tilde
Q' \rangle 
\\ = & - \langle v, \Lop Q'' \rangle 
-4\langle Q^3 Q'' , v \rangle +   (c^2 - \dot y ) \langle v, Q'' \rangle +
\frac{\dot c}{c} \langle v, \tilde Q' \rangle.
\end{split} 
\] 
Moreover,
\[ (\partial_t + c^2 \partial_x - \partial_x \Lop ) \tilde Q 
= \frac{\dot c}{c} \tilde{\tilde Q} + (c^2 -\dot y) \tilde Q' + 2 c^2 Q' 
\]
and
\[  (\partial_t + c^2 \partial_x - \partial_x \Lop )  Q' 
= \frac{\dot c}{c} \tilde Q' + (c^2 -\dot y) Q'' .
\] 
We write $ G = \alpha \tilde Q + \beta Q' + g$, 
where - using again \eqref{eqn:tildeqq} and \eqref{eqn:qprime} 
\[ \langle Q, \tilde Q \rangle \alpha = -\frac{\dot c}{c} \langle v, \tilde Q
\rangle - (c^2-\dot y) \langle v, Q' \rangle +\frac{\dot c}{3c}  
\langle v, Q \rangle +  \langle v,(Q^4)_x  \rangle,
 \] 
\[ 
\langle Q' , Q' \rangle \beta = - \frac{\dot c}{c} 
\langle v, \tilde Q' \rangle  - (c^2 -\dot y) \langle v, Q'' \rangle 
+ \frac{10\dot c}{3c} \langle v, Q' \rangle + 4 \langle v, Q^3 Q'' \rangle 
+ \langle v, \Lop_c Q''_c \rangle - 2 c^2 \langle v, Q' \rangle 
\]
and 
\[ g = - 4 \partial_x (Q^3 v)  - \frac{\langle v, Q \rangle}{ \langle Q, \tilde Q
  \rangle } 
 (\frac{\dot c}{c} \tilde{\tilde Q} + (c^2 -\dot y) \tilde Q')
- \frac{ \langle v, Q' \rangle}{ \langle Q', Q' \rangle } (\frac{\dot c}{c} 
\tilde Q'  + (c^2 - \dot y) Q'' ). 
\]
By Lemma \ref{initialu}, we have  $ \Vert g \Vert_{Y^0} \lesssim J$.   
The difference $w= \nu-u$ satisfies (abusing the notation slightly by denoting 
by $\alpha$ and $\beta$ new quantities) 
\[ w_t+ c^2 w_x - \partial_x \Lop w = \alpha \tilde Q + \beta Q'  -  g \]
with initial data $w(0)=0$ and again  by Lemma \ref{initialu} 
\[ \Vert w \Vert_{X^0\cap V^2} \lesssim \Vert g \Vert_{Y^0} \lesssim J. \] 
 We rewrite the equation for $\nu$ as 
\[ \nu_t + \nu_{xxx} = -\partial_x (4 Q \nu) + \alpha \tilde Q + \beta Q' =: F. \]
Decompose $\nu= u+ w$. We recall that 
\begin{equation} 
\label{eqn:nuorth}  \langle Q, \tilde Q \rangle \alpha  = - \frac{\dot c}{c} \langle \nu,
\tilde Q \rangle,
\end{equation} 
hence 
\[ \Vert \alpha \Vert_{L^1} + \Vert F \Vert_{\dot Y^{-\frac16}_{\infty} } \lesssim  (\Vert \frac{\dot c}{c} \Vert_{L^2} 
+ \Vert \dot y-c^2 \Vert_{L^2}) J . \]
 The $L^2$ bound for $\beta$
is simpler.  The estimates for the linear equation imply now \eqref{eqn:nu-v}.
\end{proof}

As it will be used in the sequel, we note the following simple
consequence of Lemma \ref{lem:nu}.  Namely, we have
\begin{equation} \label{eqn:Jnu} J_{c,y}(\nu) \lesssim J(v) ,\end{equation} 
where we denote by $J_{c,y}$ the quantity analogous to $J$, but for the given
path dictated by the $c,y$ modulation parameters.  After this nontrivial preliminary step we continue with the proof of
Proposition \ref{prop:linear-ap}. The strategy is to write the
equation in terms of  
\[ u = \Psi - Q_{c(t),y(t)} - \nu \]
and expand the nonlinearity. In the next step we study multilinear estimates, 
which in the last step are combined with Proposition \ref{littlewood} to obtain
the a priori estimates.

\subsection{Multilinear Estimates}
\label{sec:multlinest}

We proceed as for the initial value problem and bound multilinear expressions. 
In this section we collect nonlinear estimates in terms of the $V^2$ spaces in
order to prove Proposition  \ref{prop:linear-ap}. 

To begin, we have the following
\begin{lem} 
\label{L2com} 
Let $u$ be a tempered distribution and $u_\lambda$ its  frequency localization. Let $\phi$ be a Schwartz function.
Then,
\[ \Vert \phi  u_\lambda \Vert_{L^2} \lesssim \min\{ \lambda^{\frac12-\varepsilon}, \lambda^{-1}\} \left( \Vert   u_\lambda \Vert_{L^2(\gamma')}+ 
\Vert  \p_x u_\lambda \Vert_{L^2(\gamma')} \right).
\]
Here $\varepsilon$ is the constant of \eqref{eqn:gamma0}. 
\end{lem}

\begin{proof}
We begin with the case $\lambda \ge 1$, in which case we prove the stronger estimate where we replace $\phi$ by $\gamma'$ as defined in Section \ref{sec:linkdv}. Let $\chi \in C^\infty_0$ be supported
in $\{ \xi: \frac12 \le |\xi| \le 2\}$.
Then,
\[
\begin{split} 
 \sqrt{\gamma'} u_\lambda  = &  \sqrt{\gamma' } \p_x^{-1} \chi \left( \frac{\p_x}{\lambda} \right)   \p_x  u_\lambda =  \lambda^{-1} \sqrt{\gamma'} \left( \frac{\p_x}{\lambda} \right)^{-1}  \chi \left(
 \frac{\p_x}{\lambda} \right)  \p_x  u_\lambda
\\ = & \lambda^{-1}  \left( \frac{\p_x}{\lambda} \right)^{-1}  \chi \left(
 \frac{\p_x}{\lambda} \right) \sqrt{\gamma'} \p_x u_\lambda  
   + \lambda^{-1} \left[ \sqrt{\gamma'},  \left( \frac{\p_x}{\lambda} \right)^{-1} 
\chi \left( \frac{\p_x}{\lambda} \right) \right]   \p_x   u_\lambda ,
\end{split}
\]
where $ (\p_x/\lambda)^{-1}  \chi \left( \frac{\p_x}{\lambda} \right)   $ is
an $L^2$ bounded Fourier multiplier.  As a result,
\[ \left\Vert  \lambda^{-1}  \left( \frac{\p_x}{\lambda} \right)^{-1} \chi \left(
 \frac{\p_x}{\lambda} \right) \sqrt{\gamma'} \p_x u_\lambda  \right\Vert_{L^2} 
\lesssim \lambda^{-1} \Vert u_\lambda \Vert_{L^2(\gamma')}. 
\]
We estimate the second term on the right hand side using the adjoint of 
Lemma \ref{lem:commutator} with $a=1$ and $s=0$
\[ \left\Vert \lambda^{-1} \left[ \sqrt{\gamma'},  \left( \frac{\p_x}{\lambda} \right)^{-1} 
\chi \left( \frac{\p_x}{\lambda} \right) \right] (\gamma')^{-\frac12} \sqrt{\gamma'}  \p_x   u_\lambda \right\Vert_{L^2} 
\lesssim \lambda^{-1} \Vert u_{\lambda} \Vert_{L^2(\gamma')}. 
\]

We turn to $\lambda <1$. Clearly, 
\[ \Vert \phi u_\lambda \Vert_{L^2} \le \Vert \phi (\gamma')^{-\frac12}
\Vert_{L^2} \Vert \sqrt{\gamma'} u_\lambda \Vert_{L^\infty}.
\]
Let $\tilde{\chi} = \frac{\sin(x)}{x}$, which is the inverse Fourier
transform (up to a constant factor) of the characteristic function of the interval $[-1,1]$. Let
$x_0 \in \mathbb{R}$. We define
\[ g_\lambda (x) = u_\lambda (x) \tilde{\chi} \left( \frac{\lambda (
    x-x_0)}{100} \right).  \] 
Then, $g_\lambda$ satisfies roughly the
same frequency localization as $v_\lambda$, and it coincides with
$u_\lambda$ at $x_0$. Thus, by  Bernstein's inequalities
\[ |\sqrt{\gamma' (x_0)} u_\lambda(x_0)| \le c \lambda^{\frac12} \sqrt{\gamma'
(x_0)} \Vert g_\lambda \Vert_{L^2} \le c \lambda^{\frac12} \sup_{x,x_0} \frac{\sqrt{\gamma' (x_0)}
  \tilde{\chi} \left( \frac{\lambda (x-x_0)}{100} \right)}{\sqrt{\gamma'
  (x)}} \Vert \sqrt{\gamma'} u_\lambda \Vert_{L^2}.
\]
Now the elementary estimate
\[
\sup_{x,x_0}  \sqrt{\gamma' (x_0)} \frac{| \tilde{\chi} \left( \frac{\lambda (x-x_0)}{100} \right) | }{\sqrt{\gamma' (x)}}
\le c \lambda^{-\varepsilon} \]
completes the proof.
 \end{proof}

We proceed to prove the necessary multilinear estimates.  To begin, we
prove the following
\begin{lem}  
\label{nonlinear1}
Let $c$, $y$ satisfy \eqref{eqn:cysmall}, 
 $u \in \dot X^{-\frac16}_{\infty}$  and $\nu$, $Q$ be as in Proposition \ref{prop:linear-ap}.  Then, the following estimates hold:
\begin{equation} 
\label{cubic1}
\Vert  \p_x(u_1 u_2 u_3  Q) \Vert_{\dot Y^{-\frac16}_\infty}  \lesssim  \prod_{j=1}^3 \Vert u_j \Vert_{\dot X^{-\frac16}_\infty },
\end{equation} 
\begin{equation}\label{quad1}
\Vert \p_x(u_1 u_2 Q^2) \Vert_{\dot Y^{-\frac16}_{\infty} }
\lesssim  \prod_{j=1}^2  \Vert u_j \Vert_{\dot X^{-\frac16}_\infty},
\end{equation}
\begin{equation} 
\label{cubic2}
\Vert  \partial_x (\nu^2 u  Q) \Vert_{\dot Y^{-\frac16}_\infty}  \lesssim  
J^{\frac12} \Vert \nu \Vert_{\dot X^{-\frac16}_{\infty}}^{\frac32}    
\Vert u \Vert_{\dot X^{-\frac16}_\infty}
\end{equation} 
and
\begin{equation}\label{quad2}
\Vert \partial_x (\nu u  Q^2) \Vert_{\dot Y^{-\frac16}_{\infty} }
\lesssim  J^{\frac12} \Vert \nu \Vert^{\frac12}_{\dot X^{-\frac16}_{\infty}}  \Vert u   \Vert_{\dot X^{-\frac16}_{\infty}} .
\end{equation}
\end{lem} 

\begin{proof} 

We begin with the dual Strichartz estimate 
\[ \Vert f_\lambda  \Vert_{DV^2} \lesssim \lambda^{-\frac14} \Vert f_\lambda \Vert_{L^{\frac43}L^1}. \]
By construction, spatial Fourier multipliers in $V^p$, $U^p$, $DU^p$ and $DV^p$ are bounded by the supremum of the multiplier, hence  
\[ \Vert P_\lambda \p_x (Q^2 u_{1, \lambda_1} u_{2, \lambda_2}) \Vert_{DU^2} 
\lesssim \lambda^{\frac34} \Vert Q^2 u_{1, \lambda_1} u_{2, \lambda_2} \Vert_{L^{\frac43} L^1} 
\]
and 
\[
\begin{split} 
  \Vert  Q^2 u_{1, \lambda_1} u_{2, \lambda_2} \Vert_{L^{\frac43} L^1} 
\le & \Vert Q u_{1,\lambda_1} \Vert_{L^2 L^2} \Vert  Q u_{2,\lambda_2} 
\Vert_{L^2 L^2}^\frac12 \Vert u_{ 2,\lambda_2} \Vert_{L^\infty L^2}^{\frac12} 
\\  \lesssim & \min\{1,\lambda^{-1}_1\} \min \{1,\lambda_2^{-\frac1{2}}\} \lambda_1^{\frac16} \lambda_2^{\frac16} \Vert u_1 \Vert_{\dot X^{-\frac16}_\infty}  \Vert u_2 \Vert_{\dot X^{-\frac16}_\infty}. 
\end{split} 
\]
This is summable for $\lambda_j \in \Gamma$ and we obtain the desired estimate 
for $\lambda \le 1$. Assume now that $\lambda\ge 1$. Then, 
using H\"older and Bernstein  and $|Q'| \lesssim Q$
\[ \begin{split} \lambda^{-\frac16} \Vert P_\lambda \p_x Q^2 u_{1,\lambda_1} u_{2,\lambda_2} \Vert_{L^2 \gamma H^{-1}} 
\lesssim & \lambda^{-\frac16} \Vert Qu_{1,\lambda_1} \Vert_{L^\infty}^{\frac12} 
 \Vert Qu_{2,\lambda_2} \Vert_{L^\infty}^{\frac12}
\Vert Qu_{1,\lambda_1} \Vert_{L^2}^{\frac12}  
 \Vert Qu_{2,\lambda_2} \Vert_{L^2}^{\frac12}  
\\ \lesssim  & \lambda^{-\frac16} \lambda_1^{\frac5{12}} \min\{ 1,
\lambda_1^{-\frac12} \}   
\lambda_2^{\frac5{12}} \min\{ 1, \lambda_2^{-\frac12}\}
\Vert u_1 \Vert_{\dot X^{-\frac16}_{\infty}} \Vert u_2 \Vert_{\dot X^{-\frac16}_\infty} ,
\end{split} 
\] 
which can easily be summed over $\lambda_1$ and $\lambda_2$ if $\lambda \ge 1$.
This implies \eqref{quad1} and also \eqref{quad2}. 

We approach estimate \eqref{cubic1} similarly: We expand $u_j$, observe
that the expressions are symmetric and hence it suffices to sum over
$\lambda_1 \le \lambda_2 \le \lambda_3$. If $\lambda_1 \lesssim 1$ we
argue as above and  estimate  $u_{1,\lambda_1}$ in $L^\infty$, followed by
Bernstein's inequality. So we restrict to the case $\lambda_1 \gg 1$. 

Then, using that $Q$ is integrable, 
\[ \lambda \Vert P_\lambda Q u_{1,\lambda_1} u_{2,\lambda_2} u_{3,\lambda_3}\Vert_{DU^2}
\lesssim \lambda^{\frac34}  \prod_{j=1}^3 \Vert  u_{j,\lambda_j} \Vert_{L^4L^\infty} 
\lesssim \lambda^{\frac34} (\lambda_1 \lambda_2 \lambda_3)^{-\frac1{12}}    \prod_{j=1}^3 \Vert  u_{j,\lambda_j} \Vert_{\dot X^{-\frac16}_{\infty,T} }  ,  
\]
which is easily summable if $\lambda \lesssim 1 \lesssim \lambda_1, \lambda_2, \lambda_3$. If $\lambda >1$ we argue differently.
To simplify the argument we assume that the Fourier transform of $Q$ is supported in $[-1,1]$ - handling the tail is 
straight forward but technical.
 Instead of bounding $\lambda \Vert P_\lambda Q u_{1,\lambda_1} u_{2,\lambda_2}  u_{3,\lambda_3} \Vert_{DU^2}$, we employ duality 
and study 
\[ I= \left| \int Q u_{1,\lambda_1}  u_{2,\lambda_2}  u_{3,\lambda_3}  u_{4,\lambda_4}  dx dt \right| \]
assuming that $1 \ll \lambda_1 \le \lambda_2 \le \lambda_3$. 
Then, we have
\[
 I \le \Vert Q u_{3,\lambda_3} \Vert_{L^2} \Vert u_{1,\lambda_1} \Vert_{L^6} \Vert u_{2,\lambda_2} \Vert_{L^6} \Vert u_{4,\lambda_4} \Vert_{L^6}
\lesssim     \lambda_3^{-\frac56}\lambda_4^{-\frac16} 
 \prod_{j=1}^3   \Vert u_{j,\lambda_j} \Vert_{\dot X^{-\frac16}_{\infty,T}}
\Vert u_{4,\lambda_4} \Vert_{V^2}.
\]
The  factor  $  \lambda_3^{-\frac56} \lambda_4^{\frac56} $
is summable for fixed $\lambda_4$ over $1 \le \lambda_1 \le \lambda_2 \le \lambda_3 $, $1 \le \lambda_4 \lesssim \lambda_3$ - this suffices since $I=0$ if $\lambda_4$ is much larger than $\lambda_3$.
As a result, we have proven estimate \eqref{cubic1} and, after checking the proof, \eqref{cubic2}. 
 \end{proof} 

We turn to bounds for inner products occurring as inner products of the right hand side of \eqref{deviation} with $Q$ and $Q'$, and at the right hand side of \eqref{mod1} and \eqref{mod2}.  

\begin{lem}
\label{nonlinear2} 
Let $u \in \dot X^{-\frac16}_\infty$, and $v$, $Q$ be as in Proposition \ref{prop:linear-ap}.  In addition, let $\psi_0(t)$ be a one parameter family of Schwartz function parametrized by $t$ with 
uniformly bounded seminorms and $\psi(x,t) = \psi_0(t,x-y(t))$. 
Then for all $1\le p < \frac32$
\begin{equation} \label{quartic3} 
 \Vert  \langle \p_x (u_1 u_2 u_3 u_4) , \psi \rangle \Vert_{L^p}  \lesssim  \prod_{j=1}^4
\Vert u_j \Vert_{\dot X^{-\frac16}_\infty}, 
 \end{equation} 
where we consider  the $L^p$ norm with respect to time  and 
\begin{equation} \label{quartic4} 
 \Vert  \langle \p_x (v^3 u) , \psi \rangle \Vert_{L^p}  \lesssim  J  \Vert v
 \Vert^2_{\dot X^{-\frac16}_\infty} 
\Vert u \Vert_{\dot X^{-\frac16}_\infty}.
 \end{equation} 
For all $1\le p < 2$, we have
\begin{equation} \label{cubic3} 
 \Vert  \langle \p_x (u_1 u_2 u_3 Q) , \psi \rangle \Vert_{L^p}  \lesssim  \prod_{j=1}^3
\Vert u_j \Vert_{\dot X^{-\frac16}_\infty}  
\end{equation} 
and
\begin{equation} \label{cubic4} 
 \Vert  \langle \p_x (v^2u Q) , \psi \rangle \Vert_{L^p}  \lesssim  J \Vert v \Vert_{\dot X^{-\frac16}_\infty} \Vert u \Vert_{\dot X^{-\frac16}_\infty}.
\end{equation} 
For all $1\le p <3$, we have 
\begin{equation} \label{quadratic3} 
 \Vert  \langle \p_x (u_1 u_2  Q^2) , \psi \rangle \Vert_{L^p}  \lesssim  \prod_{j=1}^2
\Vert u_j \Vert_{\dot X^{-\frac16}_\infty} 
\end{equation}
and
\begin{equation} \label{quadratic4} 
 \Vert  \langle \p_x (v u  Q^2) , \psi \rangle  \Vert_{L^p}  \lesssim  J
\Vert u \Vert_{\dot X^{-\frac16}_\infty}.  
\end{equation}
\end{lem}
\begin{proof} 
We expand the terms in \eqref{quartic3} and we consider 
\[ I_p:=\Vert \langle u_{1,\lambda_1} u_{2,\lambda_2} u_{3,\lambda_3} u_{4,\lambda_4},
\psi \rangle \Vert_{L^p}. \]
By symmetry it suffices to look at the case $\lambda_1 \le \lambda_2 \le
\lambda_3 \le \lambda_4$. If $p=1$ we bound the terms using H\"older's and Bernstein's inequality as above: 
\[ 
\begin{split} 
I_1 \lesssim & \Vert u_{1,\lambda_1} \Vert_{L^\infty} \Vert u_{2,\lambda_2}
\Vert_{L^\infty} \Vert |\psi|^{\frac12} u_{3,\lambda_3} \Vert_{L^2} 
\Vert |\psi|^{\frac12} u_{4,\lambda_4} \Vert_{L^2} \\
\lesssim & \lambda_1^{\frac23} \lambda_2^{\frac23} \min\{\lambda_3^{\frac16},
\lambda_3^{-\frac56} \}  \min\{\lambda_4^{\frac16},
\lambda_4^{-\frac56} \} \prod_{j=1}^4 \Vert u_j \Vert_{\dot X^{-\frac16}_{\infty}} ,
\end{split} 
\]
which is easily summable.  
We obtain by H\"older's inequality 
\[ I_{\frac32} \lesssim \prod_{j=1}^4 \Vert u_{j,\lambda_j} \Vert_{L^6} 
\lesssim \prod_{j=1}^4 \Vert u_j \Vert_{\dot X^{-\frac16}_\infty} , 
\]
which we use if $1\le \lambda_1 \le \lambda_4$. If $\lambda_1\le 1$ we
estimate the corresponding term in $L^\infty$, apply Bernstein's
inequality, and argue as in the next case.  Interpolation with the
$L^1$ estimate yields a summable expression as long as $p< \frac32$.

We turn to estimate \eqref{cubic3}, denote again the $p$-norms by $I_p$  and expand again
\[ I_1  \lesssim   \Vert u_{1,\lambda_1}  \Vert_{L^\infty} 
\Vert Q u_{2,\lambda_2}  \Vert_{L^2} \Vert (\p_x \psi) u_{3,\lambda_3}\Vert_{L^2}  
\lesssim \lambda_1^{\frac23} \min\{\lambda_2^{\frac16},\lambda_2^{-\frac56}\}
  \min\{\lambda_3^{\frac16},\lambda_3^{-\frac56}\} \prod \Vert u_j \Vert_{\dot
    X^{-\frac16}_\infty} ,
\]
which again is easily summable over $\lambda_1\le \lambda_2 \le
\lambda_3$. Also
\[ I_2 \lesssim  \prod_{j=1}^3  \Vert u_{j, \lambda_j} \Vert_{L^6} \lesssim 
\prod_{j=1}^3 \Vert u_{j,\lambda_j} \Vert_{\dot X^{-\frac16}_\infty}, 
\]
which is almost summable, and by  interpolation  we obtain the bounds for any
$p <2$. 

The estimate \eqref{quadratic3} with $p=1$ follows by the same arguments as
above. It is even simpler. Again we may restrict ourselves to $\lambda_1 \le
\lambda_2$. For $p=3$ we put estimate $u_{j,\lambda_j}$ into $L^6$ and again
the full statement follows by interpolation. 
A simple check of the proof reveals that the arguments above imply
\eqref{quadratic4},  \eqref{cubic4} and \eqref{quartic4}.  
 \end{proof}

The right hand sides of \eqref{mod1} and \eqref{mod2} are  functions of $t$, 
for which  we have bounds in $L^p$ 
for $1\le  p < \frac32$ in terms  of $\Vert w \Vert_{\dot
  X^{-\frac16}_{\infty,T}}$. In the second equation \eqref{mod2} the term $\langle w, \Lop Q_{xx} \rangle$ plays a special role: It is in $L^q$ for $2\le q \le \infty$, 
but not in $L^p$ for any $p<2$ in general. In particular we cannot control the deviation of $y$ from the linear movement.

Equation \eqref{mod1} and \eqref{mod2} can be considered as scalar linear
ordinary differential equations for $\langle w, Q\rangle$ and $\langle w,
Q'\rangle$.  The kernel for the fundamental solution is uniformly bounded in
$L^p$ in the first case, for all $p$, and in the second case it is bounded in $L^1$ by $\frac1\kappa$ whereas the $L^\infty$ norm is $1$.

We collect the consequences in the following
\begin{lem} \label{prods} 
Suppose that $w$ solves \eqref{deviation} with $\langle w(0), Q\rangle = \langle w (0), Q' \rangle =0$ and $w=v+u$ where $v$ solves \eqref{linhom} with initial data $w (0)$. 
Then, 
\begin{equation}  
\sup_t  |\langle w(t),Q \rangle | \lesssim  
(J+ \Vert u \Vert_{\dot X^{-\frac16}_\infty})^2 (1+ \Vert w \Vert_{\dot X^{-\frac16}_\infty})^2 
 \end{equation} 
and
\begin{equation} 
\sup_t |\langle w(t),Q' \rangle | \lesssim  (J+ \Vert u \Vert_{\dot X^{-\frac16}_\infty})^2 (1+ \Vert w \Vert_{\dot X^{-\frac16}_\infty})^2
+ \kappa^{-\frac12 } \Vert w \Vert_{\dot X^{-\frac16}_{\infty,T}}.
\end{equation} 
Moreover, if $1\le p < \frac32$, then 
\[
\left\Vert \frac{d}{dt} \langle w(t),Q \rangle  \right\Vert_{L^p(0,\infty)} 
\lesssim  (J+ \Vert u \Vert_{\dot X^{-\frac16}_\infty})^2 (1+ \Vert w \Vert_{\dot X^{-\frac16}_\infty})^2.
\]
We may write $\frac{d}{dt} \langle w(t),Q' \rangle=\gamma_1+ \gamma_2$ 
such that  
\[
\Vert \gamma_1 \Vert_{L^p(0,\infty)} \lesssim 
(J+ \Vert u \Vert_{\dot X^{-\frac16}_\infty})^2 (1+ \Vert w \Vert_{\dot X^{-\frac16}_\infty})^2
\]
and 
\[ 
\Vert \gamma_2 \Vert_{L^2(0,\infty)} \lesssim  ( J+  \Vert  u \Vert_{\dot
  X^{-\frac16}_{\infty}}).
\]
Finally, it follows that
\[ \sup_t ( |c(t)-1| + |\dot c|)+ \Vert \dot c \Vert_{L^1}  \lesssim (J+ \Vert u \Vert_{\dot X^{-\frac16}_\infty})^2 (1+ \Vert w \Vert_{\dot X^{-\frac16}_\infty})^2 ,
\] 
\[
 \sup_t  |\dot y-1|  \lesssim (J+ \Vert u \Vert_{\dot X^{-\frac16}_\infty})^2 (1+ \Vert w \Vert_{\dot X^{-\frac16}_\infty})^2 + \kappa^{-\frac12} 
( J + \Vert u \Vert_{\dot X^{\frac16}_{\infty} } ) 
\] 
and
\begin{equation} 
\Vert \dot y - c^2 \Vert_{L^2} \lesssim (J+ \Vert u \Vert_{\dot X^{-\frac16}_{\infty}} ) 
(1+ \Vert w \Vert_{\dot X^{-\frac16}_\infty})^3 .
\end{equation} 
\end{lem} 

\begin{proof}
  This is an immediate consequence of Lemma \ref{nonlinear1} an basic
  properties of the simple ordinary differential equations.
\end{proof}

The estimate of this subsection remain true if we consider a time
integral instead of $(0,\infty)$.

\subsection{Global Bounds and Scattering Near the Soliton}
\label{sec:apps}

In this subsection  we complete the proof of  Proposition \ref{prop:linear-ap}.

\begin{proof}
  By the local existence result there exists a local solution in a
  neighborhood of the soliton. 

The decomposition $ \psi =
  Q_{c(t),y(t)} + w$ together with the modal equations \eqref{modal1} and \eqref{modal2}  
implies existence
  of $C^1$ functions $c(t)$ and $y(t)$ which satisfy \eqref{deviation}, \eqref{modal1} and \eqref{modal2} up to
  fixed time.  We recall  that after rescaling and shifting 
$\langle   w_0,Q \rangle = \langle w_0, Q'\rangle =  0 $
and $c(0)=1$, $y(0)=0$.

As in the first step  we denote the solution to the linear  equation with initial data $w(0)$ by $\nu$.  
It satisfies the estimates of Lemma \ref{lem:nu} and \eqref{eqn:Jnu}  provided 
\eqref{eqn:cysmall} is satisfied.

We suppose that  $(\psi,c,y)$ is a solution
up to time $T$, such that
\[  u = \psi - Q_{c(t),y(t)} - \nu \]
satisfies for some $k_1$, $k_2$ to be chosen later
\begin{equation}\label{aprioricondition} 
\Vert P_{Q'}^\perp  u \Vert_{\dot
    X^{-\frac16}_{\infty,T}} \le 2 k_1 J^{\frac12}  , \quad 
\Vert u \Vert_{\dot
    X^{-\frac16}_{\infty,T}} \le 2 k_2 J^{\frac12}.  \end{equation}
We shall see that there exist $\delta$, $k_1$ and $k_2$  such that if in addition $J  \le \delta$, then
\begin{equation} \label{improved} 
 \Vert P_{Q'}^\perp u \Vert_{\dot X^{-\frac16}_{\infty,T}} \le k_1 J^{\frac12}  ,
\quad \Vert u \Vert_{\dot X^{-\frac16}_{\infty,T}} \le k_2 J^{\frac12} .  
\end{equation}  
This implies the estimate conditionally depending on
\eqref{aprioricondition}. 
Observe that by Proposition \ref{prods} control of the norms implies 
validity of \eqref{eqn:cysmall} if $\delta$ is sufficiently small. 
In particularly the estimates on the linear equations hold.

On the other hand, if we fix $C$ and $\delta$ we can apply a
continuity argument with the initial data $\tau w_0$. The estimate
clearly holds for small $\tau$ and the norms depend (for finite time)
continuously on $\tau$. This implies the a priori estimate uniformly
for all $T$. The scattering statement  is an immediate consequence since functions in $V^2$  are left continuous at infinity combined with a frequency envelope argument as above. 
It remains to derive \eqref{improved} from \eqref{aprioricondition} 
for suitably chosen $k_1$, $k_2$ and $\delta$.

We formulate the crucial estimate in the following 
\begin{lem} \label{Qp} Let $C$ be given and $v$, $Q$ be as in Proposition \ref{prop:linear-ap}. There exist $k_1$, $k_2$ and $\delta$ 
such that, if \eqref{aprioricondition} holds, 
$\Vert w_0 \Vert_{\dot B^{-\frac16,2}_\infty} \le C$ and $J_{(0,T)} (v) \le \delta$ hold, then 
\[ 
\Vert P_{Q'}^{\perp} u \Vert_{\dot X^{-\frac16}_{\infty,T}} 
\le c_3 \left( \Vert u \Vert_{\dot X^{-\frac16}_{\infty,T}}^2 + J_{(0,T)}^{\frac12} (v) \Vert w_0
\Vert_{\dot B^{-\frac16,2}_\infty}^{\frac32} + J_{(0,T)} (v) \right) .
\]
\end{lem} 

We postpone the proof of Lemma \ref{Qp}. 
Clearly 
\[ \langle u, Q'\rangle = \langle w,Q'\rangle \]
and the same is true for its derivatives. By Proposition \ref{prods} and simple
properties of odes, we have with implicit constants depending on the size of the initial data
\[ \Vert \langle u, Q' \rangle \Vert_{L^1+L^2} +  \Vert \frac{d}{dt} \langle u,Q' \rangle \Vert_{L^1+L^2} \lesssim  
\Vert u \Vert_{\dot X^{-\frac16}_\infty}^2 + J + \Vert \langle w, Q'' \rangle \Vert_{L^2} 
\]
and 
\[  \Vert \langle w, Q'' \rangle \Vert_{L^2} \le \Vert \langle \nu , Q''\rangle  \Vert_{L^2} 
+  \Vert \langle  P^\perp_{Q'} u , Q'' \rangle \Vert_{L^2} .\]
Hence, 
\[ 
\Vert \langle u,Q' \rangle \Vert_{L^1+L^2} +  \Vert \frac{d}{dt} \langle u,Q' \rangle \Vert_{L^1+L^2} \lesssim  
\Vert u \Vert_{\dot X^{-\frac16}_\infty}^2 + J  + \Vert P^\perp_{Q'} u \Vert_{\dot X^{-\frac16}_\infty}.
\]
The crucial point is that the right hand side only contains the projection of
$u$, not $u$ itself. 
We obtain easily 
\[ 
\Vert (\partial_t + \partial_x^3) \gamma(t)  Q'   \Vert_{\dot Y^{-\frac16}_{\infty,T}}
\lesssim \Vert \gamma \Vert_{L^2} + \Vert \gamma' \Vert_{L^1+L^2} .
\]
As a result, using estimates similar to those in Lemma \ref{prods} we have
\begin{equation} \label{Q2} 
\left\Vert \frac{ \langle u,Q' \rangle}{\langle Q',Q' \rangle} Q' 
\right\Vert_{X^{-\frac16}_{\infty,T}} \le k_4  
 \left(\Vert  u \Vert_{\dot X^{-\frac16}_{\infty,T}}^2 + J + \Vert P^\perp u \Vert_{\dot X^{-\frac16}_{\infty,T}}\right).
\end{equation}   

By Lemma \ref{Qp} and \eqref{aprioricondition}, 
\[ 
\Vert P_{Q'}^\perp u \Vert_{\dot X^{-\frac16}_{\infty,T}} \le 
c_3 (\Vert u \Vert_{\dot X^{-\frac16}_{\infty,T}}^2 + J_{(0,T)}^{\frac12} (v)
 + J ) 
\le c_3 (4k_2^2+1) J +
c_3   J^{\frac12} 
\] 
using, as we may, 
  $\Vert u \Vert_{\dot X^{-\frac16}_{\infty,T}}\le 1$,
and, by the estimate \eqref{Q2} and \eqref{aprioricondition} 
\[ 
\Vert u \Vert_{\dot X^{-\frac16}_{\infty,T}} \le k_4 (J  +  
+    c_3 (8k_2^2+1) J (v) + c_3   J^{\frac12} )). 
\]
We choose first $ k_1$, then  $k_2$ and finally $\delta$ small
to complete the proof. 
\end{proof}

It remains to prove Lemma \ref{Qp}. 
\begin{proof}  
We write the equation for $u=w-\nu$, $u_t +c^2 u_x - \partial_x \Lop_{c,y} u =: G$
\[ G =  (\frac{\dot c}{c} + \alpha_1)  \tilde Q + (\dot y-c^2 + \beta_1 ) Q' 
 - \partial_x ( 6 Q_c^2 (u+\nu)^2 + 4 Q_c'(u+\nu)^3 + (u+\nu)^4   )   ,
\]
where $\alpha_1$, $\beta_1$ ensure the orthogonality conditions for
$\nu$, i.e. \eqref{eqn:nuorth}.  We recall that they satisfy
\[ \Vert \alpha_1 \Vert_{L^1} \lesssim  \Vert \dot c \Vert_{L^2} J \]
and
\[ \Vert \alpha_1 \Vert_{L^2} + \Vert \beta_1 \Vert_{L^2} \lesssim  J. \]
In order to apply Proposition \ref{littlewood} we have to project $u$. This
leads to a calculation similar to Lemma \ref{lem:nu}. Let $\mu= P^\perp_{Q'} \tilde P
u$ and $\mu_t + c^2 \partial_x \mu - \partial_x \Lop \mu =: H$. Then, using 
$\langle u, Q \rangle = \langle w,Q \rangle$ and $\langle u,Q' \rangle =
\langle w, Q' \rangle$,  
\begin{eqnarray*} 
H  & = &  G - 
\left(\frac{d}{dt} \frac{\langle w, Q \rangle}{\langle Q, \tilde Q \rangle}
  \right) \tilde Q 
-  \left(\frac{d}{dt} \frac{\langle w, Q'
      \rangle}{\langle Q', Q' \rangle} \right)  Q'  \\ 
  &  & - \frac{\langle w,Q \rangle}{\langle Q,\tilde Q \rangle} 
\left( \frac{\dot c}{c} \tilde{ \tilde Q} + (c^2-\dot y) \tilde Q' + 2 Q'
\right) - \frac{\langle w,Q' \rangle}{\langle Q', Q' \rangle} 
\left( \frac{\dot c}{c} \tilde Q' + (c^2-\dot y)  Q''
\right) 
\\  & = & \alpha \tilde Q + \beta Q' + g, 
\end{eqnarray*}
where 
\[ 
\begin{split} 
-g= &  \partial_x ( 6 Q_c^2 (u+\nu)^2 + 4 Q_c'(u+\nu)^3 + (u+\nu)^4   )
+ \frac{\langle w,Q \rangle}{\langle Q,\tilde Q \rangle}\frac{\dot c}{c}
\tilde{\tilde{Q}} 
\\ & + \left( \frac{\langle w,Q \rangle}{\langle Q,\tilde Q \rangle} (c^2-\dot y)
+ \frac{\langle w,Q' \rangle}{\langle Q', Q' \rangle}  \frac{\dot c}{c} \right) 
\tilde Q' 
 +  \frac{\langle w,Q' \rangle}{\langle Q', Q' \rangle} (c^2-\dot y) Q''.
\end{split} 
\]
By construction $u(0)=0$. We apply Proposition \ref{littlewood} 
\[ 
 \Vert  u \Vert_{\dot X^s_{\infty,T}} 
\lesssim \Vert g \Vert_{\dot Y^{-\frac16}_{\infty,T}} 
+ \Vert \langle g,Q \rangle \Vert_{L^1} +
\Vert \langle g^+,Q \rangle \Vert_{L^2} + \Vert \langle g,Q' \rangle
\Vert_{L^2}.
\]  
By Lemma \ref{lem:wp},  Lemma \ref{nonlinear1} and Lemma \ref{nonlinear2} 
we get 
\[ \Vert g \Vert_{\dot Y^{-\frac16}_\infty} \lesssim \Vert u \Vert_{\dot
  X^{-\frac16}_\infty}^2 + J \Vert w_0 \Vert_{\dot B^{-\frac16,2}_\infty}, 
\]
and by Proposition \ref{prods} 
\[ \Vert \langle  g, Q \rangle \Vert_{L^1} + \Vert \langle  g^+, Q \rangle \Vert_{L^2} 
+ \Vert \langle g, Q \rangle \Vert_{L^2} \lesssim  \Vert u \Vert_{\dot
  X^{-\frac16}_\infty}^2 + J^{\frac12}  \Vert w_0 \Vert^{\frac32}_{\dot
  B^{-\frac16,2}_\infty}.
\] 
Together, we have 
\[ 
\Vert \tilde P^* P^\perp  u \Vert_{\dot X^{-\frac16}_{\infty,T}}
\lesssim  \Vert u \Vert_{\dot
  X^{-\frac16}_{\infty,T}}^2 + J^{\frac12}  \Vert w_0 \Vert^{\frac32}_{\dot B^{-\frac16,2}_\infty} .
\]
\end{proof}

It is straight forward to generalize Proposition \ref{prop:linear-ap} 
to smaller function spaces in the style of  Subsection \ref{variants}.

\subsection{An Almost Inverse Wave Operator Result}
\label{sec:inverse}

In this section we will construct solutions with given asymptotic behavior, proving Theorem \ref{thm:inverse}.  This is a partial converse statement to Proposition \ref{prop:linear-ap}.

\begin{rem}
  Theorem \ref{thm:inverse} is quite satisfactory in several respects. It
  shows which asymptotic properties may characterize a solution. The
  main missing piece is uniqueness of the solution $\Psi$.  It implies
  existence of a solution for small scattering data, and, for arbitrary
  scattering states, existence of a solution with given scattering
  data for large $t$.
\end{rem}

\begin{proof}

We turn  to the - time reversed - equation
\begin{equation}
\label{equ}
\p_t w + \p_x (\p_x^2 w + 4 Q_c^3 w)  =  (\dot y -c^2) \langle w,Q_{xx} \rangle + \frac{\dot c}{c}  \langle w, \tilde Q' \rangle  +   6Q^2 w^2+ 4Qw^3+w^4,Q_{xx} , \notag
\end{equation}
with
\[
 \frac{\dot{c}}{c} \langle Q_c, \tilde{Q}_c \rangle  = -\langle w,
Q_c \rangle, \qquad (\dot y -c^2) \langle Q_c',Q_c' \rangle  = 
\kappa \langle w , Q_c' \rangle .
\]

Let $v$ be the solution to the Airy equation with initial data
$v_0$.  We may and do assume that $y_0=0$.  By Proposition
\ref{smalllimit} we know that $\lim_{t\to \infty} J_{[t,\infty)}  (v) = 0$. 
 Given $S>0$
and $y^S$ satisfying $|y^S(S) - c_\infty^2 S| < \hat\delta S$,
we solve the backwards initial value problem
\[ \Psi(S) = v(S) + Q_{c_\infty,y^S}.  \] 
We choose $1>> \hat \delta >> \delta$  to ensure that $|\dot y - c^\infty| \le \hat \delta$ for the solutions under consideration. 
The arguments of the
previous section allow to do that down to a largest  time $t^{S,y^S}$ for which
\[ |y(t^{S,y^S}) - c_\infty^2 t^{S,y^S}| = \hat \delta t^{S,y^S} .\]  
We want to show that  the infimum of the $t^{S,y^S}$ as a function of 
$y^S$ is attained for some $y^S$ 
and it is equal to zero if $\hat \delta$ is sufficiently small.  
Suppose not, and denote the infimum by $\tau>0$. By
continuous dependence on $y^S$, given $\varepsilon>0$, 
there exists an interval $[a,b]$ so
that the solution exists down to a time smaller than 
 $(1+\varepsilon) \tau$, and 
$   y^{S,a}((1+\varepsilon) \tau)  = (c_\infty^2-\delta)  (1+\varepsilon) \tau$ ,
$y^{S,b}((1+\varepsilon) \tau)  = (c_\infty^2+\delta)  (1+\varepsilon) \tau$ . 
Hence, there exists $y^{S,\varepsilon}$ with
\[  y^{S,a}((1+\varepsilon) \tau)  = c_\infty^2(1+\varepsilon) \tau. \]
But then, if $\hat \delta$ is sufficiently small, we see that a positive infimum is
not possible, and moreover this construction gives a limit which is a solution
denoted again by $(\Psi^S, y^S)$ with $y^S(0)=0$.

We consider the limit $S \to \infty$. Since $\dot y^S-c_\infty^2$ and $\dot c^S$ are small
there exists a converging subsequence $y^{S_j}$, $c^{S_j}$, $S_j \to \infty$
which converges to $c$ and $y$. There are corresponding solutions
$\Psi_j$, $u_j$ and $w_j$ of the corresponding equation. We extend $w_j$
beyond $S_j$ by $v$. 
By the stability  result, given $\delta>0$ we find
$T> 0$ such that
\[  \Vert w_j - v  \Vert_{\dot X^{-\frac16}_{\infty,[T,\infty)}} \le \delta. \]
Using a frequency envelope there exists $\Lambda$ so that
\[ \lambda^{-\frac16} \Vert (w_j)_\lambda \Vert_{V^2(T,\infty)} \lesssim
\delta, \] whenever $\lambda> \Lambda$ or $\lambda^{-1} > \Lambda$.

In particular, 
\[ \Vert (w_j - w_l)(t)  \Vert_{\dot{B}^{-\frac16,2}_{\infty}} \le \delta 
\]
for $t \ge T(\delta)$ and $j,l$ sufficiently big. Again, using $J$ small 
we are able to deduce that $(w_j)$ is a Cauchy sequence in $\dot X^{-\frac16}_T$ 
and the limit is the desired solution. 
\end{proof}

\begin{appendix}

\section{Set-up and Properties of the $U^p$, $V^p$ Spaces for the Linear KdV Equation}
\label{app:u2v2}

To define the function spaces $U^2$, $V^2$, we summarize Section $2$ of Hadac-Herr-Koch \cite{HHK}, where we suggest the reader look for further details.  Let $\mathcal{Z}$ \index{$\mathcal{Z}$} be the set of finite partitions
$ -\infty <t_0 <t_1<\ldots<t_K=\infty$.  In the
following, we consider functions taking values in $L^2:=L^2(\RR^d;\CC)$,
but in the general part of this section $L^2$ may be replaced by an
arbitrary Hilbert space.
\begin{defn}
\label{def:u}
  Let $1\leq p <\infty$. For $\{t_k\}_{k=0}^K \in \mathcal{Z}$ and
  $\{\phi_k\}_{k=0}^{K-1} \subset L^2$ with
  $\sum_{k=0}^{K-1}\|\phi_k\|_{L^2}^p=1$  we call the
  function $a:\RR \to L^2$ given by
  \begin{equation*}
    a=\sum_{k=1}^{K} \chi_{[t_{k-1},t_{k})}\phi_{k-1}
  \end{equation*}
  a $U^p$-atom where $\chi_I$ is the standard cut-off function to interval
 $I$. Furthermore, we define the atomic space
  \begin{equation*}
    U^p:=\left\{u=\sum_{j=1}^\infty \lambda_j a_j \;\Big|\; a_j \ 
\text{a $U^p$-atom},\;
      \lambda_j\in \CC \text{ s.th. } \sum_{j=1}^\infty |\lambda_j|<\infty
\right\}
  \end{equation*}
  with norm
  \begin{equation}
\label{eq:norm_u}
    \|u\|_{U^p}:=\inf \left\{\sum_{j=1}^\infty |\lambda_j|
      \;\Big|\; u=\sum_{j=1}^\infty \lambda_j a_j,
      \,\lambda_j\in \CC,\; a_j \ \text{a $U^p$-atom}\right\}.
  \end{equation}
\end{defn}
Atoms are bounded in the supremum norm, and hence every convergence here
implies uniform convergence. 

\begin{prop}
\label{prop:u}
  Let $1\leq p< q < \infty$.
  \begin{enumerate}
  \item\label{it:u_banach} The expression $\| . \|_{U^p}$ is a norm. The space
    $U^p$ is complete and hence a Banach space.
  \item\label{it:u_emb} The embeddings $U^p\subset U^q$ have norm $1$. 
  \item\label{it:u_right_cont} For $u \in U^p$ all one sided limits exist, 
including at $\pm \infty$,  $u$ is continuous from the right, and the limit at 
$-\infty$ is zero.   
  \item\label{it:u_c} The subspace of continuous functions  $U^p_c$
is closed.  \index{$U^p$!definition}
  \end{enumerate}
\end{prop}

\begin{defn}
\label{def:v}

  Let $1\leq p<\infty$. We define $V^p$ as the normed space of all
  functions $v:\RR \to L^2$ for which the norm
  \begin{equation}
\label{eq:hom_norm_v}
    \|v\|_{V^p}:=\sup_{\{t_k\}_{k=0}^K \in \mathcal{Z}}
    \left(\sum_{k=1}^{K}
      \|v(t_{k})-v(t_{k-1})\|_{L^2}^p\right)^{\frac{1}{p}}
  \end{equation}
  is finite.  Here we understand $v(\infty)$ as zero. 
Let $V^p_-$ denote the subspace of all right continuous functions with limit 
$0$ at $-\infty$.  \index{$V^p$!definition}
\end{defn}
Taking the partition $\{t,\infty\}$ one sees that the supremum norm is not
larger than the $V^p$ norm.

\begin{prop}
\label{prop:v}
Let $1\leq p<q <\infty$.
  \begin{enumerate}
   \item The expression $\| . \|_{V^p}$ is a norm and $V^p$ is complete. 
  \item\label{it:v_limits} For $v \in V^p$ all one sided limits including at
    $\pm \infty$ exist. 
\item\label{it:v_spaces} The subspace $V^p_-$ is closed. 
\item\label{it:v_emb1} The embedding $U^p\subset V_{-}^p$ is
  continuous and $\Vert u \Vert_{V^p} \le  2^{1/p}  \Vert u \Vert_{U^p}$. 
\item\label{it:v_emb2} The embeddings $V^p\subset V^q$ are continuous and 
$\Vert v \Vert_{V^q} \le \Vert v \Vert_{V^p} $.
\end{enumerate}
\end{prop}

From the proof of Proposition $2.17$ of Hadac-Herr-Koch \cite{HHK}, we have the following
\begin{lem}
\label{lem:UVlem}
Let $f \in V^p_-$, $q > p$.  Then, given $\delta > 0$ and $m>1$, there exist $f_1 \in U^p$ and $f_2 \in U^q$ such that $f = f_1 + f_2$ and
\begin{eqnarray*}
m^{-1} \| f_1 \|_{U^p} + e^{ \delta m} \| f_2 \|_{U^q} \lesssim \| f \|_{V^p}.
\end{eqnarray*}
\end{lem}

The following corollary is obvious. 
\begin{cor}
The space $V^p_{-}$ is continuously embedded in $U^q$ for $q>p$.
\end{cor}

There is a bilinear map, $B$, which for $1/p+1/q=1$, $1< p,q < \infty$ can formally be written as 
\[B(f,g) = -\int f_t g dt, \]
for $ f \in V^p$, $g \in U^q$.
It satisfies 
\[ |B(f,g)| \le \Vert f \Vert_{V^p} \Vert g \Vert_{U^q} , \] 
which is natural if we replace $g$ by an atom. The map
\[ V^p \ni f \to ( g \to B(f,g)) \in (U^q)^* \] 
is an isometric bijection. Moreover, 
\[ \Vert u \Vert_{U^p} = \sup \{ B(u,v): v \in C(\mathbb{R}), \Vert v
\Vert_{V^p}= 1 \}.   \]
If $v \in V^p_-$, then 
\[ \Vert v \Vert_{V^q} = \sup \{ B(u,v): u \in C(\mathbb{R}), \Vert u
\Vert_{U^p}= 1 \}  .\]
If the distributional derivative of $u$ is in $L^1$ and $v \in V^p$, then 
\[ B(u,v) = -\int u_t v dt. \]
Given $f\in L^1$, then $F(t)= \int_{-\infty}^t f ds  \in V^p$ for all $ p \ge
1$, and  hence in $F \in U^p$. Moreover, $\Vert f \Vert_{DU^p} := \Vert F \Vert_{U^p}
\le \Vert f \Vert_{L^1}$. We denote by $DU^p$ the metric completion of $L^1$
in the norm given by the duality pairing. Similarly we define $DV^q$. 

 There is a close relation to Besov spaces, namely
\begin{equation} \label{embed}
 B^{\frac1p,p}_1 \subset U^p \subset V^p \subset B^{\frac1p,p}_\infty 
\end{equation} 
with continuous embeddings. These embeddings  clarify the relation to 
$X^{s,b}$ spaces below. 

We claim that the convolution with an $L^1$  function $\eta$ 
  defines a bounded operator on $U^p$ and $V^p$ with norm $\le \Vert \eta
  \Vert_{L^1} $.   Because of the duality statement it suffices to verify boundedness on $U^p$.
  We approximate the characteristic function by a sum of Dirac measures.  The
  convolution with an atom clearly has norm at most $1$.  Convergence in $U^1$
  to the convolution with the characteristic function is immediate. The full
  statement is an immediate consequence, as well as the boundedness of the
  convolution by a Schwarz function on $U^p$ and $V^p$.  In particular smooth
  projections on high and low frequencies are bounded.

 Following Bourgain's strategy  for the Fourier restriction spaces we  define the adapted function spaces
\begin{eqnarray*}
U^p_{KdV} & = & S(-t) U^p, \\
V^p_{KdV} & = &  S(-t) V^p
\end{eqnarray*}
and similarly $DU^p$ and $DV^p$.

Again, we define a bilinear map $B_{KdV}$ such that for $u \in V^p_{KdV} , v
\in U^q_{KdV}$, we have for function $u$ with $(\p_t + \p_{x}^3)u \in L^1L^2$ 
\begin{eqnarray*}
B_{KdV} (u,v) =- \int \langle  (\p_t + \p_x^3)u, v \rangle dt .
\end{eqnarray*}
Note, this bilinear map is well-defined and gives a duality relation.
Hence,
\begin{eqnarray*}
\| u \|_{DV^p_{KdV}} = \sup_{\| f \|_{U^q_{KdV}} \leq 1} \int u f dx dt , \\
\| u \|_{DU^p_{KdV}} = \sup_{\| f \|_{V^q_{KdV}} \leq 1} \int u f dx dt .
\end{eqnarray*} 
Moreover, we may restrict $f$ to suitable subspaces.  More details on
how the construction of such atomic spaces allows us to put $u_t$ in
the dual space are included in Hadac-Herr-Koch \cite{HHK}.

By the construction of our spaces we obtain for a solution $u$ of the
linear KdV equation
\begin{eqnarray}
\label{lin:kdv}
\left\{ \begin{array}{c}
u_t + u_{xxx} = f, \\
u(0,x) = u_0 (x) ,
\end{array} \right.
\end{eqnarray}
the estimates
\begin{eqnarray}
\label{eqn:linkdv1}
\| u \|_{V^2_{KdV}} \lesssim \| u_0 \|_{L^2} + \| f \|_{DV^2_{KdV}} 
\end{eqnarray}
and
\begin{eqnarray}
\label{eqn:linkdv2}
\| u \|_{U^2_{KdV}} \lesssim \| u_0 \|_{L^2} + \| f \|_{DU^2_{KdV}} ,
\end{eqnarray}
which follow trivially from the construction of the $V^2_{KdV}, DV^2_{KdV}$
and $U^2_{KdV}$, $DU^2_{KdV}$ spaces.

Spatial Fourier multipliers act on $U^p$, $V^q$, $DU^p$, $DV^q$ in the obvious
way and their operator norm is bounded by the supremum of the multiplier. 

Let $(p,q)$ be a Strichartz pair. Then, 
\[ \Vert  u \Vert_{L^pL^q} \le c \Vert |D|^{-\frac{1}{p}}  u \Vert_{U^p} \]
and the dual estimate
\[ \Vert f \Vert_{DV^{p'}} \le c \Vert |D|^{-\frac{1}{p}} f \Vert_{L^{p'}L^{q'}} \]
  hold. The first estimate is not hard to check on atoms. Since convergence in
  $U^p$ and in $L^pL^q$ both imply pointwise convergence for subsequences we
  obtain the full estimate. The second estimate follows by duality.

Similarly the local smoothing estimates carry over to $U^p$ spaces and to
$DV^q$. Let $c(t)$ and $y(t)$ satisfy \eqref{eqn:cysmall} 
Then 
\[ \Vert u \Vert_{L^2 H^1_{\sqrt{\gamma'}})} \le c \Vert u \Vert_{U^2} \]
and 
\[ \Vert f \Vert_{DV^2} \le c \Vert f \Vert_{L^2 \sqrt{\gamma'} H^{-1}}. \] 
In the same fashion the bilinear estimates for solutions to the free
equation imply bilinear estimates for functions in   $U^2$.

The smooth decomposition into high and low modulation (i.e. the smooth
projection of the frequencies to $\tau-\xi^3$ large respectively small) is
bounded in $U^2$ and $V^2$, and the $L^2$ norm of the high modulation part gains the
inverse of square root of the truncation as factor by the embeddings \eqref{embed}.

\end{appendix}

\printindex


\begin{thebibliography}{XX}

\bibitem{AAR} G.E. Andrews, R. Askey, R. Roy, {\em Special Functions}, Cambridge University Press, Cambridge (1999).

\bibitem{BeLi} {\sc H. Berestycki and P. L. Lion}.  {Nonlinear scalar field equations, I: Existence of a ground state,}  {\it Arch. Rational Mech. Anal.}, {\bf 82}, no. 4, 313-345 (1983).

\bibitem{Cote} R. C\^ote, {\em Construction of solutions to the subcritical gKdV equation with a given asymptotical behavior}, J. of Func. Anal., {\bf 241}, No. 1 (2006), 143-211.

\bibitem{Fed} M. Fedoryuk, {\em Asymptotic Analysis}, Springer-Verlag (1993).

\bibitem{Gr} A. Gr\"unrock, {\em A bilinear Airy estimate with application to gKdV-3}, Differential Integral Equations, {\bf 18}, No. 12 (2005), 1333-1339.

\bibitem{HHK} M. Hadac, S. Herr and H. Koch, {\em Well-posedness and
    scattering for the KP-II equation in a critical space}, Ann. de
  l'Inst. H. Poincar\'e - Analyse non Lin\'eaire, {\bf 26}, No. 3
  (2009), 917-941.  Erratum published at {\it http://dx.doi.org/10.1016/j.anihpc.2010.01.006}.

\bibitem{Kato} T. Kato, {\em On the Cauchy problem for the (generalized) KdV equation}. Stud. Appl. Math. Adv. Math. Suppl. Stud., {\bf 8} (1983), 93-128.

\bibitem{KPV} C. Kenig, G. Ponce and L. Vega, {\em Well-posedness and
    scattering for the generalized Korteweg-de Vries equation via the
    contraction principle}.  Comm. Pure Appl. Math., {\bf 46}, No. 4
  (1993), 527-620.

\bibitem{KT1} H. Koch and D. Tataru.  {\em Dispersive estimates for principally normal pseudodifferential operators},  Commun. Pure Appl. Math., {\bf 58}, No. 2 (2005), 217-284.

\bibitem{KT2} H. Koch and D. Tataru.  {\em A priori bounds for the $1D$ cubic NLS in Negative Sobolev Spaces},  Int. Math. Res. Not. (2007), Article ID rnm053.

\bibitem{Lamb} G.L. Lamb.  {\em Elements of Soliton Theory},  John Wiley \& Sons, New York (1980).

\bibitem{Martel} Y. Martel, {\em Linear problems related to asymptotic stability of solitons of the generalized KdV equations,} SIAM J. Math. Anal. {\bf 38}, No. 3 (2006), 759-781. 

\bibitem{MM-Gafa} Y. Martel and F. Merle, {\em Instability of solitons for the critical generalized Korteweg-de Vries equation,} Geom. Funct. Anal. {\bf 11} (2001), 74-123.

\bibitem{MM1} Y. Martel and F. Merle, {\em Asymptotic stability of solitons for subcritical generalized KdV equations,} Arch. Rational Mech. Anal. {\bf 157} (2001), 219-254. 

\bibitem{MM2} Y. Martel and F. Merle, {\em Asymptotic stability of solitons for subcritical gKdV equations revisited,} Nonlinearity {\bf 18} (2005), 55-80.

\bibitem{MM3} Y. Martel and F. Merle, {\em Asymptotic stability of solitons of the gKdV equations with general nonlinearity,} Math. Ann. {\bf 341}, No. 2 (2008), 391--427.

\bibitem{PW} R. Pego and M. Weinstein, {\em Asymptotic stability of solitary waves}, Comm. Math. Physics {\bf 164}, No. 2 (1994), 305-349.

\bibitem{Ste} E. Stein.  {\it Harmonic Analysis: real-variable methods, orthogonality, and oscillatory integrals}, Princeton Mathematical Series, {\bf 43}. Monographs in Harmonic Analysis, III.  Princeton University Press, Princeton, NJ (1993).

\bibitem{Tao} T. Tao, {\em Scattering for the quartic generalised Korteweg-de Vries equation}, J. Diff. Eq. {\bf 3} (2006), 623-651.

\bibitem{Tataru} D. Tataru, {\em Carleman estimates, unique continuation and applications}, Lecture Notes.

\bibitem{Titchmarsh} E.C. Titchmarsh, {\em Eigenfunction Expansions Associated with Second-order Differential Equations}, Oxford at the Clarendon Press, Oxford (1962).

\bibitem{Wat} G. N. Watson, {\em A Treatise on the Theory of Bessel Functions}, Cambridge University Press, Cambridge (1995).

\bibitem{W} M. Weinstein, {\em Modulational stability of ground states of NLS}, SIAM J. Math. Anal. {\bf 16}, No. 3 (1985), 472-491.

\bibitem{Wiener24}
N.~Wiener.
\newblock {The quadratic variation of a function and its Fourier coefficients.}
\newblock In Pesi~Rustom Masani, editor, {\em {Collected works with
  commentaries. Volume II: Generalized harmonic analysis and Tauberian theory;
  classical harmonic and complex analysis}}, volume~15 of {\em Mathematicians
  of Our Time}. {Cambridge, Mass. - London: The MIT Press. XIII, 969 p. }, 1979
  (1924).

\end{thebibliography}
\end{document}